\input amstex\documentstyle{amsppt}  
\pagewidth{12.5cm}\pageheight{19cm}\magnification\magstep1
\topmatter
\title Truncated convolution of character sheaves\endtitle
\author G. Lusztig\endauthor
\address{Department of Mathematics, M.I.T., Cambridge, MA 02139}\endaddress
\thanks{Supported in part by National Science Foundation grant 1303060.}\endthanks   
\endtopmatter   
\document
\define\rk{\text{\rm rk}}

\define\Irr{\text{\rm Irr}}

\define\Bpq{\Bumpeq}

\define\dL{\dot L}

\define\bco{\bar{\co}}

\define\unz{\un{\z}}
\define\unc{\un{\c}}
\define\unb{\un{\bul}}
\define\unst{\un{*}}

\define\mpb{\medpagebreak}

\define\bQ{\bar Q}

\define\da{\dagger}

\define\frl{\forall}

\define\si{\sim}

\define\qua{\quad}

\define\bZ{\bar Z}
\define\lb{\linebreak}

\define\op{\oplus}
   
\define\part{\partial}
\define\emp{\emptyset}
\define\imp{\implies}

\define\n{\notin}
\define\iy{\infty}
\define\m{\mapsto}
\define\do{\dots}

\define\lra{\leftrightarrow}

\define\sub{\subset}    
\define\bxt{\boxtimes}
\define\T{\times}
\define\ti{\tilde}
\define\nl{\newline}
\redefine\i{^{-1}}

\define\un{\underline}
\define\ov{\overline}
\define\ot{\otimes}
\define\bbq{\bar{\QQ}_l}

\define\Hom{\text{\rm Hom}}

\define\tr{\text{\rm tr}}

\define\supp{\text{\rm supp}}

\redefine\spa{\spadesuit}

\define\a{\alpha}
\redefine\b{\beta}
\redefine\c{\chi}
\define\g{\gamma}
\redefine\d{\delta}
\define\e{\epsilon}

\define\io{\iota}
\redefine\o{\omega}
\define\p{\pi}
\define\ph{\phi}
\define\ps{\psi}
\define\r{\rho}
\define\s{\sigma}
\redefine\t{\tau}
\define\th{\theta}

\define\z{\zeta}
\define\x{\xi}

\define\vt{\vartheta}

\redefine\D{\Delta}

\define\Si{\Sigma}

\redefine\L{\Lambda}
\define\Ph{\Phi}

\redefine\aa{\bold a}

\define\boc{\bold c}

\define\kk{\bold k}

\define\pp{\bold p}

\define\ww{\bold w}

\redefine\AA{\bold A}
\define\BB{\bold B}
\define\CC{\bold C}
\define\DD{\bold D}
\define\EE{\bold E}
\define\FF{\bold F}

\define\HH{\bold H}

\define\JJ{\bold J}

\define\LL{\bold L}

\define\NN{\bold N}

\define\QQ{\bold Q}
\define\RR{\bold R}

\define\TT{\bold T}

\define\ZZ{\bold Z}
\define\XX{\bold X}

\define\ca{\Cal A}
\define\cb{\Cal B}
\define\cc{\Cal C}
\define\cd{\Cal D}
\define\ce{\Cal E}

\define\cg{\Cal G}
\define\ch{\Cal H}

\define\ck{\Cal K}

\define\cm{\Cal M}

\define\co{\Cal O}

\define\cu{\Cal U}

\define\cw{\Cal W}
\define\cz{\Cal Z}
\define\cx{\Cal X}
\define\cy{\Cal Y}

\define\fc{\frak c}

\define\fD{\frak D}
\define\fE{\frak E}
\define\fF{\frak F}

\define\fL{\frak L}

\define\fS{\frak S}
\define\fT{\frak T}

\define\fZ{\frak Z}

\define\tih{\ti h}

\define\ty{\ti y}

\define\tC{\ti C}

\define\tL{\ti L}

\define\tP{\ti P}

\define\sh{\sharp}

\define\bp{\bar p}

\define\bul{\bullet}

\define\cir{\bul}

\define\prq{\preceq}

\define\Rep{\text{\rm Rep}}
\define\Reg{\text{\rm Reg}}
\define\tLL{\ti\LL}
\define\bvt{\bar{\vt}}

\define\BBD{BBD}
\define\BN{BN}
\define\BOF{BFO1}
\define\BFO{BFO2}
\define\ENO{ENO}
\define\EW{EW}
\define\GI{Gi}
\define\GR{Gr}
\define\JS{JS}
\define\KL{KL}
\define\LGR{L1}
\define\CSI{L2}
\define\CSII{L3}
\define\CSIII{L4}
\define\COX{L5}
\define\EXO{L6}
\define\CEL{L7}
\define\HEC{L8}
\define\CON{L9}
\define\CDGVII{L10}
\define\CDGIX{L11}
\define\CDGX{L12}
\define\CLE{L13}
\define\MV{MV}
\define\MU{Mu}
\define\MA{Ma}
\define\MG{Mg}
\define\SO{S}

\head Introduction\endhead
\subhead 0.1\endsubhead
Let $\kk$ be an algebraically closed field. Let $G$ be a reductive connected group over 
$\kk$. Let $W$ be the Weyl group of $G$ and let $\boc$ be a two-sided cell of $W$. Let $\cc^\boc G$ the 
category of perverse sheaves on $G$ which are direct sums of unipotent character sheaves whose associated 
two-sided cell (see 1.5) is $\boc$ and let $\cc^\boc\cb^2$ be the category of semisimple 
$G$-equivariant perverse sheaves on $\cb^2$ (the product of two copies of the flag manifold) which belong to
$\boc$. Now $\cc^\boc\cb^2$ has a structure of monoidal category (truncated convolution) introduced in
\cite{\CEL} such that the induced ring structure on the Grothendieck group is the $J$-ring attached to 
$\boc$, see \cite{\HEC, 18.3}. In this paper, we define and study a structure of braided monoidal category 
(truncated convolution) on $\cc^\boc G$ in the case where

(a) $\kk$ is an algebraic closure of a finite field $\FF_q$,
\nl
thus proving a conjecture in \cite{\CON}. In the case where $\kk$ has characteristic zero such a monoidal 
structure was defined by Bezrukavnikov, Finkelberg and Ostrik \cite{\BFO} (in the language of $D$-modules),
who also proved in that case

(b) the existence of an equivalence between $\cc^\boc G$ and the centre of the monoidal category 
$\cc^\boc\cb^2$
\nl
and conjectured that (b) holds without restriction on the characteristic. Note that (b) is made plausible by
the fact that, as a consequence of a conjecture in the last paragraph of \cite{\CEL, 3.2} and of the 
classification of unipotent character sheaves in \cite{\CSIII}, the simple objects of the centre of 
$\cc^\boc\cb^2$ should be in bijection with the simple objects of $\cc^\boc G$. (The idea that the derived 
category of character sheaves with unspecified $\boc$ is equivalent to the centre of the derived category of 
$G$-equivariant sheaves on $\cb^2$ with unspecified $\boc$, appeared in Ben-Zvi and Nadler's paper 
\cite{\BN} and in \cite{\BFO}, again in characteristic zero; we refer to this case as the ``untruncated'' 
case.) 

In this paper we prove (b) in the case where $\kk$ is as in (a), see Theorem 9.5. (In the remainder of this 
paper we assume that $\kk,\FF_q$ are as in (a).) 
Much of the proof involves the definition and study of truncated versions $\unc,\unz,\unst$ of several 
known functors $\c,\z,*$ in the untruncated case. Here $\c$ is the known induction functor from complexes on
$\cb^2$ to complexes on $G$ which I used in the 1980's in the definition of character sheaves; $\z$ is an 
adjoint of $\c$ which I used in the late 1980's to characterize the character sheaves (see 2.5); $*$ is the 
convolution of complexes of sheaves on $G$ defined by Ginzburg \cite{\GI}. The truncated version $\unc$ of 
$\c$ has been already used (but not named) in \cite{\CSIII}. Note that our definition of the truncated 
convolution $\unst$ and truncated restriction $\unz$ involves in an essential way the weight filtrations; it
is not clear how these operations are related to the corresponding operations in characteristic zero 
considered in \cite{\BFO} where weight filtrations do not appear. (In our definition of $\unc$ the 
consideration of weight filtrations is not necessary.) Much of this paper is concerned with establishing
various connections between $\unc,\unz,\unst$. One of these connections, the adjointness of $\unc$ and 
$\unz$ (of which the untruncated version holds by definition) is here surprisingly complicated. We first 
prove a weak form of it (\S8) which we use in the proof of Theorem 9.5 and we then use Theorem 9.5 to 
prove its full form (Theorem 9.8). 

In \S10 we discuss the possibility of a noncrystallographic extension of some of our results, making use of
\cite{\EW}.

Throughout this paper we assume that we have a fixed split $\FF_q$-structure on $G$.

This paper contains several references to results in \cite{\CSIII} which in {\it loc.cit.} are conditional 
on the cleanness of character sheaves; these references are justified since cleanness is now available (see 
\cite{\CLE} and its references). This paper also contains several references to \cite{\HEC,\S14}; these are
justified by the results in \cite{\HEC,\S15}.

We will show elsewhere that the methods and results of this paper extend to
non-unipotent character sheaves on $G$ (at least when the centre of $G$ is connected).

I wish to thank Victor Ostrik for some useful comments.

\subhead 0.2\endsubhead
{\it Notation.} Let $\cb$ be the variety of Borel subgroups of $G$, with the $\FF_q$-structure 
inherited from $G$. Let $\nu=\dim\cb$, $\D=\dim(G)$, $\r=\rk(G)$. We shall view $W$ as an indexing set for 
the orbits of $G$ acting on $\cb^2:=\cb\T\cb$ by simultaneous conjugation; let $\co_w$ be the orbit 
corresponding to $w\in W$ and let $\bco_w$ be the closure of $\co_w$ in $\cb^2$. Note that 
$\co_w,\bco_w$ are naturally defined over $\FF_q$. For $w\in W$ we set $|w|=\dim\co_w-\nu$ (the length
of $w$). Define $w_{max}\in W$ by the condition $|w_{max}|=\nu$.

For $B\in\cb$, let $U_B$ be the unipotent radical of $B$. Then $B/U_B$ is independent of $B$; it is ``the'' 
maximal torus $T$ of $G$. It inherits a split $\FF_q$-structure from $G$. Let $\cx$ be the group of 
characters of $T$.

Let $\Rep W$ be the category of finite dimensional representations of $W$ over $\QQ$; let $\Irr W$ be a set 
of representatives for the isomorphism classes of irreducible objects of $\Rep W$. For any $E\in\Irr W$ we 
denote by $E^\da$ the object of $\Irr W$ which is isomorphic to the tensor product of $E$ and the sign 
representation.

For an algebraic variety $X$ over $\kk$ we denote by $\cd(X)$ the bounded derived category of constructible 
$\bbq$-sheaves on $X$ ($l$ is a fixed prime number invertible in $\kk$); let $\cm(X)$ be the subcategory of 
$\cd(X)$ consisting of perverse sheaves on $X$. If $X$ has a fixed $\FF_q$-structure $X_0$, we denote by 
$\cd_m(X)$ what in \cite{\BBD, 5.1.5} is denoted by $\cd_m^b(X_0,\bbq)$. Note that any object $K\in\cd_m(X)$
can be viewed as an object of $\cd(X)$ which will be denoted again by $K$. For $K\in\cd(X)$ and $i\in\ZZ$ 
let $\ch^iK$ be the $i$-th cohomology sheaf of $K$, $\ch^i_xK$ its stalk at $x\in X$, and let $K^i$ be the 
$i$-th perverse cohomology sheaf of $K$. For $K\in\cd(X)$ (or $K\in\cd_m(X)$) and $n\in\ZZ$ we write 
$K[[n]]=K[n](n/2)$ where $[n]$ is a shift and $(n/2)$ is a Tate twist; we write $\fD(K)$ for the Verdier 
dual of $K$. Let $\cm_m(X)$ be the subcategory of $\cd_m(X)$ whose objects are in $\cm(X)$. If 
$K\in\cm_m(X)$ and $j\in\ZZ$ we denote by $\cw^jK$ the subobject of $K$ which has weight $\le j$ and is such 
that $K/\cw^jK$ has weight $>j$, see \cite{\BBD, 5.3.5}; let $gr_jK=\cw^jK/\cw^{j-1}K$ be the associated pure
perverse sheaf of weight $j$. For $K\in\cd_m(X)$ we shall often write $K^{\{i\}}$ instead of 
$gr_i(K^i)(i/2)$.

If $K\in\cm(X)$ and $A$ is a simple object of $\cm(X)$ we denote by $(A:K)$ the multiplicity of $A$ in a 
Jordan-H\"older series of $K$. 

For $i\in\ZZ$ and $K\in\cd_m(X)$ let $\t_{\le i}K\in\cd_m(X)$ be what in \cite{\BBD} is denoted by 
${}^p\t_{\le i}K$.

Assume that $C\in\cd_m(X)$ and that $\{C_i;i\in I\}$ is a family of objects of $\cd_m(X)$. We shall write
$C\Bpq\{C_i;i\in I\}$ if the following condition is satisfied: there exist distinct elements
$i_1,i_2,\do,i_s$ in $I$, objects $C'_j\in\cd_m(X)$ ($j=0,1,\do,s$) and distinguished triangles 
$(C'_{j-1},C'_j,C_{i_j})$ for $j=1,2,\do,s$ such that $C'_0=0$, $C'_s=C$; moreover, $C_i=0$ unless $i=i_j$ 
for some $j\in[1,s]$. (See \cite{\CDGVII, 32.15}.)

We will denote by $\pp$ the variety consisting of one point. For any variety $X$ let 
$\fL_X=\a_!\bbq\in\cd_mX$ where $\a:X\T T@>>>X$ is the obvious projection. We sometimes write $\fL$ instead 
of $\fL_X$.

Let $v$ be an indeterminate. For any $\ph\in\QQ[v,v\i]$ and any $k\in\ZZ$ we write $(k;\ph)$ for the
coefficient of $v^k$ in $\ph$. Let $\ca=\ZZ[v,v\i]$.

\head Contents\endhead
1. Preliminaries and truncated induction.

2. Truncated restriction.

3. Truncated convolution on $\cb^2$.

4. Truncated convolution on $G$.

5. Truncated convolution and truncated restriction.

6. Analysis of the composition $\unz\unc$.

7. Analysis of the composition $\unz\unc$ (continued).

8. Adjunction formula (weak form).

9. Equivalence of $\cc^\boc G$ with the centre of $\cc^\boc\cb^2$.

10. Remarks on the noncrystallographic case.

\head 1. Preliminaries and truncated induction\endhead
\subhead 1.1\endsubhead
For $y\in W$ let $L_y\in\cd_m(\cb^2)$ be the constructible sheaf which is $\bbq$ (with the standard mixed 
structure of pure weight $0$) on $\co_y$ and is $0$ on $\cb^2-\co_y$; let $L_y^\sh\in\cd_m(\cb^2)$ be its
extension to an intersection cohomology complex of $\bco_y$ (equal to $0$ on $\cb^2-\bco_y$). Let 
$\LL_y=L_y^\sh[[|y|+\nu]]\in\cd_m(\cb^2)$.

Let $r\ge1$. For $\ww=(w_1,w_2,\do,w_r)\in W^r$ we set $|\ww|=|w_1|+\do+|w_r|$.

For any $i<i'$ in $[1,r]$ let $p_{i,i'}:\cb^{r+1}@>>>\cb^2$ be the projection to the $i,i'$ factors. From 
the definitions we see that
$$L^{[1,r]}_\ww:=
p_{01}^*L_{w_1}^\sh\ot p_{12}^*L_{w_2}^\sh\ot\do\ot p_{r-1,r}^*L_{w_r}^\sh\in\cd_m(\cb^{r+1})$$
is the intersection cohomology complex of the projective variety
$$\co^{[1,r]}_\ww=\{(B_0,B_1,\do,B_r)\in\cb^{r+1};(B_{i-1},B_i)\in\bco_{w_i}\frl i\in[1,r]\}$$
extended by $0$ on $\cb^{r+1}-\co^{[1,r]}_\ww$ (it has the standard mixed structure of pure weight $0$). For
any $J\sub[1,r]$ we set 
$$\align&\co^J_\ww
=\{(B_0,B_1,\do,B_r)\in\cb^{r+1};(B_{i-1},B_i)\in\bco_{w_i}\frl i\in J,(B_{i-1},B_i)\in\co_{w_i} \\&
\frl i\in[1,r]-J\}.\endalign$$
Let $i_J:\co^J_\ww@>>>\co^{[1,r]}_\ww$ (resp. $i'_J:\co^{[1,r]}_\ww-\co^J_\ww@>>>\co^{[1,r]}_\ww$) be the 
obvious open (resp. closed) imbedding and let $L_\ww^J\in\cd_m(\cb^{r+1})$ (resp. 
$\dL_\ww^J\in\cd_m(\cb^{r+1})$) be $i_J^*L_\ww^{[1,r]}$ (resp. $i'_J{}^*L_\ww^{[1,r]}$) extended by $0$ on 
$\cb^{r+1}-\co_\ww^J$ (resp.  $\cb^{r+1}-(\co_\ww^{[1,r]}-\co_\ww^J)$); we have a distinguished triangle
$$(L_\ww^J,L_\ww^{[1,r]},\dL_\ww^J)\tag a$$
in $\cd_m(\cb^{r+1})$. We have the following result.

(b) {\it For any $h\in\ZZ$, any composition factor of $(\dL_\ww^J)^h\in\cm(\cb^{r+1})$ is of the form
$L_{\ww'}^{[1,r]}[|\ww'|+\nu]$ for some $\ww'=(w'_1,w'_2,\do,w'_r)\in W^r$ such that $w_i=w'_i$ for all 
$i\in J$.}
\nl
By a standard argument this can be reduced to the case where $r=1$. We then use the fact that
$(\dL_\ww^J)^h\in\cm(\cb^2)$ is equivariant for the diagonal $G$-action and all $G$-equivariant simple
perverse sheaves on $\cb^2$ are of the form $\LL_y$ for some $y\in W$.

We show:

(c) {\it $(\dL_\ww^J[|\ww|+\nu-1])^j=0$ for any $j>0$.}
\nl
It is enough to show that $\dim\supp\ch^h(\dL_\ww^J[|\ww|+\nu-1])\le-h$ for any $h\in\ZZ$. Assume first that
$h\le-|\ww|-\nu$. Since $L_\ww^{[1,r]}$ is an intersection cohomology complex, we have 
$$\dim\supp\ch^{h-1}(L_\ww^{[1,r]}[|\ww|+\nu])<-h+1,$$
hence 
$$\dim\supp\ch^{h-1}(\dL_\ww^J[|\ww|+\nu])<-h+1,$$ 
hence 
$$\dim\supp\ch^{h-1}(\dL_\ww^J[|\ww|+\nu])\le-h,$$ 
hence 
$$\dim\supp\ch^h(\dL_\ww^J[|\ww|+\nu-1])\le-h.$$
Next we assume that $h=-|\ww|-\nu+1$. Then  
$$\dim\supp\ch^{h-1}(\dL_\ww^J[|\ww|+\nu])\le\dim(\co_\ww^{[1,r]}-\co_\ww^J)\le|\ww|+\nu-1=-h,$$
hence $\dim\supp\ch^h(\dL_\ww^J[|\ww|+\nu-1])\le-h$. Now assume that $h\ge-|\ww|-\nu+2$. Then 
$\ch^{h-1}(L_\ww^{[1,r]}[|\ww|+\nu])=0$ hence $\ch^{h-1}(\dL_\ww^J[|\ww|+\nu])=0$ hence 
$\ch^h(\dL_\ww^J[|\ww|+\nu-1])=0$. This proves (c). 

\subhead 1.2\endsubhead
For ${}^1L,{}^2L,\do,{}^rL$ in $\cd_m(\cb^2)$ we set
$${}^1L\cir{}^2L\cir\do\cir{}^rL=p_{0r!}(p_{01}^*{}^1L\ot p_{12}^*{}^2L\ot\do\ot p_{r-1,r}^*{}^rL)
\in\cd_m(\cb^2).$$ 
If $\ww=(w_1,w_2,\do,w_r)$ is as in 1.1 we set
$$L_\ww^\cir=p_{0r!}L^{[1,r]}_\ww=L_{w_1}^\sh\cir L_{w_2}^\sh\cir\do\cir L_{w_r}^\sh\in\cd_m(\cb^2).$$
If $J$ is as in 1.1, then
$$p_{0r!}L^J_\ww={}^1L\cir {}^2L\cir\do\cir{}^rL\in\cd_m(\cb^2)\tag a$$
where ${}^iL=L_{w_i}^\sh$ for $i\in J$, ${}^iL=L_{w_i}$ for $i\in [1,r]-J$.

Using the decomposition theorem \cite{\BBD} for the proper map $p_{0r}$, we see that
$$L_\ww^\cir[|\ww|]\cong\op_{w\in W,k\in\ZZ}(L_w^\sh[k+|w|])^{\op N(w,k)}\tag b$$
in $\cd(\cb^2)$ where $N(w,k)\in\NN$.

\subhead 1.3\endsubhead
Let $\HH$ be the Hecke algebra of $W$ (see \cite{\HEC, 3.2} with $L(w)=|w|$) over $\ca$ 
and let $\{c_w;w\in W\}$ be the ``new'' basis of $\HH$, see \cite{\HEC, 5.2}. As in 
\cite{\HEC, 13.1}, for $x,y\in W$ we write $c_xc_y=\sum_{z\in W}h_{x,y,z}c_z$ where 
$h_{x,y,z}\in\ca$. For $x,z\in W$ we write $z\prq x$ if there exists $\x\in\HH c_x\HH$ such that 
$c_z$ appears with $\ne0$ coefficent in the expansion of $\x$ in the new basis. This is a preorder on $W$.
Recall that the two-sided cells of $W$ are the equivalence classes associated to this preorder. For 
$x,y\in W$ we write $x\si y$ if $x,y$ belong to the same two-sided cell, that is $x\prq y$ and $y\prq x$. 
For $x,y\in W$ we write $x\si_Ly$ if $x,y$ belong to the same left cell of $W$, see \cite{\HEC, 8.1}. If 
$\boc$ is a two-sided cell and $w\in W$ we write $w\prq\boc$ (resp. $\boc\prq w$) if $w\prq w'$ (resp. 
$w'\prq w$) for some $w'\in\boc$; we write $w\prec\boc$ (resp. $\boc\prec w$) if $w\prq\boc$ (resp. 
$\boc\prq w$) and $w\n\boc$. If $\boc,\boc'$ are two-sided cells we write $\boc\prq\boc'$ (resp. 
$\boc\prec\boc'$) if $w\prq w'$ (resp. $w\prec w'$) for some $w\in\boc,w'\in\boc'$. Let $\aa:W@>>>\NN$ be 
the $\aa$-function in \cite{\HEC, 13.6}. 

If $\boc$ is a two-sided cell, then for all $w\in\boc$ we have $\aa(w)=\aa(\boc)$ where $\aa(\boc)$ is a 
constant. Note that the numbers $N(w,k)$ in 1.2(b) satisfy:
$$c_{w_1}c_{w_2}\do c_{w_r}=\sum_{w\in W}\ph_wc_w\text{ where }\ph_w=\sum_{k\in\ZZ}N(w,k)v^k.\tag a$$
If $x,y,z\in W$ then 
$$h_{x,y,z}=h_{x,y,z}^*v^{-\aa(z)}+\text{higher powers of $v$},$$
$$h_{x,y,z}=h_{x,y,z}^*v^{\aa(z)}+\text{lower powers of $v$}$$
where $h_{x,y,z}^*\in\NN$; moreover, if $h_{x,y,z}\ne0$ then $\aa(x)\le\aa(z)$, $\aa(y)\le\aa(z)$ (see 
\cite{\HEC, P4}); if $h_{x,y,z}^*\ne0$ then $x\si y\si z$ (see \cite{\HEC, P8}) hence 
$\aa(x)=\aa(y)=\aa(z)$. 

If $\boc$ is a two-sided cell of $W$ then the subquotient
$$(\op_{w\in W;w\prq\boc}\QQ c_w)/(\op_{w\in W;w\prec\boc}\QQ c_w)$$ 
of the group algebra $\QQ[W]$ is naturally an object $[\boc]$ of $\Rep W$. 
If $\L$ is a left cell of $W$ then the subquotient
$$(\op_{w\in W;w\in\L\text{ or }w\prec\boc}\QQ c_w)/(\op_{w\in W;w\prec\boc}\QQ c_w)$$ 
of the group algebra $\QQ[W]$ is naturally an object $[\L]$ of $\Rep W$. 

For $E\in\Irr W$, there is a unique two-sided cell $\boc_E$ of $W$ such that $[\boc_E]$  contains $E$. 
(This differs from the usual definition of two -sided cell attached to $E$ by multiplication on the left or 
right by $w_{max}$.)

{\it Until the end of \S9 we fix a two-sided cell $\boc$ of $W$ and we set $a=\aa(\boc)$.}
\nl
Since for $w\in\boc$ we have (in (a)):
$$\ph_w=\sum_{z_2,z_3,\do,z_{r-1}\text{ in }W}h_{w_1,w_2,z_2}h_{z_2,w_3,z_3}\do h_{z_{r-1},w_r,w},$$
we see that 
$$\align&N(w,k)\ne0\imp k\ge-(r-1)a;N(w,-(r-1)a)\ne0\imp w_i\in\boc\text{ for all }i.\tag b\endalign$$ 
In addition, if $w_1,w_2,\do,w_r$ are in $\boc$, then 
$$N(w,-(r-1)a)=\sum h_{w_1,w_2,z_2}^*h_{z_2,w_3,z_3}^*\do h_{z_{r-1},w_r,w}^*\tag c$$
where the sum is taken over all $z_2,z_3,\do,z_{r-1}$ in $\boc$.

\mpb

Let $\JJ$ be the free $\ZZ$-module with basis $\{t_z;z\in W\}$. It is known (see \cite{\HEC, 18.3}) that 
there is a well defined structure of associative ring (with $1$) on $\JJ$ such that if $x,y\in W$ then 
$t_xt_y=\sum_{z\in W}h_{x,y,z}^*t_z$. For each two-sided cell $\boc'$ let $\JJ^{\boc'}$ be the subgroup of 
$\JJ$ generated by $\{t_z;z\in\boc'\}$. Then $\JJ^{\boc'}$ is a subring of $\JJ$ and we have
$\JJ=\op_{\boc'}\JJ^{\boc'}$ (as rings) where $\boc'$ runs over the two-sided cells of $W$.

If $w_1,w_2,\do,w_r$ above are in $\boc$ then clearly,
$$t_{w_1}t_{w_2}\do t_{w_r}=\sum_{w\in\boc}N(w,-(r-1)a)t_w\tag d$$
where $N(w,-(r-1)a)$ is as in (c).

The unit element of $\JJ^{\boc'}$ is $\sum_{d\in\DD_{\boc'}}t_d$ where 
$\DD_{\boc'}$ is the set of distinguished involutions of $\boc'$. Let $\DD=\cup_{\boc'}\DD_{\boc'}$.
We define $\ps:\HH@>>>\ca\ot\JJ$ by $\ps(c_w)=\sum_{z\in W,d\in\DD;\aa(d)=\aa(z)}h_{x,d,z}t_z$.
From \cite{\HEC, 18.8} we see that $\ps$ is a homomorphism of $\ca$-algebras with $1$. Specializing $v$ to 
$1$ we get a ring homomorphism $\ps_1:\ZZ[W]@>>>\JJ$ where $\ZZ[W]$ is the group algebra of $W$. This becomes
an isomorphism $\ps_1^\QQ$ after tensoring by $\QQ$ (see \cite{\HEC, 20.1}). For $E\in\Irr W$ we denote by 
$E_\iy$ the simple $\QQ\ot\JJ$-module which corresponds to $E$ under $\ps_1^\QQ$. Now the 
$\QQ(v)\ot\JJ$-module $\QQ(v)\ot_\QQ E_\iy$ can be regarded as a $\QQ(v)\ot_\ca\HH$-module $E(v)$ via the 
algebra homomorphism (actually an isomorphism) $\QQ(v)\ot_\ca\HH@>>>\QQ(v)\ot\JJ$ induced by $\ps$.

Let $\Irr_\boc W=\{E\in\Irr W;\boc_E=\boc\}$. Let $E\in\Irr W$. From the definitions we see that we have 
$E\in\Irr_\boc W$ if and only if $E_\iy$ is a simple $\QQ\ot\JJ^\boc$-module. From the definitions, for any 
$z\in\boc$ we have
$$\tr(c_z,E(v))=\tr(t_z,E_\iy)v^a+\text{ lower powers of }v.\tag e$$   

\proclaim{Lemma 1.4} Let $r\ge1$ and let $\ww=(w_1,w_2,\do,w_r)\in W^r$.

(a) Assume that $w_i\in\boc$ for some $i\in[1,r]$ and that $w\in W,k\in\ZZ$ are such that $N(w,k)$ in 1.2(b)
is $\ne0$. Then either $w\in\boc$, $k\ge-(r-1)a$ or $w\prec\boc$; if $w\in\boc$ and $k=-(r-1)a$, then 
$w_j\in\boc$ for all $j\in[1,r]$. 

(b) Assume that $w_i\in\boc$ for some $i\in[1,r]$. If $j\in\ZZ$ (resp. $j>\nu+(r-1)a$) then 
$(L_\ww^\cir[|\ww|])^j$ is a direct sum of simple perverse sheaves of the form $\LL_z$ where $z\in W$ 
satisfies $z\prq\boc$ (resp. $z\prec\boc$).

(c) Assume that $w_i\prec\boc$ for some $i\in[1,r]$ and that $w\in W$, $k\in\ZZ$ are such that $N(w,k)$ in 
1.2(b) is $\ne0$. Then $w\prec\boc$.

(d) Assume that $w_i\prec\boc$ for some $i\in[1,r]$. If $j\in\ZZ$ then $(L_\ww^\cir[|\ww|])^j$ is a direct 
sum of simple perverse sheaves of the form $\LL_z$ where $z\in W$ satisfies $z\prec\boc$.
\endproclaim
We prove (a). If $r=1$ the result is obvious (we have $k=0$). We now assume that $r\ge2$. From the 
definitions we see that there exists a permutation $w'_1,w'_2,\do,w'_r$ of $w_1,w_2,\do,w_r$, a sequence 
$z_1,z_2,\do,z_r$ in $W$ and a sequence $f_2,\do,f_r$ in $\NN[v,v\i]$ such that (i) $z_1=w'_1\in\boc$, 
$z_r=w$, (ii) for any $i\in[2,r]$, $c_{z_i}$ appears with coefficient $f_i$ in $c_{z_{i-1}}c_{w'_i}$ or in 
$c_{w'_i}c_{z_{i-1}}$, (iii) $(k;f_2f_3\do f_n)\ne0$ (see 0.2). From the definition of $\prq$ we have 
$z_r\prq z_{r-1}\prq\do\prq z_2\prq z_1$. Hence $z_r\prq\boc$, 
$\aa(z_r)\ge\aa(z_{r-1})\ge\do\ge\aa(z_2)\ge\aa(z_0)=a$ (see \cite{\HEC, P4}), 
$v^{\aa(z_i)}f_i\in\NN[v]$ for $i=2,\do,r$. Hence $v^{\aa(z_2)+\do+\aa(z_r)}f_2f_3\do f_r\in\NN[v]$ so that 
$k+\aa(z_2)+\do+\aa(z_r)\ge0$. We also see that if $z_r\in\boc$ so that $\aa(z_r)=a$ then
$\aa(z_r)=\aa(z_{r-1})=\do=\aa(z_2)=a$ and $k+(r-1)a\ge0$. Now assume that $z_r\in\boc$ and $k=-(r-1)a$. 
Then $v^{\aa(z_2)+\do+\aa(z_r)}f_2f_3\do f_r\in\NN[v]$ has $\ne0$ constant term hence 
$v^{\aa(z_i)}f_i\in\NN[v]$ has $\ne0$ constant term for $i=2,\do,r$. Using \cite{\HEC, P8}, we deduce that 
$z_i,z_{i-1},w'_i$ are in the same two-sided cell for $i=2,\do,r$. Hence $w'_2,\do,w'_r$ are in $\boc$. 
Since $w'_1\in\boc$, we see that $w_1,\do,w_r$ are in $\boc$. This proves (a).

Note that in (a) we have necessarily $w\prq\boc$. Replacing $\boc$ in (a) by the two-sided cell containing
$w_i$ in (c) we deduce that (c) holds. 

We prove (b). By 1.2(b) we have 
$$(L_\ww^\cir[|\ww|])^j\cong\op_{w\in W,k\in\ZZ}((\LL_w)^{j+k-\nu})^{\op N(w,k)}
=\op_{w\in W}(\LL_w)^{\op N(w,\nu-j)}.\tag e$$
Hence if $\LL_z$ appears as a summand in the last direct sum then $N(z,\nu-j)\ne0$. Using (a) we see that 
$z\prq\boc$ and that $z\prec\boc$ if $\nu-j<-(r-1)a$. This proves (b). The same proof, using (c) instead of 
(a) yields (d).

\subhead 1.5\endsubhead
We consider the maps $\cb^2@<f<<X@>\p>>G$ where 
$$X=\{(B,B',g)\in\cb\T\cb\T G;gBg\i=B'\},f(B,B',g)=(B,B'),\p(B,B',g)=g.$$ 
Now $L\m\c(L)=\p_!f^*L$ defines a functor $\cd_m(\cb^2)@>>>\cd_m(G)$.
For $i\in\ZZ,L\in\cd_m(\cb^2)$ we write $\c^i(L)$ instead of $(\c(L))^i$.

The functor $\c$ is the main tool used in the definition \cite{\CSI} of (unipotent) character sheaves. For 
any $z\in W$ we set $R_z=\c(L_z^\sh)\in\cd_m(G)$. A {\it unipotent character sheaf} is a simple perverse
sheaf $A\in\cm(G)$ such that $(A:R_z^j)\ne0$ for some $z\in W,j\in\ZZ$. Let $CS(G)$ be a set of 
representatives for the isomorphism classes of unipotent character sheaves. 

By \cite{\CSIII, 14.11}, for any $A\in CS(G)$, any $z\in W$ and any $j\in\ZZ$ we have
$$(A:R_z^j)=(j-\D-|z|;(-1)^{j-\D}\sum_{E\in\Irr W}c_{A,E}\tr(c_z,E(v)))\tag a$$   
where $E(v)$ is as in 1.3 and $c_{A,E}$ are uniquely defined rational numbers; for $E'\in Rep(W)$ we set
$$c_{A,E'}=\sum_{E\in\Irr W}(\text{multiplicity of $E$ in $E'$})c_{A,E}.$$
Moreover, given $A$, there is a unique 
two-sided cell $\boc_A$ of $W$ such that $c_{A,E}=0$ whenever $E\in\Irr W$ satisfies $\boc_E\ne\boc_A$; 
this differs from the two-sided cell associated to $A$ in \cite{\CSIII, 16.7} by multiplication on the left 
or right by $w_{max}$. Note that 

(b) {\it $(A:R_z^j)\ne0$ for some $z\in\boc_A,j\in\ZZ$ and conversely, if $(A:R_z^j)\ne0$ for 
$z\in W,j\in\ZZ$, then $\boc_A\prq z$;}
\nl
see \cite{\CDGIX, 41.8}, \cite{\CDGX, 44.18}. 

For example, if $G=GL_2(\kk)$ and $W=\{1,s\}$, we have $CS(G)=\{A_0,A_1\}$ with $A_1\not\cong A_0=\bbq[\D]$,
and $R_1=A_0[-\D]\op A_1[-\D]$, $R_s=A_0[-\D]\op A_0[-\D-2]$. Thus $R_1^\D\cong A_0\op A_1$, $R_1^j=0$ if 
$j\ne\D$ and $R_s^\D\cong A_0$, $R_s^{\D+2}\cong A_0$, $R_s^j=0$ if $j\n\{\D,\D+2\}$. We have 
$\Irr W=\{E_0,E_1\}$ where $E_0$ is the unit representation, $E_1$ is the sign representation and
$$\tr(c_1,E_0(v))=1,\tr(c_1,E_1(v))=1,\tr(c_s,E_0(v))=v+v\i,\tr(c_s,E_1(v))=0.$$ 
It follows that $c_{A_i,E_j})=\d_{ij}$ for $i,j\in\{0,1\}$. Hence $\boc_{A_0}=\boc_{E_0}=\{s\}$ (resp. 
$\boc_{A_1}=\boc_{E_1}=\{1\}$).

We return to the general case. For $A\in CS(G)$ let $a_A$ be the value of the $\aa$-function on $\boc_A$. If
$z\in W,E\in\Irr(W)$ satisfy $\tr(c_z,E(v))\ne0$ then $\boc_E\prq z$; if in addition we have $z\in\boc_E$, 
then 
$$\tr(c_z,E(v))=\g_{z,E}v^{a_E}+\text{lower powers of }v$$ 
where $\g_{z,E}\in\ZZ$ and $a_E$ is the value of the $\aa$-function on $\boc_E$. Hence from (a) we see that

(c) {\it $(A:R_z^j)=0$ unless $\boc_A\prq z$ and, if $z\in\boc_A$, then}
$$\align&(A:R_z^j)\\&=(-1)^{j+\D}(j-\D-|z|;(\sum_{E\in\Irr W;\boc_E=\boc_A}c_{A,E}\g_{z,E})v^{a_A}
+\text{ lower powers of }v))\endalign$$
{\it which is $0$ unless $j-\D-|z|\le a_A$.}
\nl
For $Y=G$ or $Y=\cb^2$ let $\cm^\spa Y$ be the category of perverse sheaves on $Y$ whose composition factors 
are all of the form $A\in CS(G)$, when $Y=G$, or of the form $\LL_z$ with $z\in W$ (when $Y=\cb^2$). 
Let $\cm^\prq Y$ (resp. $\cm^\prec Y$) be the category of perverse sheaves on $Y$ 
whose composition factors are all of the form $A\in CS(G)$ with $\boc_A\prq\boc$ (resp. $\boc_A\prec\boc$),
when $Y=G$, or of the form $\LL_z$ with $z\prq\boc$ (resp. $z\prec\boc$) when $Y=\cb^2$. 
Let $\cd^\spa Y$ (resp. $\cd^\prq Y$ or $\cd^\prec Y$) be the 
category of all $K\in\cd(Y)$ such that $K^i\in\cm^\spa Y$ (resp. $K^i\in\cm^\prec Y$ or
$K^i\in\cm^\prec Y$) for all $i\in\ZZ$. Let $\cm_m^\spa Y$ (or $\cm^\prq_mY$, or $\cm^\prec_mY$) be the
category of all $K\in\cm_mY$ which are also in $\cm^\spa Y$ (or $\cm^\prq Y$ or $\cm^\prec Y$). Let 
$\cd_m^\spa Y$ (or $\cd^\prq_mY$, or $\cd^\prec_mY$) be the category of all $K\in\cd_mY$ which are also 
in $\cd^\spa Y$ (or $\cd^\prq Y$ or $\cd^\prec Y$).
From (c) we deduce:

(d) {\it If $z\prq\boc$ then $R_z^j\in\cm^\prq G$ for all $j\in\ZZ$ and. If $z\in\boc$ and $j>a+\D+|z|$, then
$R_z^j\in\cm^\prec G$. If $z\prec\boc$ then $R_z^j\in\cm^\prec G$ for all $j\in\ZZ$.}

\proclaim{Lemma 1.6}Let $r\ge1$, $J\sub[1,r]$, $J\ne\emp$ and $\ww=(w_1,w_2,\do,w_r)\in W^r$. Let 
$\fE=\D+ra$.

(a) Assume that $w_i\in\boc$ for some $i\in[1,r]$. If $j\in\ZZ$ (resp. $j>\fE$) then 
$\c^j(p_{0r!}L_\ww^{[1,r]}[|\ww|])$ is in $\cm^\prq G$ (resp. $\cm^\prec G$).   

(b) Assume that $w_i\in\boc$ for some $i\in J$. If $j\in\ZZ$ (resp. $j\ge\fE$) then 
$\c^j(p_{0r!}\dL_\ww^J[|\ww|])$ is in $\cm^\prq G$ (resp. $\cm^\prec G$).

(c) Assume that $w_i\in\boc$ for some $i\in J$. If $j\ge\fE$ then the cokernel of the map 
$$\c^j(p_{0r!}L_\ww^J[|\ww|])@>>>\c^j(p_{0r!}L_\ww^{[1,r]}[|\ww|])$$
associated to 1.1(a) is in $\cm^\prec G$.

(d) Assume that $w_i\in\boc$ for some $i\in J$. If $j\in\ZZ$ (resp. $j>\fE$) then 
$\c^j(p_{0r!}L_\ww^J[|\ww|])$ is in $\cm^\prq G$ (resp. $\cm^\prec G$).

(e) Assume that $w_i\prec\boc$ for some $i\in J$. If $j\in\ZZ$ then 
$\c^j(p_{0r!}L_\ww^{[1,r]}[|\ww|])\in\cm^\prec G$ and $\c^j(p_{0r!}L_\ww^J[|\ww|])\in\cm^\prec G$.
\endproclaim  
We prove (a). Let $A$ be a simple perverse sheaf on $G$ and let $j\in\ZZ$ be such that $A$ is a composition
factor of $\c^j(p_{0r!}L_\ww^{[1,r]}[|\ww|])=\c^{j+|\ww|}(L_\ww^\cir)$. Then there exists $h'$ such that 
$(A:\c^{j+|\ww|-h'}((L_\ww^\cir)^{h'}))\ne0$. By 1.2(b) we have 
$$\align&(L_\ww^\cir)^{h'}\cong\op_{w\in W,k\in\ZZ}((\LL_w[k-|\ww|-\nu])^{h'})^{\op N(w,k)}\\&=
\op_{w\in W,k\in\ZZ}((\LL_w)^{h'+k-|\ww|-\nu})^{\op N(w,k)}=\op_{w\in W}(\LL_w)^{\op N(w,|\ww|+\nu-h')}.
\endalign$$
Hence $A$ is a composition factor of 
$$\op_{w\in W}(\c^{j+|\ww|-h'}(\LL_w))^{\op N(w,|\ww|+\nu-h')}.$$
Thus there exists $z\in W$ such that $N(z,|\ww|+\nu-h')\ne0$ and $(A:\c^{j+|\ww|-h'}(\LL_z))\ne0$. From 
$N(z,|\ww|+\nu-h')\ne0$ and 1.4(a) we see that $z\prq\boc$. We also see that $A\in CS(G)$ and $\boc_A\prq z$,
see 1.5(b); hence $\boc_A\prq\boc$. If $z\prec\boc$ or if $\boc_A\prec z$ then $\boc_A\prec\boc$. Assume now
that $z\in\boc$ and $z\in\boc_A$ so that $\boc_A=\boc$. From 1.4 we see that $|\ww|+\nu-h'\ge-(r-1)a$ that 
is $h'\le|\ww|+\nu+(r-1)a$.

We have $(A:R_z^{j-h'+|\ww|+\nu+|z|})\ne0$ hence by 1.5(c) we have 
$$j-h'+|\ww|+\nu+|z|-\D-|z|\le a_A$$ 
that is 
$j-h'+|\ww|+\nu-\D\le a$. Combining this with the inequality $h'\le|\ww|+\nu+(r-1)a$ we obtain $j\le\D+ra$. 
This proves (a).

We prove (b). Let $A$ be a simple perverse sheaf on $G$ and $j\in\ZZ$ be such that
$(A:\c^j(p_{0r!}\dL_\ww^J[|\ww|]))\ne0$. There exists $h$ such that
$(A:\c^j(p_{0r!}(\dL_\ww^J[|\ww|])^h)[-h]))\ne0$. We have $(\dL_\ww^J[|\ww|])^h\ne0$ hence 
$(\dL_{\ww}^J[|\ww|+\nu-1])^{h-\nu+1}\ne0$ hence by 1.1(c), $h-\nu+1\le0$. From 1.1(b) we see that there 
exists $\ww'=(w'_1,w'_2,\do,w'_r)\in W^r$ such that $w_i=w'_i$ for all $i\in J$ and such that $A$ is a 
composition factor of 
$$\c^j(p_{0r!}(L_{\ww'}^{[1,r]}[|\ww'|+\nu])[-h])=\c^{j+\nu-h}(p_{0r!}(L_{\ww'}^{[1,r]}[|\ww'|]))
=\c^{j+\nu-h}(L_{\ww'}^\cir[|\ww'|]).$$
From (a) (for $\ww'$ instead of $\ww$) we see that $A\in CS(G)$, $\boc_A\prq\boc$ and that 
$\boc_A\prec\boc$ if $j+\nu-h>\D+ra$ that is, if $j>h+\D-\nu+ra$. If $j\ge\D+ra$ then using $h-\nu+1\le0$ 
(that is $0>h-\nu$) we see that we have indeed $j>h+\D-\nu+ra$. This proves (b).

We prove (c). From 1.1(a) we get a distinguished triangle
$$(\c(p_{0r!}L_\ww^J[[|\ww|]]),\c(p_{0r!}L_\ww^{[1,r]}[[|\ww|]]),\c(p_{0r!}\dL_\ww^J[[|\ww|]]))$$
in $\cd_m(G)$. This gives rise for any $j$ to an exact sequence
$$\align&\c^{j-1}(p_{0r!}\dL_\ww^J[[|\ww|]])@>>>\c^j(p_{0r!}L_\ww^J[[|\ww|]])
@>>>\c^j(p_{0r!}L_\ww^{[1,r]}[[|\ww|]])\\&@>>>\c^j(p_{0r!}\dL_\ww^J[[|\ww|]])\tag f\endalign$$
in $\cm_m(G)$. Using this and (b) we see that (c) holds.

Now (d) follows from the previous exact sequence using (a),(b).

Replacing $\boc$ in (a) and (d) by the two-sided cell containing $w_i$ in (e) we deduce that (e) holds. 

\subhead 1.7\endsubhead
Let $CS_\boc=\{A\in CS(G);\boc_A=\boc\}$. For any $z\in\boc$ we set $n_z=a+\D+|z|$. Let $A\in CS_\boc$ and 
let $z\in\boc$. We have:
$$(A:R_z^{n_z})=(-1)^{a+|z|}\sum_{E\in\Irr_\boc W}c_{A,E}\tr(t_z,E_\iy).\tag a$$
Indeed, from 1.5(a) we have
$$(A:R_z^{n_z})=(-1)^{a+|z|}\sum_{E\in\Irr_\boc W}c_{A,E}(a;\tr(c_z,E(v)))$$
and it remains to use 1.3(e). We show:

(b) {\it For any $A\in CS_\boc$ there exists $z\in\boc$ such that $(A:R_z^{n_z})\ne0$.}
\nl
Assume that this is not so. Then, using (a), we see that 
$$\sum_{E\in\Irr_\boc W}c_{A,E}\tr(t_z,E_\iy)=0$$
for any $z\in\boc$. This shows that the linear functions $t_z\m\tr(t_z,E_\iy)$ on $\JJ^\boc$ (for various 
$E$ as above) are linearly dependent. (It is known that $c_{A,E}\ne0$ for some $E\in\Irr_\boc W$.) This 
is a contradiction since the $E_\iy$ form a complete set of simple modules for the semisimple algebra 
$\QQ\ot\JJ^\boc$.

\mpb

Let $\boc^0=\{z\in\boc;z\si_Lz\i\}$. If $z\in\boc-\boc^0$ and $E\in\Irr_\boc W$, then $\tr(t_z,E_\iy)=0$ 
(see \cite{\HEC, 24.2}). From this and (a) we deduce

(c) {\it If $z\in\boc-\boc^0$, then $R_z^{n_z}=0$.}

\subhead 1.8\endsubhead 
For $Y=G$ or $\cb^2$ let $\cc^\spa Y$ be the subcategory of $\cm^\spa Y$ consisting of semisimple objects; 
let $\cc^\spa_0Y$ be the subcategory of $\cm_mY$ consisting of those $K\in\cm_mY$ such that $K$ is pure of 
weight $0$ and such that as an object of $\cm(Y)$, $K$ belongs to $\cc^\spa Y$. Let $\cc^\boc Y$ be the 
subcategory of $\cm^\spa Y$ consisting of objects which are direct sums of objects in $CS_\boc$ (if $Y=G$) 
or of the form $\LL_z$ with $z\in\boc$ (if $Y=\cb^2$). Let $\cc^\boc_0Y$ be the subcategory of $\cc^\spa_0Y$
consisting of those $K\in\cc^\spa_0Y$ such that as an object of $\cc^\spa Y$, $K$ belongs to $\cc^\boc Y$. 
For $K\in\cc^\spa_0Y$, let $\un{K}$ be the largest subobject of $K$ such that as an object of $\cc^\spa Y$, 
we have $\un{K}\in\cc^\boc Y$.   

\proclaim{Proposition 1.9} (a) If $L\in\cd^\prq\cb^2$ then $\c(L)\in\cd^\prq G$. If $L\in\cd^\prec\cb^2$ 
then $\c(L)\in\cd^\prec G$.

(b) If $L\in\cm^\prq\cb^2$ and $j>a+\nu+\r$ then $\c^j(L)\in\cm^\prec G$.
\endproclaim
It is enough to prove the proposition assuming in addition that $L=\LL_z$ where $z\prq\boc$. Then (a) 
follows from 1.6(a),(e). Im the setup of (b) we have $\c^j(\LL_z)=\c^{j+\nu}(L_z^\sh[[|z|]])$ and this is in
$\cm^\prec G$ since $j+\nu>\D+a$, see 1.6(a).

\subhead 1.10\endsubhead
For $L\in\cc^\boc_0\cb^2$ we set 
$$\unc(L)=\un{(\c^{a+\nu+\r}(L)}((a+\nu+\r)/2)=\un{(\c(L))^{\{a+\nu+\r\}}}\in\cc^\boc_0G.$$
The functor $\unc:\cc^\boc_0\cb^2@>>>\cc^\boc_0G$ is called {\it truncated induction}. For $z\in\boc$ we 
have 
$$\unc(\LL_z)=\un{R_z^{n_z}}(n_z/2).\tag a$$
Indeed,
$$\align&\unc(\LL_z)=\un{\c^{a+\nu+\r}(\LL_z)}((a+\nu+\r)/2)
=\un{(\c(L_z^\sh[[|z|+\nu]])^{a+\nu+\r}}((a+\nu+r)/2)\\&
=\un{\c^{a+\D+|z|}(L_z^\sh)}((a+\D+|z|)/2)=\un{\c^{n_z/2}(L_z^\sh)}(n_z/2).\endalign$$
We shall denote by $\t:\JJ^\boc@>>>\ZZ$ the group homomorphism such that $\t(t_z)=1$ if $z\in\DD_\boc$ and 
$\t(t_z)=0$, otherwise. For $z,u\in\boc$ we show:
$$\dim\Hom_{\cc^\boc G}(\unc(\LL_z),\unc(\LL_u))=\sum_{y\in\boc}\t(t_{y\i}t_zt_yt_{u\i}).\tag b$$
Using (a) and the definitions we see that the left hand side of (b) equals
$$\sum_{A\in CS_\boc}(A:R_z^{n_z})(A:R_u^{n_u}),$$
hence, using 1.7(a) it equals
$$\sum_{E,E'\in\Irr_\boc W}
(-1)^{|z|+|u|}\sum_{A\in CS_\boc}c_{A,E}c_{A,E'}\tr(t_z,E_\iy)\tr(t_u,E'_\iy).$$
Replacing in the last sum $\sum_{A\in CS_\boc}c_{A,E}c_{A,E'}$ by $1$ if $E=E'$ and by $0$ if 
$E\ne E'$ (see \cite{\CSIII, 13.12}) we obtain 
$$\sum_{E\in\Irr_\boc W}(-1)^{|z|+|u|}\tr(t_z,E_\iy)\tr(t_u,E_\iy).$$
This is equal to $(-1)^{|z|+|u|}$ times the trace of the operator $\x\m t_z\x t_{u\i}$ on $\JJ^\boc\ot\CC$.
The last trace is equal to the sum over $y\in\boc$ of the coefficient of $t_y$ in $t_zt_yt_{u\i}$; this 
coefficient is equal to $\t(t_{y\i}t_zt_yt_{u\i})$ since for $y,y'\in\boc$, $\t(t_{y'}t_y)$ is $1$ if 
$y'=y\i$ and is $0$ if $y'\ne y\i$ (see \cite{\HEC, 20.1(b)}). Thus we have
$$\dim\Hom_{\cc^\boc G}(\unc(\LL_z),\unc(\LL_u))=(-1)^{|u|+|z|}\sum_{y\in\boc}\t(t_{y\i}t_zt_yt_{u\i}).$$
Since $\dim\Hom_{\cc^\boc G}(\unc(\LL_z),\unc(\LL_u))\in\NN$ and 
$\sum_{y\in\boc}\t(t_{y\i}t_zt_yt_{u\i})\in\NN$, it follows that (b) holds.

\subhead 1.11\endsubhead
A version of the following result (at the level of Grothendieck groups) appears in \cite{\CSI}.

(a) {\it Let $L,L'\in\cd_m(\cb^2)$. Assume that $L'$ is a $G$-equivariant perverse sheaf. We have 
canonically $\c(L\cir L')=\c(L'\cir L)$.}
\nl
Let $Z=\cb^2\T G$. Define $c:Z@>>>\cb^2\T\cb^2\T G$ by 
$$c((B_1,B_2),g)=((B_1,B_2),(B_2,gB_1g\i),g)$$ 
and $d:Z@>>>G$ by $d((B_1,B_2),g)=g$. Define $c':Z@>>>\cb^2\T\cb^2\T G$ by 
$$c'((B_1,B_2),g)=((B_2,gB_1g\i),(B_1,B_2),g).$$
We have 
$$\c(L\cir L')=d_!c^*(L\bxt L'\bxt\bbq),\qua \c(L'\cir L)=d_!c'{}^*(L\bxt L'\bxt\bbq).$$
Define $t:Z@>>>Z$, $u:\cb^2\T\cb^2\T G$ by 
$$t((B_1,B_2),g)=((B_2,gB_1g\i),g),$$
$$u((B_1,B_2),(B_3,B_4),g)=((B_1,B_2),(gB_3g\i,gB_4g\i),g).$$
We have $ct=uc'$, $dt=d$. Since $L'$ is $G$-equivariant we have canonically \lb
$u^*(L\bxt L'\bxt\bbq)=L\bxt L'\bxt\bbq$. Hence
$$\align&d_!c^*(L\bxt L'\bxt\bbq)=d_!t_!t^*c^*(L\bxt L'\bxt\bbq)=d_!c'{}^*u^*(L\bxt L'\bxt\bbq)\\&=
d_!c'{}^*(L\bxt L'\bxt\bbq).\endalign$$
This proves (a).

\mpb

We will not use (a) in this paper; a characteristic zero analogue of (a) plays a role in \cite{\BFO}.

\proclaim{Lemma 1.12} Let $Y_1,Y_2$ be among $G,\cb^2$ and let $\XX\in\cd_m^\prq Y_1$. Let $c,c'$ be 
integers and let $\Ph:\cd_m^\prq Y_1@>>>\cd_m^\prq Y_2$ be a functor which takes distinguished triangles to 
distinguished triangles, commutes with shifts, maps $\cd^\prec_mY_1$ into $\cd^\prec_mY_2$ and maps 
complexes of weight $\le i$ to complexes of weight $\le i$ (for any $i$). Assume that (a),(b) below hold:
$$(\Ph(\XX_0))^h\in\cm_m^\prec Y_2\text{ for any }\XX_0\in\cm_m^\prq Y_1\text{ and any }h>c;\tag a$$
$$\XX\text{ has weight $\le0$ and }\XX^i\in\cm^\prec Y_1\text{ for any }i>c'.\tag b$$
Then 
$$(\Ph(\XX))^j\in\cm^\prec Y_2\text{ for any }j>c+c',\tag c$$
and we have canonically
$$\un{(\Ph(\un{\XX^{\{c'\}}}))^{\{c\}}}=\un{(\Ph(\XX))^{\{c+c'\}}}.\tag d$$
\endproclaim
From the distinguished triangle $(\t_{<i}\XX,\t_{\le i}\XX,\XX^i[-i])$ we get a distinguished triangle 
$(\Ph(\t_{<i}\XX),\Ph(\t_{\le i}\XX),\Ph(\XX^i[-i])$; hence we have an exact sequence
$$(\Ph(\XX^i))^{h-1}@>>>(\Ph(\t_{<i}\XX))^{i+h}@>>>(\Ph(\t_{\le i}\XX))^{i+h}
@>>>(\Ph(\XX^i))^h@>>>(\Ph(\t_{<i}\XX))^{i+h+1}.$$
From this and (a),(b) we see by induction on $i$ that 

$(\Ph(\t_{\le i}\XX))^{i+h}\in\cm^\prec Y_2$ if $i+h>c+c'$ (in particular, 
$(\Ph(\XX))^k\in\cm^\prec Y_2$ if $k>c+c'$ so that (c) holds);

$(\Ph(\t_{\le c'}\XX))^{c+c'}@>\b>>(\Ph(\XX^{c'}))^c$ has kernel and cokernel in $\cm^\prec Y_2$;

$(\Ph(\t_{\le i}\XX))^{c+c'}@>\b'>>(\Ph(\t_{\le i+1}\XX))^{c+c'}$ has kernel and cokernel in 
$\cm^\prec Y_2$ for $i\ge c'$.
\nl
Here the maps $\b,\b'$ come from the previous exact sequence. Now $\b,\b'$ are strictly compatible with the 
weight filtrations (see \cite{\BBD, 5.3.5}); we deduce that the maps
$$gr_{c+c'}(\Ph(\t_{\le c'}\XX))^{c+c'}@>>>gr_{c+c'}(\Ph(\XX^{c'}))^c,$$
$$gr_{c+c'}(\Ph(\t_{\le i}\XX))^{c+c'}@>>>gr_{c+c'}(\Ph(\t_{\le i+1}\XX))^{c+c'}\text{ (for }i\ge c')$$
induced by $\b,\b'$ have kernel and cokernel in $\cm^\prec Y_2$. Since these are maps between semisimple 
perverse sheaves we see that they induce isomorphisms
$$\un{gr_{c+c'}(\Ph(\t_{\le c'}\XX))^{c+c'}}@>\si>>\un{gr_{c+c'}(\Ph(\XX^{c'}))^c},$$
$$\un{gr_{c+c'}(\Ph(\t_{\le i}\XX))^{c+c'}}@>\si>>\un{gr_{c+c'}(\Ph(\t_{\le i+1}\XX))^{c+c'}}\text{ (for }
i\ge c').$$
By composition we get a canonical isomorphism 
$$\un{gr_{c+c'}(\Ph(\XX^{c'}))^c}@>\si>>\un{gr_{c+c'}(\Ph(\XX))^{c+c'}};\tag e$$
(note that $\un{gr_{c+c'}(\Ph(\XX))^{c+c'}}=\un{gr_{c+c'}(\Ph(\t_{\le i}\XX))^{c+c'}}$ for $i\gg0$).

For any $j$ we have an exact sequence 
$$0@>>>\cw^{j-1}(\XX^{c'})@>>>\cw^j(\XX^{c'})@>>>gr_j\XX^{c'}@>>>0$$ 
hence a distinguished triangle 
$$(\Ph(\cw^{j-1}(\XX^{c'})),\Ph(\cw^j(\XX^{c'})),\Ph(gr_j\XX^{c'}))$$ 
which gives rise to an exact sequence 
$$\align&(\Ph(gr_j\XX^{c'}))^{c-1}@>>>(\Ph(\cw^{j-1}(\XX^{c'}))^c@>>>(\Ph(\cw^j(\XX^{c'})))^c\\&
@>>>(\Ph(gr_j\XX^{c'}))^c@>>>(\Ph(\cw^{j-1}(\XX^{c'})))^{c+1}\endalign$$
and to an exact sequence
$$\align&gr_{c+c'}(\Ph(gr_j\XX^{c'}))^{c-1}@>>>gr_{c+c'}(\Ph(\cw^{j-1}(\XX^{c'})))^c@>>>
gr_{c+c'}(\Ph(\cw^j(\XX^{c'})))^c\\&@>>>gr_{c+c'}(\Ph(gr_j\XX^{c'}))^c@>>>
gr_{c+c'}(\Ph(\cw^{j-1}(\XX^{c'}))^{c+1}.\endalign$$
Now $\Ph(\cw^j(\XX^{c'}))$ is mixed of weight $\le j$ hence $(\Ph(\cw^j(\XX^{c'})))^c$ is mixed of weight 
$\le c+j$ so that $gr_{c+c'}(\Ph(\cw^j(\XX^{c'})))^c=0$ if $j<c'$. Moreover 
$gr_{c+c'}(\Ph(gr_j\XX^{c'}))^c=0$ if $j>c'$ since $\XX^{c'}$ is mixed of weight $\le c'$. Thus we have an 
exact sequence
$$0@>>>gr_{c+c'}(\Ph(\cw^{c'}(\XX^{c'})))^c@>>>gr_{c+c'}(\Ph(gr_{c'}\XX^{c'}))^c@>>>
gr_{c+c'}(\Ph(\cw^{c'-1}(\XX^{c'})))^{c+1}$$
and we have
$$gr_{c+c'}(\Ph(\cw^{c'}(\XX^{c'})))^c=gr_{c+c'}(\Ph(\cw^{c'+1}(\XX^{c'})))^c
=gr_{c+c'}(\Ph(\cw^{c'+2}(\XX^{c'}))^c=\do.$$
Thus we have an exact sequence 
$$0@>>>gr_{c+c'}(\Ph(\XX^{c'}))^c@>>>gr_{c+c'}(\Ph(gr_{c'}\XX^{c'}))^c@>>>
gr_{c+c'}(\Ph(\cw^{c'-1}(\XX^{c'})))^{c+1}.$$
By (a) we have $(\Ph(\cw^{c'-1}(\XX^{c'})))^{c+1}\in\cm^\prec Y_2$ hence
$$gr_{c+c'}(\Ph(\cw^{c'-1}(\XX^{c'})))^{c+1}\in\cm^\prec Y_2.$$
Thus $gr_{c+c'}(\Ph(\XX^{c'}))^c$ is a subobject of $gr_{c+c'}(\Ph(gr_{c'}\XX^{c'}))^c$ and the quotient is 
in $\cm^\prec Y_2$. Since $gr_{c+c'}(\Ph(gr_{c'}\XX^{c'}))^c$ is semisimple in $\cm(Y_2)$ it follows that
$$\un{gr_{c+c'}(\Ph(\XX^{c'}))^c}=\un{gr_{c+c'}(\Ph(gr_{c'}\XX^{c'}))^c}.$$ 
This, together with (e) gives
$$\un{gr_{c+c'}(\Ph(gr_{c'}\XX^{c'}))^c}=\un{gr_{c+c'}(\Ph(\XX))^{c+c'}}.$$
It follows that
$$\un{gr_{c+c'}(\Ph(\un{gr_{c'}\XX^{c'}}))^c}=\un{gr_{c+c'}(\Ph(\XX))^{c+c'}}$$
so that (d) holds.

\subhead 1.13\endsubhead
Let $L\in\cc^\boc_0\cb^2$. We have clearly $\fD(L)\in\cc^\boc_0\cb^2$. We show that we have canonically:
$$\unc(\fD(L))=\fD(\unc(L)).\tag a$$
By the relative hard Lefschetz theorem \cite{\BBD, 5.4.10} applied to the projective morphism $\p$ and to
$f^*L[[\nu+\r]]$ (a perverse sheaf of pure weight $0$ on $X$, notation of 1.5) we have canonically for any 
$i$:
$$(\p_!f^*L[[\nu+\r]])^{-i}=(\p_!f^*L[[\nu+\r]])^i(i).\tag b$$ 
We have used that $f$ is smooth with fibres of dimension $\nu+\r$. This also shows that 
$$\fD(\c(\fD(L)))=\c(L)[[2\nu+2\r]].\tag c$$
Using (b),(c) we have
$$\align&\fD(\unc(\fD(L)))=\fD((\c(\fD(L)))^{a+\nu-r}((a+\nu+\r)/2))\\&=
(\fD(\c(\fD(L))))^{-a-\nu-\r}((-a-\nu-\r)/2))\\&=
(\c(L)[[2\nu+2\r]])^{-a-\nu-\r}((-a-\nu-\r)/2))=(\c(L)[[\nu+\r]])^{-a}(-a/2))\\&=
(\c(L)[[\nu+\r]])^{a}(a/2))=(\c(L))^{a+\nu+\r}((a+\nu+\r)/2))=\unc L.\endalign$$
This proves (a).

\subhead 1.14\endsubhead
Let $d\in\DD_\boc$ and let $\L_d$ be the left cell containing $d$. We show:
$$(A:\unc(\LL_d))=(-1)^{a+|d|}c_{A,[\L_d]}\text{ for any }A\in CS_\boc.\tag a$$
For any $E\in\Irr_\boc W$ we have $\tr(t_d,E_\iy)=\text{ multiplicity of $E$ in }[\L_d]$.
Hence, using 1.10(a) and 1.7(a), we have
$$\align&
(A:\unc(\LL_d))=(-1)^{a+|d|}\sum_{E\in\Irr_\boc W}c_{A,E}(\text{ multiplicity of $E$ in $[\L_d]$})\\&=
(-1)^{a+|d|}c_{A,[\L_d]}.\endalign$$
It remains to show that 
$$|d|=a\mod2.\tag b$$
If $p_{1,d}\in\ZZ[v\i]$ is as in \cite{\HEC, 5.3}, then $v^{-a}$ appears with nonzero coefficient in
$p_{1,d}$, see \cite{\HEC, 14.1} hence by \cite{\HEC, 5.4(b)} we have $-a=|d|-|1|\mod2$. This proves (b)
hence (a).

\subhead 1.15\endsubhead
Let $\p_1:\{(B,g)\in\cb\T G;g\in B\}@>>>G$ be the first projection. Let $\Si:=\p_{1!}\bbq[[\D]]$. As observed
in \cite{\LGR}, $\p_1$ is small, so that $\Si$ is a perverse sheaf on $G$; moreover, $\Si$ has a natural 
$W$-action so that $\Si=\op_{E\in\Irr W}E\ot A_E$ where $A_E=\Hom_W(E,\Si)$ is a simple perverse sheaf. 
Since $\Si=\c(\LL_1)[[\nu+\r]]$ we have $A_E\in\cc^\spa_0G$ for any $E$. It is known that $A_E\in\cm^\prq G$
if and only if $\boc_E\prec\boc$ and $A_E\in\cc^\boc G$ if and only if $\boc_E=\boc$. There is a unique 
$E_\boc\in\Irr_\boc W$ such that $E_\boc^\da$ is a special representation of $W$.  

We show:

(a) {\it Assume that $(A_{E_\boc}:\unc(\LL_d))\le1$ for any $d\in\DD_\boc$. Then for any $d\in\DD_\boc$ we 
have $(A_{E_\boc}:\unc(\LL_d))=1$.}
\nl
For any $d\in\DD_\boc$ we set $\d(d)=c_{A_{E_\boc},[\L_d]}$. By 1.14 our assumption is that 
$\d(d)\in\{0,1\}$ for any $d\in\DD_\boc$ and we must prove that $\d(d)=1$ for any $d\in\DD_\boc$. Since 
$c_{A_{E_\boc},[\boc]}=\sum_{d\in\DD_\boc}\d(d)$, it is enough to show that
$c_{A_{E_\boc},[\boc]}=|\DD_\boc|$. Since $\boc_{A_{E_\boc}}=\boc$ we have
$c_{A_{E_\boc},[\boc']}=0$ for any two-sided cell $\boc'\ne\boc$. Hence it is enough to show that
$\sum_{\boc'}c_{A_{E_\boc},[\boc']}=|\DD_\boc|$ where $\boc'$ runs over the two-sided cells in $W$.
Let $\Reg$ be the regular representation of $W$. We have $\sum_{\boc'}c_{A_{E_\boc},[\boc']}=
c_{A_{E_\boc},\Reg}$ hence it is
enough to show that $c_{A,\Reg}=|\DD_\boc|$ where $A=A_{E\boc}$. From 1.5(a) we have
$$(A:R_1^\D)=(0;\sum_{E\in\Irr W}c_{A,E}\dim(E))=\sum_{E\in\Irr W}c_{A,E}\dim(E)=c_{A,\Reg}$$
hence it is enough to show that $(A:R_1^\D)=|\DD_\boc|$.
Recall that $R_1[\D]=\Si$ hence it is enough to show that $(A:\Si)=|\DD_\boc|$. We have
$(A:\Si)=\dim E_\boc$. It remains to show that $\dim(E_\boc)=|\DD_\boc|$. This is a well known property of 
special representations.

\mpb

We will see in 6.4 that the assumption of (a) is in fact satisfied.

\head 2. Truncated restriction\endhead
\subhead 2.1\endsubhead
The following result and its proof are similar to 1.6. 

\proclaim{Lemma 2.2}Let $r\ge1$, $J\sub[1,r]$, $J\ne\emp$ and $\ww=(w_1,w_2,\do,w_r)\in W^r$. Let 
$\fF=\nu+(r-1)a$.

(a) Assume that $w_i\in\boc$ for some $i\in[1,r]$. If $j\in\ZZ$ (resp. $j>\fF$) then 
$(p_{0r!}L_\ww^{[1,r]}[|\ww|])^j$ is in $\cm^\prq\cb^2$ (resp. $\cm^\prec\cb^2$).

(b) Assume that $w_i\in\boc$ for some $i\in J$. If $j\in\ZZ$ (resp. $j\ge\fF$) then 
$(p_{0r!}\dL_\ww^J[|\ww|])^j$ is in $\cm^\prq\cb^2$ (resp. $\cm^\prec\cb^2$).

(c) Assume that $w_i\in\boc$ for some $i\in J$. If $j\ge\fF$ then the cokernel of the map 
$$(p_{0r!}L_\ww^J[|\ww|])^j@>>>(p_{0r!}L_\ww^{[1,r]}[|\ww|])^j$$
associated to 1.1(a) is in $\cm^\prec\cb^2$.

(d) Assume that $w_i\in\boc$ for some $i\in J$. If $j\in\ZZ$ (resp. $j>\fF$) then 
$(p_{0r!}L_\ww^J[|\ww|])^j$ is in $\cm^\prq\cb^2$ (resp. $\cm^\prec\cb^2$).

(e) Assume that $w_i\prec\boc$ for some $i\in J$. If $j\in\ZZ$ then 
$(p_{0r!}L_\ww^{[1,r]}[|\ww|])^j\in\cm^\prec\cb^2$ and $(p_{0r!}L_\ww^J[|\ww|])^j\in\cm^\prec\cb^2$.
\endproclaim  
We prove (a). Let $L=\LL_z$, $z\in W$ and $j\in\ZZ$ be such that $L$ is a composition factor of 
$(p_{0r!}L_\ww^{[1,r]}[|\ww|])^j=(L_\ww^\cir[\ww])^j$. By 1.2(b) we have
$$(L_\ww^\cir[|\ww|])^j\cong\op_{w\in W,k\in\ZZ;j+k-\nu=0}(\LL_w)^{\op N(w,k)}$$
hence $N(z,\nu-j)=0$. From $N(z,\nu-j)\ne0$ and 1.4(a) we see that $z\prq\boc$. Assume now that $z\in\boc$. 
From 1.4 we see that $\nu-j\ge-(r-1)a$ that is $j\le\fF$.

We prove (b). Let $L=\LL_z$, $z\in W$ and $j\in\ZZ$ be such that $L$ is a composition factor of 
$(p_{0r!}\dL_\ww^J[|\ww|])^j$. There exists $h$ such that $L$ is a composition factor of 
$(p_{0r!}(\dL_\ww^J[|\ww|])^h)[-h])^j$. We have $(\dL_\ww^J[|\ww|])^h\ne0$ hence 
$(\dL_{\ww}^J[|\ww|+\nu-1])^{h-\nu+1}\ne0$ hence by 1.1(c), $h-\nu+1\le0$. From 1.1(b) we see that there 
exists $\ww'=(w'_1,w'_2,\do,w'_r)\in W^r$ such that $w_i=w'_i$ for all $i\in J$ and such that $L$ is a 
composition factor of 
$$(p_{0r!}(L_{\ww'}^{[1,r]}[|\ww'|+\nu])[-h])^j=(p_{0r!}(L_{\ww'}^{[1,r]}[|\ww'|]))^{j+\nu-h}
=(L_{\ww'}^\cir[|\ww'|])^{j+\nu-h}.$$
From (a) (for $\ww'$ instead of $\ww$) we see that $z\prq\boc$ and that $z\prec\boc$ if $j+\nu-h>\fF$
that is, if $j>h+\fF-\nu$. If $j\ge\fF$ then using $h-\nu+1\le0$ (that is $0>h-\nu$) we see that we have 
indeed $j>h+\fF-\nu$. This proves (b).

We prove (c). From 1.1(a) we get a distinguished triangle
$$(p_{0r!}L_\ww^J[[|\ww|]],p_{0r!}L_\ww^{[1,r]}[[|\ww|]],p_{0r!}\dL_\ww^J[[|\ww|]])$$
in $\cd_m(\cb^2)$. This gives rise for any $j$ to an exact sequence
$$(p_{0r!}\dL_\ww^J[[|\ww|]])^{j-1}@>>>(p_{0r!}L_\ww^J[[|\ww|]])^j
@>>>(p_{0r!}L_\ww^{[1,r]}[[|\ww|]])^j@>>>(p_{0r!}\dL_\ww^J[[|\ww|]])^j\tag f$$
in $\cm_m(\cb^2)$. Using this and (b) we see that (c) holds.

Now (d) follows from the previous exact sequence using (a),(b).

Replacing $\boc$ in (a) and (d) by the two-sided cell containing $w_i$ in (e) we deduce that (e) holds. 

\subhead 2.3\endsubhead
Let $r\ge1$ and let $x_1,x_2,\do,x_r$ be elements of $W$ such that at least one of them is in $\boc$. We 
show:

(a) {\it If $\un{(\LL_{x_1}\cir\LL_{x_2}\cir\do\cir\LL_{x_r})^{\{(r-1)(a-\nu)\}}}\ne0$ then $x_i\in\boc$
for all $i\in[1,r]$.}
\nl
By assumption we have 
$$\un{(L_{x_1}^\sh[[|x_1|]]\cir L_{x_2}^\sh[[|x_2|]]\cir\do\cir L_{x_r}^\sh[[|x_r|]])^{\{\nu+(r-1)a\}}}\ne0.
$$ 
Using 1.2(b) we see that there exists $w\in\boc$ such that $\LL_w$ appears with nonzero multiplicity in
$$\sum_{z\in W,n\in\ZZ}((L_z^\sh[n+|z|])^{(r-1)a+\nu})^{\op N_y(z,n)}$$
(that is, $N_y(w,-(r-1)a)\ne0$) where $N_y(z,n)\in\NN$ are given by the following identity in $\HH$:
$$c_{x_1}c_{x_2}\do c_{x_r}=\sum_{z\in W,n\in\ZZ}N_y(z,n)v^nc_z.$$
From $N_y(w,-(r-1)a)\ne0$ we see using 1.3(b) that $x_i\in\boc$ for all $i$. 

\subhead 2.4\endsubhead
In the setup of 2.3, we see using 1.3(d), that 

(a) $N_y(w',-(r-1)a)$ is the coefficient of $t_{w'}$ in $t_{w_1}t_{w_2}\do t_{w_r}$. 

\subhead 2.5\endsubhead
Let $\p,f$ be as in 1.5. Now $K\m\z(K)=f_!\p^*K$ defines a functor $\cd_m(G)@>>>\cd_m(\cb^2)$. For 
$i\in\ZZ,K\in\cd_m(G)$ we write $\z^i(K)$ instead of $(\z(K))^i$.

A functor closely related to $\z$ (in which a complex $K$ on $G$ was integrated over the cosets of the 
unipotent radical of a Borel subgroup, rather than over the cosets of a Borel subgroup as in $\z$) was 
introduced in \cite{\MV} and by the author in 1987 (unpublished, but mentioned in \cite{\MV, \S5} and 
\cite{\GI, \S0}) when I found a criterion for $K$ to be a character sheaf in terms of the cohomology sheaves
of the image of $K$ under this functor. My proof of that criterion was based in part on something close to 
the following result, a version of which (at the level of Grothendieck groups) appears also in 
\cite{\GR, (3.3.1)}.

\proclaim{Proposition 2.6} For any $L\in\cd_m(\cb^2)$ we have
$$\z(\c(L))\Bpq\{\op_{y\in W;|y|=k}L_y\cir L\cir L_{y\i}\ot\fL[[2k-2\nu]];k\in\NN\},\tag a$$
$$\align&\z(\c(L))\Bpq
\\&\{\op_{y\in W;|y|=k}L_y\cir L\cir L_{y\i}\ot\fL[[2k-2\nu-2\r]]\ot\L^d\cx[[d]](d/2);k\in\NN,d\in[0,\r]\},
\tag b\endalign$$
where $\fL,\cx$ are as in 0.2.
\endproclaim
Let
$$Y=\{(B_1,B_2,B_3,B_4,g)\in\cb\T\cb\T\cb\T\cb\T G;gB_1g\i=B_4,gB_2g\i=B_3\}.$$
For $ij=14$ or $23$ we define $h'_{ij}:Y@>>>X$ by $(B_1,B_2,B_3,B_4,g)\m(B_i,B_j,g)$ and 
$h_{ij}:Y@>>>\cb^2$ by $(B_1,B_2,B_3,B_4,g)\m(B_i,B_j)$. We have $\p^*\p_!=h'_{14!}h'_{23}{}^*$ hence 
$$\z(\c(L))=f_!\p^*\p_!f^*(L)=f_!h'_{14!}h'_{23}{}^*f^*(L)=h_{14!}h_{23}^*L.$$
For $k\in\NN$ let $Y^k=\cup_{y\in W;|y|=k}Y_y$ where
$$Y_y=\{(B_1,B_2,B_3,B_4,g)\in Y;(B_1,B_2)\in\co_y,(B_3,B_4)\in\co_{y\i}\}$$
and let $Y^{\ge k}:=\cup_{k';k'\ge k}Y^{k'}$, an open subset of $Y$;
let $h_{ij}^k:Y^k@>>>\cb^2$, $h_{ij}^{\ge k}:Y^{\le k}@>>>\cb^2$ be the restrictions of $h_{ij}$. 
For any $k\in\NN$ we have a distinguished triangle
$$(h_{14!}^{\ge k+1}h_{23}^{\ge k+1*}L),h_{14!}^{\ge k}h_{23}^{\ge k*}L,h_{14!}^kh_{23}^{k*}L).$$
It follows that we have
$$\z(\c(L))\Bpq\{h_{14!}^kh_{23}^{k*}L;k\in\NN\}.$$
For $k\in\NN$ let $Z^k=\cup_{y\in W;|y|=k}Z_y$ where
$$Z_y=\{(B_1,B_2,B_3,B_4)\in\cb^4;(B_1,B_2)\in\co_y,(B_3,B_4)\in\co_{y\i}\};$$
for $i,j\in[1,4]$ we define $\tih_{ij}^k:Z^k@>>>\cb^2$ and $\tih_{ij}^y:Z_y@>>>\cb^2$
by $(B_1,B_2,B_3,B_4)\m(B_i,B_j)$.
We have an obvious morphism $u:Y^k@>>>Z^k$ whose fibres are isomorphic to $\kk^{\nu-k}$ times the
$\r$-dimensional torus $T$. We have a commutative diagram
$$\CD
\cb^2@<h_{23}^k<<Y^k@>h_{14}^k>>\cb^2\\
@V1VV     @VuVV   @V1VV \\
\cb^2@<\tih_{23}^k<<Z^k@>\tih_{14}^k>>\cb^2\endCD$$
We have
$$h_{14!}^kh_{23}^{k*}L=\tih_{14!}^ku_!u^*\tih_{23}^{k*}L=\tih_{14!}^k(\tih_{23}^{k*}L\ot u_!\bbq)=
(\tih_{14!}^k\tih_{23}^{k*}L)\ot\fL[[-2\nu+2k]].$$  
We deduce that
$$\z(\c(L))\Bpq\{(\tih_{14!}^k\tih_{23}^{k*}L)\ot\fL[[-2\nu+2k]];k\in\NN\}.$$
Since $Z^k$ is the union of open and closed subvarieties $Z_y,|y|=k$, we have
$$\tih_{14!}^k\tih_{23}^{k*}L=\op_{y\in W;|y|=k}\tih_{14!}^y\tih_{23}^{y*}L.$$
From the definitions we have
$$\tih_{14!}^y\tih_{23}^{y*}L=L_y\cir L\cir L_{y\i}.$$
This completes the proof of (a). Now (b) follows from (a) using 
$$\fL[[2\r]]\Bpq\{\bbq\ot\L^d\cx[[d]](d/2);d\in[0,\r]\}\tag c$$
which follows from the definitions.

\proclaim{Proposition 2.7}Let $w\in W$ and let $j\in\ZZ$. We set $S=\z(R_w)[[2\r+2\nu+|w|]]\in\cd_m(\cb^2)$.

(a) If $w\prq\boc$ then $S^j\in\cm^\prq\cb^2$.

(b) If $w\in\boc$ and $j>\nu+2a$ then $S^j\in\cm^\prec\cb^2$.

(c) If $w\prec\boc$ then $S^j\in\cm^\prec\cb^2$.

(d) $S^j$ is mixed of weight $\le j$.

(e) If $j\ne\nu+2a$ and $w\in\boc$ then $gr_{\nu+2a}S^j\in\cm^\prec\cb^2$.

(f) If $k>\nu+2a$ and $w\in\boc$ then $gr_kS^j\in\cm^\prec\cb^2$.
\endproclaim
Let $J=\{2\}\sub[1,3]$. Using 2.5 and 1.2(a) with $r=3$ we have
$$S\Bpq\{p_{03!}L^J_{y,w,y\i}[[|w|+2|y|]])\ot\L^d\cx[[d]](d/2);d\in[0,\r],y\in W\}.\tag g$$
Using this and the definitions we see that to prove (a) it is enough to show that for any $y,d$ as above we
have
$$(p_{03!}L^J_{y,w,y\i}[[|w|+2|y|]]\ot\L^d\cx[[d]](d/2))^j\in\cm^\prq\cb^2;\tag h$$
this follows from 2.2(d),(e). This proves (a). 

At the same time we see that to prove (d) it is enough to show that for any $y,d$ as above, (h) is mixed of
weight $\le j$. Since $\bbq[[d]](d/2)$ is pure of weight $-d\le0$, to prove the last statement it is enough 
to show that $p_{03!}L^J_{y,w,y\i}[[|w|+2|y|]]$ is mixed of weight $\le0$.
Note that $L^J_{y,w,y\i}[[|w|+2|y|]]$ is obtained by $()_!$ under an open imbedding from
$L^{[1,3]}_{y,w,y\i}[[|w|+2|y|]]$ which is pure of weight $0$ hence it is
mixed of weight $\le0$ (see \cite{\BBD, 5.1.14}), hence 
$p_{03!}L^J_{y,w,y\i}[[|w|+2|y|]]$ is mixed of weight $\le0$ (see \cite{\BBD, 5.1.14}).
This proves (d).

We prove (b). It is again enough to show that for any $y,d$ as above
$$(p_{03!}L^J_{y,w,y\i}[[|w|+2|y|]]\ot\L^d\cx{d}[[d]](d))^j$$
is in $\cm^\prec\cb^2$ if $j>\nu+2a$. This follows from 2.2(d) since $j>\nu+2a$, $d\ge0$ implies 
$j+d>\nu+2a$.

Now (c) follows from (a) by replacing $\boc$ by the two-sided cell containing $w$.

We prove (e). If $j>\nu+2a$ this follows from (b). If $j<\nu+2a$ we have $gr_{\nu+2a}S^j=0$ by (a). This 
proves (e).

We prove (f). If $j<k$ we have $gr_kS^j=0$ by (d). If $j\ge k$ we have $j>\nu+2a$ so that 
$S^j\in\cm^\prec\cb^2$ by (b). This proves (f).

\proclaim{Proposition 2.8} (a) If $K\in\cd^\prq G$ then $\z(K)\in\cd^\prq\cb^2$. If $K\in\cd^\prec G$ then 
$\z(K)\in\cd^\prec\cb^2$.

(b) If $K\in\cm^\prq G$ and $j>\r+\nu+a$ then $\z^j(K)\in\cm^\prec\cb^2$.
\endproclaim
It is enough to prove the proposition assuming in addition that $K=A\in CS(G)$. By 1.7(b) we can find 
$w\in\boc_A$ such that $(A:R_w^{n_w})\ne0$. Then $A[-n_w]$ is a direct summand of $R_w$. Hence $\z(A)$ is a 
direct summand of $\z(R_w)[\D+a+|w|]$ and $\z^j(A)$ is a direct summand of
$\z^{j+\D+a+|w|}(R_w)=\z^{j-\r+a}(R_w[2\r+2\nu+|w|])$. Using 2.7 we deduce that (a) holds and that, in the
setup of (b), $\z^j(A)\in\cm^\prec\cb^2$ provided that $j-\r+a>\nu+2a$. Hence (b) holds.

\subhead 2.9\endsubhead
For $K\in\cc^\boc_0G$ we set 
$$\unz(K)=\un{(\z(K))^{\{\r+\nu+a\}}}\in\cc^\boc_0\cb^2.$$
We say that $\unz(K)$ is the {\it truncated restriction} of $K$.

\subhead 2.10\endsubhead
We note the following result, a version of which was first stated in \cite{\GI, 9.2.1}.

(a) {\it Let $K\in\cd_m(G)$ and let $L\in\cm_m(\cb^2)$ be $G$-equivariant. Then there is a canonical 
isomorphism $L\cir\z(K)@>\si>>\z(K)\cir L$.}
\nl
We have $\z(K)\cir L=c_!d^*(K\bxt L)$, $L\cir\z(K)=c'_!d'{}^*(K\bxt L)$ where  
$$Z=\{(g,B,B'',B')\in G\T\cb\T\cb\T\cb;gBg\i=B''\},$$
$$Z'=\{(g,B,B'',B')\in G\T\cb\T\cb\T\cb;g\i B'g=B''\},$$

$d:Z@>>> G\T\cb^2$ is $(g,B,B'',B')\m(g,(B'',B'))$, 

$d':Z'@>>> G\T\cb^2$ is $(g,B,B'',B')\m(g,(B,B''))$,

$c:Z@>>>\cb^2$, $c':Z'@>>>\cb^2$ are $(g,B,B'',B')\m(B,B')$.
\nl
We identify $Z,Z'$ with $G\T\cb^2$ by 
$(g,B,B'',B')\m(g,(B,B'))$. Then $d$ becomes $d_1:(g,(B,B'))\m(g,(gBg\i,B'))$, $d'$ becomes 
$d'_1:(g,(B,B'))\m(g,(B,g\i B'g))$ and $c,c'$ become $c_1:(g,(B,B'))\m(B,B')$. It is enough to show that 
$d_1^*(K\bxt L)=d'_1{}^*(K\bxt L)$. Define $u:G\T\cb^2@>>>G\T\cb^2$ by $(g,(B,B'))\m(g,(gBg\i,gB'g\i))$. By
the $G$-equivariance of $L$ we have canonically $u^*(\bbq\bxt L)=\bbq\bxt L$. We have  $d_1=ud'_1$ hence 
$d_1^*(K\bxt L)=d'_1{}^*u^*(K\bxt L)=d'_1{}^*(K\bxt L)$ and (a) follows.

\proclaim{Proposition 2.11} (a) If $L\in\cm^\prq\cb^2$ and $j>2a+2\nu+2\r$ then 
$(\z(\c(L)))^j\in\cm^\prec\cb^2$.

(b) If $L\in\cc^\boc_0\cb^2$, we have canonically 
$$\unz(\unc(L))=\un{(\z(\c(L)))^{\{2a+2\nu+2\r\}}}\in\cc^\boc_0\cb^2.$$
\endproclaim
We apply 1.12 with $\Ph=\z:\cd_m(G)@>>>\cd_m(\cb^2)$ and with $\XX=\c(L)$,
$(c,c')=(a+\nu+\r,a+\nu+\r)$, see 2.8, 1.9. The result follows.

\head 3. Truncated convolution on $\cb^2$\endhead
\subhead 3.1\endsubhead
We show that for $L,L'\in\cd^\spa\cb^2$, (a) and (b) below hold.

(a) {\it If $L\in\cd^\prq\cb^2$ or $L'\in\cd^\prq\cb^2$ then $L\cir L'\in\cd^\prq\cb^2$. If 
$L\in\cd^\prec\cb^2$ or $L'\in\cd^\prec\cb^2$ then $L\cir L'\in\cd^\prec\cb^2$.}

(b) {\it Assume that $L,L'\in\cm^\spa\cb^2$ and that either $L$ or $L'$ is in $\cm^\prq\cb^2$. If $j>a-\nu$ 
then $(L\cir L')^j\in\cm^\prec\cb^2$.}
\nl
We can assume that $L=\LL_z,L'=\LL_{z'}$ with $z\prq\boc$ or $z'\prq\boc$. Then (a) follows from 1.4(b),(c).
To prove (b) we can further assume that $z\in\boc$ or $z'\in\boc$. According to 1.4(b) we have 
$(L_z^\sh[|z|]\cir L_{z'}^\sh[|z'|])^{j'}\in\cm^\prec\cb^2$ if $j'>\nu+a$ hence \lb
$(L_z^\sh[|z|+\nu]\cir L_{z'}^\sh[|z'|+\nu])^j\in\cm^\prec\cb^2$ if $j+2\nu>\nu+a$ that is if $j>a-\nu$; 
this proves (b).

\subhead 3.2\endsubhead
For $L,L'\in\cc^\boc_0\cb^2$, we set 
$$L\unb L'=\un{(L\cir L')^{\{a-\nu\}}}\in\cc^\boc_0\cb^2.\tag a$$
Using 1.12 twice, we see that for $L,L',L''\in\cc^\boc_0\cb^2$ we have canonically
$$(L\unb L')\unb L''=\un{(L\cir L'\cir L'')^{\{2a-2\nu\}}},$$
$$L\unb(L'\unb L'')=\un{(L\cir L'\cir L'')^{\{2a-2\nu\}}}.$$
Hence 
$$(L\unb L')\unb L''=L\unb(L'\unb L'').$$
We see that $L,L'\m L\unb L'$ defines an associative tensor product structure on $\cc^\boc_0\cb^2$. (A 
closely related result appears in \cite{\CEL}.) Hence if ${}^1L,{}^2L,\do,{}^rL$ are in $\cc^\boc_0\cb^2$ 
then ${}^1L\unb{}^2L\unb\do\unb{}^rL\in\cc^\boc_0\cb^2$ is well defined. Using 1.12 repeatedly, we have
$${}^1L\unb{}^2L\unb\do\unb{}^rL=\un{({}^1L\cir{}^2L\cir\do\cir{}^rL)^{\{(r-1)(a-\nu)\}}}.\tag b$$

\subhead 3.3\endsubhead
Let $L,L'\in\cc^\boc_0\cb^2$. We show that we have canonically:
$$\fD(L\unb L')=\fD(L)\unb\fD(L').\tag a$$
We can assume that $L=\LL_{w_1},L'=\LL_{w_2}$ where $w_1,w_2\in\boc$. Let $\ww=(w_1,w_2)$. Let 
$L^{[1,2]}_\ww$ be the intersection cohomology complex of the projective variety
$$\{(B_0,B_1,B_2)\in\cb^3;(B_0,B_1)\in\bco_{w_1},(B_1,B_2)\in\bco_{w_2}\}$$
extended by $0$ on the complement to this variety in $\cb^3$ and let $p_{02}:\cb^3@>>>\cb^2$ be the map
$(B_0,B_1,B_2)\m(B_0,B_2)$. By definition we have 
$$L\cir L'=p_{02!}L^{[1,2]}_\ww[[|w_1|+|w_2|+2\nu]].$$
We must show that
$$\fD((L\cir L')^{a-\nu}((a-\nu)/2)))=(L\cir L')^{a-\nu}((a-\nu)/2)).\tag b$$
By the hard Lefschetz theorem \cite{\BBD, 5.4.10} applied to the projective morphism $p_{02}$ and to 
$L^{[1,2]}_\ww[[|w_1|+|w_2|+\nu]]$ (a perverse sheaf of pure weight $0$ on $\cb^3$) we have canonically for
any $i$:
$$(p_{02!}L^{[1,2]}_\ww[[|w_1|+|w_2|+\nu]])^{-i}=(p_{02!}L^{[1,2]}_\ww[[|w_1|+|w_2|+\nu]])^i(i)$$
that is $(L\cir L'[[-\nu]])^{-i}=(L\cir L'[[-\nu]])^i(i)$, hence
$$(L\cir L')^{-i-\nu}=(L\cir L')^{i-\nu}(i).\tag c$$
We have $\fD(L\cir L')=L\cir L'[[-2\nu]]$ hence $\fD((L\cir L')^i)=(L\cir L')^{-i-2\nu}(-\nu)$.
Thus
$\fD((L\cir L')^{a-\nu})=(L\cir L')^{-\nu-a}(-\nu)=(L\cir L')^{-\nu+a}(a-\nu)$.
(The last equality uses (c).) This proves (b), hence (a).

\mpb

The following result is a truncated version of 2.10.

\proclaim{Proposition 3.4}Let $K\in\cc^\boc_0G,L\in\cc^\boc_0\cb^2$. There is a canonical isomorphism
$$L\unb\unz(K)@>\si>>\unz(K)\unb L.\tag a$$
\endproclaim
Applying 1.12 with $\Ph:\cd_m^\prq\cb^2@>>>\cd_m^\prq\cb^2$, $L'\m L'\cir L$, $\XX=\z(K)$,
$(c,c')=(a-\nu,a+\r+\nu)$ (see 3.1, 2.8), we deduce that we have canonically 
$$\un{((\z(K))^{\{a+\r+\nu\}}\cir L)^{\{a-\nu\}}}=\un{(\z(K)\cir L)^{\{2a+\r\}}}.\tag b$$
Using 1.12 with $\Ph:\cd_m^\prq\cb^2@>>>\cd_m^\prq\cb^2$, $L'\m L\cir L'$, $\XX=\z(K)$,
$(c,c')=(a-\nu,a+\r+\nu)$ (see 3.1, 2.8), we deduce that we have canonically 
$$\un{(L\cir(\z(K))^{\{a+\r+\nu\}})^{\{a-\nu\}}}=\un{(L\cir\z(K))^{\{2a+\r\}}}.\tag c$$
We now combine (b),(c) with 2.10(a); we obtain the isomorphism (a).

\head 4. Truncated convolution on $G$\endhead
\subhead 4.1\endsubhead
Let $\mu:G\T G@>>>G$ be the multiplication map. For $K,K'\in\cd_m(G)$ we define the 
convolution $K*K'\in\cd_m(G)$ by $K*K'=\mu_!(K\bxt K')$. For $K,K',K''$ in 
$\cd_m(G)$ we have canonically $(K*K')*K''=K*(K'*K'')$ (and we denote this by $K*K'*K''$).

Note that if $K\in\cd_m(G)$ and $K'\in\cm_m(G)$ is $G$-equivariant for the conjugation action of $G$ then we
have a canonical isomorphism
$$K*K'@>\si>>K'*K.\tag a$$
Define $r:G\T G@>>>G$, $p_1:G\T G@>>>G$, $p_2:G\T G@>>>G$ by $r:(x,y)\m x\i yx$, $p_1:(x,y)\m x$,
$p_2:(x,y)\m y$. Without any assumption on $K'$ we have
$$\mu_!(p_1^*K\ot r^*K')=\mu_!(p_2^*K\ot p_1^*K')=K'*K.$$
In our case we have canonically $r^*K'=p_2^*K'$. Hence 
$$\mu_!(p_1^*K\ot r^*K')=\mu_!(p_1^*K\ot p_2^*K')=K*K'$$ 
and (a) follows.

\proclaim{Lemma 4.2} Let $K\in\cd_m(G)$, $L\in\cd_m(\cb^2)$. We have canonically $K*\c(L)=\c(L\cir\z(K))$.
\endproclaim
Let $Z=\{(g_1,g_2,B,B')\in G\T G\T\cb\T\cb;g_2Bg_2\i=B'\}$. Define $c:Z@>>>G\T\cb^2$ by
$(g_1,g_2,B,B')\m(g_1,(B,B'))$ and $d:Z@>>>G$ by $(g_1,g_2,B,B')\m g_1g_2$. From the definitions we see that
both $K*\c(L),\c(L\cir\z(K))$ can be identified with $d_!c^*(K\bxt L)$. The lemma follows.

\proclaim{Proposition 4.3} For any $L,L'\in\cd_m(\cb^2)$ we have
$$\align&\c(L)*\c(L')[[2\r+2\nu]]\\&\Bpq  
\{\c(L'\cir L_y\cir L\cir L_{y\i})[[2|y|]]\ot\L^d\cx[[d]](d/2);d\in[0,\r],y\in W\}.\endalign$$
\endproclaim
From 2.6(b) we deduce
$$\align&L'\cir\z(\c(L))[[2\nu+2\r]]\\&
\Bpq\{L'\cir L_y\cir L\cir L_{y\i}[[2|y|]]\ot\L^d\cx[[d]](d/2);y\in W,d\in[0,\r]\}\endalign$$
and
$$\align&\c(L'\cir\z(\c(L)))[[2\nu+2\r]]\\&
\Bpq\{\c(L'\cir L_y\cir L\cir L_{y\i})[[2|y|]]\ot\L^d\cx[[d]](d/2);y\in W,d\in[0,\r]\}.\endalign$$
It remains to show that $\c(L'\cir\z(\c(L)))=\c(L)*\c(L')$. This follows from 4.2 with $K,L$ replaced by
$\c(L),L'$.

\proclaim{Proposition 4.4} Let $w,w'\in W$ and let $j\in\ZZ$. We set 
$C=R_w*R_{w'}[[2\r+2\nu+|w|+|w'|]]\in\cd_m(G)$.

(a) If $w\prq\boc$ or $w'\prq\boc$ then $C^j\in\cm^\prq G$.

(b) If $j>\D+4a$ and either $w\in\boc$ or $w'\in\boc$ then $C^j\in\cm^\prec G$.

(c) If $w\prec\boc$ or $w'\prec\boc$ then $C^j\in\cm^\prec G$.

(d) $C^j$ is mixed of weight $\le j$.

(e) If $j\ne\D+4a$ and either $w\in\boc$ or $w'\in\boc$ then $gr_{\D+4a}C^j\in\cm^\prec G$.

(f) If $k>\D+4a$ and $w\in\boc$ or $w'\in\boc$ then $gr_k C^j\in\cm^\prec G$.
\endproclaim
Let $J=\{1,3\}\sub[1,4]$. Using 4.3 and 1.2(a) with $r=4$ we have
$$C\Bpq\{\c(p_{04!}L^J_{w',y,w,y\i}[[|w|+|w'|+2|y|]])\ot\L^d\cx[[d]](d/2);d\in[0,\r],y\in W\}.\tag g$$
Using this and the definitions we see that to prove (a) it is enough to show that for any $y,d$ as above,
$$\c^j(p_{04!}L^J_{w',y,w,y\i}[[|w|+|w'|+2|y|]]\ot\L^d\cx[[d]](d/2))\in\cm^\prq G;\tag h$$
this follows from 1.6(d),(e). This proves (a). At the same time we see that to prove (d) it is enough to 
show that for any $y,d$ as above, (h) is mixed of weight $\le j$. Since $\bbq\ot\L^d\cx[[d]](d/2)$ is pure 
of weight $-d\le0$, to prove the last statement it is enough to show that 
$\c(p_{04!}L^J_{w',y,w,y\i}[[|w|+|w'|+2|y|]])$ is mixed of weight $\le0$.
This follows from the fact that $p_{04!}L^J_{w',y,w,y\i}[[|w|+|w'|+2|y|]]$ is mixed of weight $\le0$ (as in
the proof of 2.7(d)). This proves (d).

We prove (b). It is again enough to show that for any $y,d$ as above
$$\c^j(p_{04!}L^J_{w',y,w,y\i}[[|w|+|w'|+2|y|]]\ot\L^d\cx[[d]](d/2))^j$$
is in $\cm^\prec G$ if $j>\D+4a$. This follows from 1.6(d) since $j>\D+4a$, $d\ge0$ implies $j+d>\D+4a$.

Now (c) follows from (a) by replacing $\boc$ by the two-sided cell containing $w$ (if $w\prec\boc$) or $w'$ 
(if $w'\prec\boc$).

We prove (e). If $j>\D+4a$ this follows from (b). If $j<\D+4a$ we have $gr_{\D+4a}C^j=0$ by (a). This proves 
(e).

We prove (f). If $j<k$ we have $gr_kC^j=0$ by (d). If $j\ge k$ we have $j>\D+4a$ so that $C^j\in\cm^\prec G$ 
by (b). This proves (f).

\proclaim{Proposition 4.5}Let $K,K'\in\cd^\spa_m(G)$.

(a) If $K\in\cd^\prq G$ or $K'\in\cd^\prq G$ then $K*K'\in\cd^\prq G$; if $K\in\cd^\prec G$ or 
$K'\in\cd^\prec G$ then $K*K'\in\cd^\prec G$.

(b) If $K\in\cm^\prq G$, $K'\in\cm^\prq G$ and $j>\r+2a$ then $(K*K')^j\in\cm^\prec G$.
\endproclaim 
It is enough to prove the proposition assuming in addition that $K=A\in CS(G)$, $K'=A'\in CS(G)$. By 1.7(b)
we can find $w\in\boc_A$, $w'\in\boc_{A'}$ such that $(A:R_w^{n_w})\ne0$, $(A':R_{w'}^{n_{w'}})\ne0$. Then 
$A[-n_w]$ is a direct summand of $R_w$ and $A'[-n_{w'}]$ is a direct summand of $R_{w'}$. Hence $A*A'$ is a 
direct summand of $R_w*R_{w'}[2\D+\aa(w)+\aa(w')+|w|+|w'|]$ and $(A*A')^j$ is a direct summand of 
$$(R_w*R_{w'}[2\r+2\nu+|w|+|w'|])^{j+\aa(w)+\aa(w')+2\nu}.$$
 Using 4.4 we deduce that (a) holds and that 
$(A*A')^j\in\cm^\prec G$ provided that $j+\aa(w)+\aa(w')+2\nu>\D+4a$. Hence (b) holds. (To prove (b) we can 
assume, by (a), that $w\in\boc,w'\in\boc$ hence $\aa(w)=\aa(w')=a$.)

\subhead 4.6\endsubhead
For $K,K'\in\cc^\boc_0G$ we set
$$K\unst K'=\un{(K*K')^{\{2a+\r\}}}\in\cc^\boc_0G.$$
We say that $K\unst K'$ is the {\it truncated convolution} of $K,K'$.
Note that 4.1(a) induces for $K,K'\in\cc^\boc_0G$ a canonical isomorphism
$$K\unst K'@>\si>>K'\unst K.\tag a$$
We have also
$$K\unst K'=\op_{j\in\ZZ}\un{gr_{2a+\r}((K*K')^j)((2a+\r)/2)}.\tag b$$
This follows from 4.4(e). 

\proclaim{Proposition 4.7} Let $K,K',K''\in\cc^\boc_0G$. There is a canonical isomorphism
$$(K\unst K')\unst K''@>\si>>K\unst (K'\unst K'').\tag a$$
\endproclaim
We use 1.12 with $\Ph:\cd_m(G)@>>>\cd_m(G)$, $K_1\m K_1*K''$, with $\XX=K*K'$, $(c,c')=(2a+\r,2a+\r)$ (see 
4.5); we deduce that we have canonically 
$$(K\unst K')\unst K''=\un{(K*K'*K'')^{\{4a+2\r\}}}.\tag b$$
Next we use 1.12 with $\Ph:\cd_m(G)@>>>\cd_m(G)$, $K_1\m K*K_1$, with $\XX=K'*K''$, $(c,c')=(2a+\r,2a+\r)$ 
(see 4.5); we deduce that we have canonically 
$$K\unst(K'\unst K'')=\un{(K*K'*K'')^{\{4a+2\r\}}}.\tag c$$
We combine (b),(c); (a) follows.

\subhead 4.8\endsubhead
An argument similar to that in 4.7 shows that the associativity isomorphism provided by 4.7 satisfies the 
pentagon property.

\head 5. Truncated convolution and truncated restriction\endhead
\subhead 5.1\endsubhead
The following proposition asserts a compatibility of truncated restriction with truncated convolution. 

\proclaim{Proposition 5.2} Let $K,K'\in\cc^\boc_0G$. There is a canonical isomorphism (in $\cc^\boc_0\cb^2$):
$$\unz(K')\unb\unz(K)@>\si>>\unz(K\unst K')$$
\endproclaim
The proof is given in 5.6.

\proclaim{Proposition 5.3}. Let $K,K'\in\cc^\boc_0G$. We have canonically
$$\unz(K')\unb\unz(K)=\un{(\z(K')\cir\z(K))^{\{3a+2\r+\nu\}}}.\tag a$$
\endproclaim
We set $L=\z(K)$, $L'=\z(K')$. 
Let ${}^0L\in\cm^\prq_m\cb^2$. Applying 1.12 with $\Ph:\cd^\prq_m\cb^2@>>>\cd^\prq_m\cb^2$, 
${}^1L\m{}^0L\cir{}^1L$, $\XX=L$, $(c,c')=(a-\nu,a+\nu+\r)$, we see that
$$({}^0L\cir L)^j\in\cm^\prec\cb^2\text{ for any ${}^0L\in\cm^\prq\cb^2$ and any }j>2a+\r.\tag b$$
Using 1.12 with $\Ph:\cd_m^\prq\cb^2@>>>\cd_m^\prq\cb^2$, ${}^1L\m {}^1L\cir L$, $\XX=L'$,
$(c,c')=(2a+\r,a+\r+\nu)$ (see (b), 2.8), we deduce that we have canonically 
$$\un{(\un{L'{}^{\{a+\r+\nu\}}}\cir L)^{\{2a+\r\}}}=\un{(L'\cir L)^{\{3a+2\r+\nu\}}}.\tag c$$
Let $L'_0=\un{L'{}^{\{a+\r+\nu\}}}$. Applying 1.12 with $\Ph:\cd_m^\prq\cb^2@>>>\cd_m^\prq\cb^2$,  
${}^1L\m L'_0\cir{}^1L$, $\XX=L$, $(c,c')=(a-\nu,a+\r+\nu)$ (see 3.1, 2.8), we deduce that we have 
canonically 
$$\un{(L'_0\cir\un{L^{\{a+\r+\nu\}}})^{\{a-\nu\}}}=\un{(L'_0\cir L)^{\{2a+\r\}}}.$$
Combining with (c) we obtain
$$\un{(L'_0\cir\un{L^{\{a+\r+\nu\}}})^{\{a-\nu\}}}=\un{(L'\cir L)^{\{3a+2\r+\nu\}}}$$
and (a) follows.

\proclaim{Proposition 5.4} Let $K,K'\in\cc^\boc_0G$. We have canonically
$$\unz(K\unst K')=\un{(\z(K*K'))^{\{3a+\nu+2\r\}}}.\tag a$$
\endproclaim
We set $\ck=K*K'$. Applying 1.12 with $\Ph:\cd_m^\prq G@>>>\cd_m^\prq\cb^2$, $K_1\m\z(K_1)$, $\XX=\ck$,
$(c,c')=(a+\r+\nu,2a+\r)$ (see 2.8, 4.5), we deduce that we have canonically 
$$\un{(\z(\un{\ck^{\{2a+\r\}}}))^{\{a+\r+\nu\}}}=\un{(\z(\ck))^{\{3a+2\r+\nu\}}}$$
and (a) follows.

A version of the following lemma goes back to \cite{\GI}.
\proclaim{Lemma 5.5} Let $K,K'\in\cd_m(G)$. There is a canonical isomorphism in $\cd_m(\cb^2)$:
$$\z(K*K')@>\si>>\z(K')\cir\z(K)\tag b$$
\endproclaim

\subhead 5.6\endsubhead
We prove Proposition 4.2. Let $K,K'\in\cc^\boc_0G$. We have canonically
$$\unz(K')\unb\unz(K)=\un{(\z(K')\cir\z(K))^{\{3a+2\r+\nu\}}}
=\un{(\z(K*K'))^{\{3a+2\r+\nu\}}}=\unz(K\unst K').$$
(These equalities comes from 5.3(a), 5.5, 5.4(a).) Proposition 5.2 follows.

\head 6. Analysis of the composition $\unz\unc$\endhead
\subhead 6.1\endsubhead
Let $e,f,e'$ be integers such that $e\le f\le e'-3$ and let $\e=e'-e+1$; we have $\e\ge4$. We set 
$$\cy=\{((B_e,B_{e+1},\do,B_{e'}),g)\in\cb^\e\T G;gB_fg\i=B_{f+3},gB_{f+1}g\i=B_{f+2}\}.$$
Define $\vt:\cy@>>>\cb^\e$ by $((B_e,B_{e+1},\do,B_{e'}),g)\m(B_e,B_{e+1},\do,B_{e'})$. 
For $i,j$ in $\{e,e+1,\do,e'\}$ let $p_{ij}:\cb^e@>>>\cb^2$ be the projection to the $i,j$ coordinate;
define $h_{ij}:\cy@>>>\cb^2$ by $h_{ij}=p_{ij}\vt$. Now $G^{\e-2}$ acts on $\cy$ by 
$$\align&(g_e,\do,g_f,g_{f+3},\do,g_{e'}):((B_e,B_{e+1},\do,B_{e'}),g)\m\\&
(g_eB_eg_e\i,g_{e+1}B_{e+1}g_{e+1}\i,\do,g_{f-1}B_{f-1}g_{f-1}\i,g_fB_fg_f\i,g_fB_{f+1}g_f\i,\\&
g_{f+3}B_{f+2}g_{f+3}\i,g_{f+3}B_{f+3}g_{f+3}\i,g_{f+4}B_{f+4}g_{f+4}\i,\do,g_{e'}B_{e'}g_{e'}\i),g_{f+3}gg_f\i);\endalign$$
this induces a $G^{\e-2}$-action on $\cb^\e$ so that $\vt$ is $G^{\e-2}$-equivariant.

Let $E=\{e,e+1,\do,e'-1\}-\{f,f+2\}$. Assume that $x_n\in\boc$ are given for $n\in E$. 
Let $P=\ot_{n\in E}p_{n,n+1}^*\LL_{x_n}\in\cd_m\cb^\e$,
$\tP=\ot_{n\in E}h_{n,n+1}^*\LL_{x_n}=\vt^*P\in\cd_m\cy$.
In 6.1-6.7 we will study 
$$h_{ee'!}\tP\in\cd_m\cb^2.$$
Setting $\Xi=\vt_!\bbq\in\cd_m\cb^\e$, we have 
$$h_{ee'!}\tP=p_{ee!}(\Xi\ot P).$$
Clearly $\Xi^j$ is $G^{\e-2}$-equivariant for any $j$. For any $y,y'$ in $W$ we set
$$Z_{y,y'}:=\{(B_e,B_{e+1},\do,B_{e'})\in\cb^\e;(B_f,B_{f+1})\in\co_y,(B_{f+2},B_{f+3})\in\co_{y'}\}.$$
These are the orbits of the $G^{\e-2}$-action on $\cb^\e$. Note that the fibre of $\vt$ over a point of 
$Z_{y,y'}$ is isomorphic to $T\T\kk^{\nu-|y|}$ if $yy'=1$ and is empty if $yy'\ne1$. Thus 

(a) $\Xi|_{Z_{y,y'}}$ is $0$ if $yy'\ne1$
\nl
and for any $y\in W$ we have  

(b) $\ch^h\Xi|_{Z_{y,y\i}}=0$ if $h>2\nu-2|y|+2\r$, $\ch^{2\nu-2|y|+2\r}\Xi|_{Z_{y,y\i}}=\bbq(-\nu+|y|-\r)$.
\nl
The closure of $Z_{y,y'}$ in $\cb^\e$ is denoted by $\bZ_{y,y'}$. We set $k_\e=\e\nu+2\r$.
We have the following result.

\proclaim{Lemma 6.2} (a) We have $\Xi^j=0$ for any $j>k_\e$. Hence, setting $\Xi'=\t_{\le k_\e-1}\Xi$, we 
have a canonical distinguished triangle $(\Xi',\Xi,\Xi^{k_\e}[-k_\e])$. 

(b) If $\x\in Z_{y,y'}$ and $i=2\nu-|y|-|y'|+2\r$, the induced homomorphism 
$\ch^i_\x\Xi@>>>\ch^{i-k_\e}_\x(\Xi^{k_\e})$ is an isomorphism.
\endproclaim
To prove (a) it is enough to show that $\dim\supp\ch^i(\Xi[k_\e])\le-i$ for any $i$. Now $\supp\ch^i\Xi$ is 
a union of $G^{\e-2}$-orbits hence of subvarieties $Z_{y,y'}$ and $\dim Z_{y,y'}=(\e-2)\nu+|y|+|y'|$. Thus it 
is enough to show that if $\ch^i_\x(\Xi[k_\e])\ne0$ with $\x\in Z_{y,y'}$ then $(\e-2)\nu+|y|+|y'|\le-i$. 
From 6.1(a),(b) we see that $y=y'$ and $i+\e\nu+2\r\le2\nu-2|y|+2\r$; the desired result follows.
    
We prove (b). We have an exact sequence
$$\ch^i_\x\Xi'@>>>\ch^i_\x\Xi@>>>\ch^i(\Xi^{k_\e}[-k_\e])@>>>\ch^{i+1}_\x\Xi'.$$
Hence it is enough to show that $\ch^{i'}_\x\Xi'=0$ if $i'\ge i$. Assume that $\ch^{i'}_\x\Xi'\ne0$ for some
$i'\ge i$. Then $Z_{y,y'}\sub\supp\ch^{i'}\Xi'$. We have $(\Xi'[k_\e-1])^h=0$ for all $h>0$ hence 
$\dim\supp\ch^{i''}(\Xi'[k_\e-1])\le-i''$ for any $i''$. Taking $i''=i'-k_\e+1$ we deduce that 
$\dim Z_{y,y'}\le-i'+k_\e-1$ hence $i'\le 2\nu-|y|-|y'|+2\r-1=i-1$. This contradicts $i'\ge i$ and proves 
(b).

\subhead 6.3\endsubhead
For any $y,y'$ in $W$ let $\fT_{y,y'}$ be the intersection cohomology complex of $\bZ_{y,y'}$ extended by 
$0$ on $\cb^\e-\bZ_{y,y'}$, to which $[[(\e-2)\nu+|y|+|y'|]]$ is applied. Note that
$$\fT_{y,y'}=p_{f,f+1}^*\LL_y\ot p_{f+2,f+3}^*\LL_{y'}[[(\e-4)\nu]].\tag a$$
We have the following result.

\proclaim{Lemma 6.4} We have canonically $gr_0(\Xi^{k_\e}(k_\e/2))=\op_{y\in W}\fT_{y,y\i}$.
\endproclaim
Since $gr_0(\Xi^{k_\e}(k_\e/2))$ is a semisimple $G^{\e-2}$-equivariant perverse sheaf of pure weight $0$, we 
have canonically $gr_0(\Xi^{k_\e}(k_\e/2))=\op_{y,y'\in W}V_{y,y'}\ot\fT_{y,y'}$ where $V_{y,y'}$ are mixed 
$\bbq$-vector spaces of pure weight $0$. Using the definition or by \cite{\BBD, 5.1.14}, $\Xi$ is mixed of 
weight $\le0$ hence $\Xi^{k_\e}(k_\e/2)$ is mixed of weight $\le0$. Hence we have an exact sequence in 
$\cm_m\cb^\e$
$$0@>>>\cw^{-1}(\Xi^{k_\e}(k_\e/2))@>>>\Xi^{k_\e}(k_\e/2)@>>>gr_0(\Xi^{k_\e}(k_\e/2))@>>>0\tag a$$
that is  
$$0@>>>\cw^{-1}(\Xi^{k_\e}(k_\e/2))@>>>\Xi^{k_\e}(k_\e/2)@>>>\op_{y,y'\in W}V_{y,y'}\ot\fT_{y,y'}@>>>0.$$
Hence for any $\ty,\ty'$ in $W$ and any $\FF_q$-rational point $\x\in Z_{\ty,\ty'}$ we have an exact 
sequence of stalks of cohomology sheaves 
$$\align&\ch^h_\x\cw^{-1}(\Xi^{k_\e}(k_\e/2))@>\a>>\ch^h_\x\Xi^{k_\e}(k_\e/2)@>>>\\&
\op_{y,y'\in W}V_{y,y'}\ot\ch^h_\x\fT_{y,y'}@>>>\ch^{h+1}_\x\cw^{-1}(\Xi^{k_\e}(k_\e/2));\tag b\endalign$$
here we take $h=-(\e-2)\nu-|\ty|-|\ty'|$.
Now the vector spaces in (b) are mixed and the maps respect the mixed structures. From 6.2(b) and 6.1 we see
that $\ch^h_\x(\Xi^{k_\e}(k_\e/2))=\ch^{h+k_\e}_\x\Xi(k_\e/2)=V_0(-h/2)$ where $V_0$ is $0$ if $\ty\ty'\ne1$
and is $\bbq$ if $\ty\ty'=1$. In particular $\ch^i_\x(\Xi^{k_\e}(k_\e/2))$ is pure of weight $h$. 
On the other hand the mixed vector space $\ch^h_\x\cw^{-1}(\Xi^{k_\e}(k_\e/2))$ has weight $\le h-1$. Hence 
the map $\a$ in (b) must be zero.

Assume that $\ch^h_\x\fT_{y,y'}\ne0$. Then $Z_{\ty,\ty'}$ is contained in the support of $\ch^h\fT_{y,y'}$ 
which has dimension $\le-h$ (resp. $<-h$ if $(y,y')\ne(\ty,\ty')$); hence $-h=\dim Z_{ty,\ty'}$ is
$\le-h$ (resp. $<-h$); 
we see that we must have $(y,y')=(\ty,\ty')$. Note also that $\ch^h_\x\fT_{y,y'}=\bbq(-h/2)$.

Assume that $\ch^{h+1}_\x\cw^{-1}(\Xi^{k_\e}(k_\e/2))\ne0$; then $Z_{\ty,\ty'}$ is contained in the support 
of $\ch^{h+1}_\x\cw^{-1}(\Xi^{k_\e})$ which has dimension $\le-h-1$; hence $-h=\dim Z_{ty,\ty'}\le-h-1$, a 
contradiction. We see that (b) becomes an isomorphism
$$V_0(-h/2)@>\si>>V_{\ty,\ty'}(-h/2).$$
It follows that we have canonically $V_{\ty,\ty'}=V_0$. The lemma is proved.

\subhead 6.5\endsubhead
Let $y,\ty\in W$. Using the definitions and 1.2(a) we have 
$$\align&p_{ee'!}(\fT_{y,\ty}\ot P[[(6-2\e)\nu]])\\&
=L_{x_1}^\sh\cir\do\cir L_{x_{f-1}}^\sh\cir L_y^\sh\cir L_{x_{f+1}}^\sh\cir L_{\ty}^\sh\cir
L_{x_{f+3}}^\sh\cir\do\cir L_{x_{e'}}^\sh[[\nu+|y|+|\ty|+\sum_{n\in E}|x_n|]].\tag a\endalign$$

\proclaim{Lemma 6.6}The map $\Xi@>>>\Xi^{k_\e}[-k_\e]$ (coming from $(\Xi',\Xi,\Xi^{k_\e}[-k_\e])$ in 
6.2(a)) induces a morphism 
$$(p_{ee'!}(\Xi\ot P))^{(\e-2)a+(6-\e)\nu+2\r}@>>>(p_{ee'!}(\Xi^{k_\e}\ot P))^{(\e-2)a+(6-\e)\nu+2\r-k_\e}$$
whose kernel and cokernel are in $\cm^\prec_m\cb^2$.
\endproclaim
It is enough to prove that
$$(p_{ee'!}(\Xi'\ot P))^j\in\cm^\prec_m\cb^2\text{ for any }j\ge (\e-2)a+(6-\e)\nu+2\r.$$
We have $\Xi'\Bpq\{(\Xi')^h[-h];h\le k_\e-1\}$ hence
$$p_{ee'!}(\Xi'\ot P)\Bpq\{p_{ee'!}(\Xi')^h\ot P)[-h];h\le k_\e-1\}$$ 
so that it is enough to show that 
$$(p_{ee'!}(\Xi')^h\ot P)[-h])^j\in\cm^\prec_m\cb^2$$
for any $j\ge (\e-2)a+(6-\e)\nu+2\r$ and any $h\le k_\e-1$.
Now $(\Xi')^h$ is $G^{\e-2}$-equivariant hence its composition factors are of the form $\fT_{y,y'}$ with 
$y,y'$ in $W$; hence it is enough to show that for any $y,y'$ in $W$ we have
$$(p_{ee'!}(\fT_{y,y'}\ot P)[-h])^j\in\cm^\prec_m\cb^2$$
for any $j\ge (\e-2)a+(6-\e)\nu+2\r$ and any $h\le k_\e-1$ or equivalently (see 6.5(a),(b)) that
$$\align&(L_{x_1}^\sh\cir\do\cir L_{x_{f-1}}^\sh\cir L_y^\sh\cir L_{x_{f+1}}^\sh\cir L_{\ty}^\sh\cir
L_{x_{f+3}}^\sh\cir\do\cir L_{x_{e'}}^\sh\\&[[(2\e-5)\nu+|y|+|y'|+\sum_n|x_n|]])^{j-h}
\in\cm^\prec_m\cb^2\endalign$$
for any $j\ge(\e-2)a+(6-\e)\nu+2\r$ and any $h\le f_\e-1$. Using 2.2(a) it is enough to show that 
$j-h+(2\e-5)\nu>\nu+(\e-2)a$. We have 
$$j-h+(2\e-5)\nu\ge (\e-2)a+(6-\e)\nu+2\r-\e\nu-2\r+1+(2\e-5)\nu=(\e-2)a+\nu+1$$
and the lemma is proved.

\proclaim{Lemma 6.7} We have canonically
$$\un{(h_{ee'!}\tP)^{\{(\e-2)a+(6-\e)\nu+2\r\}}}=\op_{y\in\boc}Q_y$$
where 
$$\align&Q_y=\un{(p_{ee'!}(\fT_{y,y\i}\ot P))^{\{(\e-2)a+(6-2\e)\nu\}}}\\&
=\LL_{x_1}\unb\do\unb\LL_{x_{f-1}}\unb\LL_y\unb\LL_{x_{f+1}}\unb\LL_{y\i}\unb\LL_{x_{f+3}}\unb\do
\unb\LL_{x_{e'}}.\endalign$$
\endproclaim
From the exact sequence 6.4(a) we deduce a distinguished triangle in $\cd_m\cb^2$:
$$(p_{ee'!}(\cw^{-1}(\Xi^{k_\e}(k_\e/2))\ot P,p_{ee'!}(\Xi^{k_\e}(k_\e/2)\ot P),  
p_{ee'!}(gr_0(\Xi^{k_\e}(k_\e/2))\ot P).$$
This induces an exact sequence in $\cm_m\cb^2$:
$$\align&(p_{ee'!}(\cw^{-1}(\Xi^{k_\e}(k_\e/2))\ot P))^{(\e-2)a+(6-2\e)\nu}\\&@>>>
(p_{ee'!}(\Xi^{k_\e}(k_\e/2)\ot P))^{(\e-2)a+(6-2\e)\nu}\\&
@>>>(p_{ee'!}(gr_0(\Xi^{k_\e}(k_\e/2))\ot P))^{(\e-2)a+(6-2\e)\nu}\\&@>>>
(p_{ee'!}(\cw^{-1}(\Xi^{k_\e}(k_\e/2))\ot P))^{(\e-2)a+(6-2\e)\nu+1}.\tag a\endalign$$
We show that 
$$(p_{ee'!}(\cw^{-1}(\Xi^{k_\e}(k_\e/2))\ot P))^{(\e-2)a+(6-2\e)\nu+1}\in\cm^\prec_m\cb^2.\tag b$$
We argue as in the proof of 6.6. Now $\cw^{-1}(\Xi^{k_\e}(k_\e/2))$ is $G^{\e-2}$-equivariant hence its 
composition factors are of the form $\fT_{y,y'}$ with $y,y'$ in $W$; hence it is enough to show that for any 
$y,y'$ in $W$ we have
$$(p_{ee'!}\fT_{y,y'}\ot P))^{(\e-2)a+(6-2\e)\nu+1}\in\cm^\prec_m\cb^2$$
or equivalently (see 6.5(a),(b)) that
$$\align&(L_{x_1}^\sh\cir\do\cir L_{x_{f-1}}^\sh\cir L_y^\sh\cir L_{x_{f+1}}^\sh\cir L_{\ty}^\sh\cir 
L_{x_{f+3}}^\sh\cir\do\cir L_{x_{e'}}^\sh\\&[[(2\e-5)\nu+|y|+|y'|+\sum_n|x_n|]])^{(\e-2)a+(6-2\e)\nu+1}
\in\cm^\prec_m\cb^2.\endalign$$
Using 2.2(a) it remains to note that $(\e-2)a+(6-2\e)\nu+1+(2\e-5)\nu>\nu+(\e-2)a$.

Next we show that 
$$gr_{(\e-2)a+(6-2\e)\nu}(p_{ee'!}(\cw^{-1}(\Xi^{k_\e}(k_\e/2))\ot P))^{(\e-2)a+(6-2\e)\nu}=0.\tag c$$
Indeed, $\cw^{-1}(\Xi^{k_\e}(k_\e))$ has weight $\le-1$, $P$ has weight $0$ hence
$\cw^{-1}(\Xi^{k_\e}(k_\e/2))\ot P$ has weight $\le-1$ and
$p_{ee'!}(\cw^{-1}(\Xi^{k_\e}(k_\e/2))\ot P)$ has weight $\le-1$ so that
$(p_{ee'!}(\cw^{-1}(\Xi^{k_\e}(k_\e/2))\ot P))^{(\e-2)a+(6-2\e)\nu}$ has weight $\le (\e-2)a+(6-2\e)\nu-1$ 
and (c) follows.

Using (b),(c) we see that (a) induces a morphism
$$\align&gr_{(\e-2)a+(6-2\e)\nu}(p_{ee'!}(\Xi^{k_\e}(k_\e/2)\ot P))^{(\e-2)a+(6-2\e)\nu}\\&@>>>
gr_{(\e-2)a+(6-2\e)\nu}(p_{ee'!}(gr_0(\Xi^{k_\e}(k_\e/2))\ot P))^{(\e-2)a+(6-2\e)\nu}\endalign$$
which has kernel $0$ and cokernel in $\cm^\prec_m\cb^2$. Hence we have an induced isomorphism
$$\align&\un{gr_{(\e-2)a+(6-2\e)\nu}(p_{ee'!}}
\un{(\Xi^{k_\e}(k_\e/2)\ot P))^{(\e-2)a+(6-2\e)\nu}}((\e-2)a+(6-2\e)\nu)/2)@>\si>>  \\&
\un{gr_{(\e-2)a+(6-2\e)\nu}}\un{(p_{ee'!}(gr_0(\Xi^{k_\e}(k_\e/2))\ot P))^{(\e-2)a+(6-2\e)\nu}}
((\e-2)a+(6-2\e)\nu)/2).\tag d\endalign$$
The left hand side of (d) can be identified (by 6.6) with   
$$\align&\un{gr_{(\e-2)a+(6-2\e)\nu}(p_{ee'!}(\Xi(k_\e)\ot P))^{(\e-2)a+(6-\e)\nu+2\r}}
((\e-2)a+(6-2\e)\nu)/2)\\&=
\un{gr_{(\e-2)a+(6-\e)\nu+2\r}(p_{ee'!}(\Xi\ot P))^{(\e-2)a+(6-\e)\nu+2\r}}(k_\e/2)((\e-2)a+(6-2\e)\nu)/2)\\&
=\un{(p_{ee'!}(\Xi\ot P))^{\{(\e-2)a+(6-\e)\nu+2\r\}}};\endalign$$
the right hand side of (d) can be identified (by 6.4 and 6.5(a),(b)) with $\op_{y\in W}Q_y$ where
$$\align&Q_y=\un{gr_{(\e-2)a+(6-2\e)\nu}(p_{ee'!}(\fT_{y,y\i}\ot P))^{(\e-2)a+(6-2\e)\nu}}
((\e-2)a+(6-2\e)\nu)/2)\\&
=\un{gr_{(\e-2)a+(6-2\e)\nu}(L_{x_1}^\sh\cir\do\cir 
L_{x_{f-1}}^\sh\cir L_y^\sh\cir L_{x_{f+1}}^\sh\cir L_{\ty}^\sh\cir L_{x_{f+3}}^\sh\cir\do\cir
 L_{x_{e'}}^\sh}
\\&\un{[[(2\e-5)\nu+|y|+|y'|+\sum_n|x_n|]])^{(\e-2)a+(6-2\e)\nu}}((\e-2)a+(6-2\e)\nu)/2).\endalign$$
Thus,
$$\align&Q_y=
\un{gr_{(\e-2)a+(6-2\e)\nu}(\LL_{x_1}\cir\do\cir\LL_{x_{f-1}}\cir\LL_y\cir\LL_{x_{f+1}}\cir\LL_{\ty}\cir}\\&
\un{\LL_{x_{f+3}}\cir\do\cir\LL_{x_{e'}})^{(\e-2)a-(\e-2)\nu}}((\e-4)\nu/2)((\e-2)a+(6-2\e)\nu)/2)\\&
=\un{gr_{(\e-2)a+(2-\e)\nu}(\LL_{x_1}\cir\do\cir\LL_{x_{f-1}}\cir\LL_y\cir\LL_{x_{f+1}}\cir\LL_{\ty}\cir}\\&
\un{\LL_{x_{f+3}}\cir\do\cir\LL_{x_{e'}})^{(\e-2)a-(\e-2)\nu}}((\e-2)a+(\e-2)\nu)/2)\\&
=\un{(\LL_{x_1}\cir\do\cir\LL_{x_{f-1}}\cir\LL_y\cir\LL_{x_{f+1}}\cir\LL_{\ty}\cir}
\un{\LL_{x_{f+3}}\cir\do\cir\LL_{x_{e'}})^{\{(\e-2)a-(\e-2)\nu\}}}.\endalign$$
Thus we have canonically
$$\un{(p_{ee'!}(\Xi\ot P))^{\{(\e-2)a+(6-\e)\nu+2\r\}}}=\op_{y\in W}Q_y$$
where
$$\align&Q_y=\un{(p_{ee'!}(\fT_{y,y\i}\ot P))^{\{(\e-2)a+(6-2\e)\nu\}}}\\&=
\un{(\LL_{x_1}\cir\do\cir\LL_{x_{f-1}}\cir\LL_y\cir\LL_{x_{f+1}}\cir\LL_{\ty}\cir} 
\un{\LL_{x_{f+3}}\cir\do\cir\LL_{x_{e'}})^{\{(\e-2)a-(\e-2)\nu\}}}.\endalign$$
The expression following the last $=$ sign is $0$ if $y\n\boc$ (see 2.3) and is
$$\LL_{x_1}\unb\do\unb\LL_{x_{f-1}}\unb\LL_y\unb\LL_{x_{f+1}}\unb\LL_{\ty}\unb\LL_{x_{f+3}}\unb
\do\unb\LL_{x_{e'}}$$
if $y\in\boc$ (see 3.2).
The lemma is proved.

\proclaim{Theorem 6.8} Let $x\in\boc$. We have canonically
$$\unz(\unc(\LL_x))=\op_{y\in\boc}\LL_y\unb\LL_x\unb\LL_{y\i}.\tag a$$
\endproclaim
In 6.1 we take $e=f=1,e'=4$ hence $\e=4$. In this case we have
$$\cy=\{((B_1,B_2,B_3,B_4),g)\in\cb^4\T G;gB_1g\i=B_4,gB_2g\i=B_3\}.$$
Let $x\in\boc$. From Lemma 6.7 we have canonically
$$(h_{14!}h_{23}^*\LL_x))^{\{2a+2\nu+2\r\}}=\op_{y\in\boc}\LL_y\unb\LL_x\unb\LL_{y\i}.\tag b$$
By the proof of 2.6 we have 
$$\z(\c(\LL_x))=h_{14!}h_{23}^*\LL_x.$$
Hence, using 2.11(b), we have 
$$\unz(\unc(\LL_x))=\un{(h_{14!}h_{23}^*\LL_x)^{\{2a+2\nu+2\r\}}}.$$
Substituting this into (b) we obtain (a).

\subhead 6.9\endsubhead
Using 2.4 we see that 6.8(a) implies
$$\unz\unc\LL_x\cong\op_{z\in\boc}(\LL_z)^{\op\ps_x(z)}\tag a$$
in $\cc^\boc\cb^2$ where $\ps_x(z)\in\NN$ are given by the following equation in $\JJ^\boc$:
$$\sum_{y\in\boc}t_yt_xt_{y\i}=\sum_{z\in\boc}\ps_x(z)t_z.$$

\head 7. Analysis of the composition $\unz\unc$ (continued)\endhead
\subhead 7.1\endsubhead
Let $\fZ=\{((\b_1,\b_2,\b_3,\b_4),g)\in\cb^4\T G;g\b_2g\i=\b_3\}$.
Define $d,d':\fZ@>>>\cb^2$ by $d(\b_1,\b_2,\b_3,\b_4),g)=(\b_1,g\i\b_4g)$,
$d'(\b_1,\b_2,\b_3,\b_4),g)=(g\b_1g\i,\b_4)$. Let $u\in\boc$. We set 
$\tLL_u=d^*\LL_u=d'{}^*\LL_u\in\cd_m(\fZ)$; the last equality follows from the $G$-equivariance of $\LL_u$.
Define $\bvt:\fZ@>>>\cb^4$ by $((\b_1,\b_2,\b_3,\b_4),g)\m(\b_1,\b_2,\b_3,\b_4)$. Now $G^2$ acts on $\fZ$ by
$$(g_1,g_2):((\b_1,\b_2,b_3,\b_4),g)\m((g_1\b_1g_1\i,g_1\b_2g_1\i,g_2b_3g_2\i,g_2\b_4g_2\i),g_2gg_1\i);$$
this induces a $G^2$-action on $\cb^4$ so that $\bar\vt$ is $G^2$-equivariant.
Note also that $G^2$ acts on $\cb^2$ by $(g_1,g_2):(B,B')\m(g_1Bg_1\i,g_1B'g_1\i)$ and that $d,d'$ are
$G^2$-equivariant. It follows that a shift of $\tLL_u$ is $G\T G$-equivariant perverse sheaf and 
$(\bvt_!\tLL_u)^j$ is $G^2$-equivariant for any $j$. 

For $i,j$ in $\{1,2,3,4\}$ let $\bp_{ij}:\cb^4@>>>\cb^2$ be the projection to the $i,j$ coordinates.

For any $y,z$ in $W$ we set
$$\fZ_{y,z}=\{(\b_1,\b_2,\b_3,\b_4)\in\cb^4;(\b_1,\b_2)\in\co_y,(\b_3,\b_4)\in\co_z\}.$$
These are the orbits of the $G^2$-action on $\cb^4$. Let $\TT_{y,z}$ be the intersection cohomology complex 
of the closure $\bar\fZ_{y,z}$ of $\fZ_{y,z}$ extended by $0$ on $\cb^4-\bar\fZ_{y,z}$, to which 
$[[2\nu+|y|+|z|]]$ has been applied. We have $\TT_{y,z}=\bp_{12}^*\LL_y\ot\bp_{34}^*\LL_z$.

We denote by ${}'\cm^\prq\cb^4$ (resp. ${}''\cm^\prq\cb^4$) the category of perverse sheaves on $\cb^4$ 
whose composition factors are all of the form $\TT_{y,z}$ with $y\prq\boc$, $z\in W$ (resp. $y\in W$, 
$z\prq\boc$). We denote by ${}'\cm^\prec\cb^4$ (resp. ${}''\cm^\prec\cb^4$) the category of perverse sheaves
on $\cb^4$ whose composition factors are all of the form $\TT_{y,z}$ with $y\prec\boc$, $z\in W$ (resp. 
$y\in W$, $z\prec\boc$). 
Let $\cm^\prq\cb^4$ (resp. $\cm^\prec\cb^4$) be the category of perverse sheaves on $\cb^4$
whose composition factors are all of the form $\TT_{y,z}$ with $y\prq\boc$, $z\prq\boc$ 
(resp. $y\prq\boc$, $z\prec\boc$ or $y\prec\boc$, $z\prq\boc$).
Let $\cd_m^\prq\cb^4$ (resp. $\cd_m^\prec\cb^4$) be the category consisting of all $K\in\cd_m\cb^4$
such that for any $j\in\ZZ$, $K^j$ belongs to $\cm^\prq\cb^4$ (resp. $\cm^\prec\cb^4$).

Let $\cc^\prq\cb^4$ be the subcategory of $\cm^\prq\cb^4$ consisting of semisimple objects; let 
$\cc^\prq_0\cb^4$ be the subcategory of $\cm_m\cb^4$ consisting of those $K\in\cm_m\cb^4$ such that $K$ is 
pure of weight $0$ and such that as an object of $\cm\cb^4$, $K$ belongs to $\cc^\prq\cb^4$. Let 
$\cc^\boc\cb^4$ be the the subcategory of $\cm^\prq\cb^4$ consisting of objects which are direct sums of 
objects of the form $\TT_{y,z}$ with $y\in\boc$, $z\in\boc$. Let $\cc^\boc_0\cb^4$ be the subcategory of 
$\cc^\prq_0\cb^4$ consisting of those $K\in\cc^\prq_0\cb^4$ such that as an object of $\cc^\prq\cb^4$, $K$ 
belongs to $\cc^\boc\cb^4$. For $K\in\cc^\prq_0\cb^4$, let $\un{K}$ be the largest subobject of $K$ such 
that as an object of $\cc^\prq\cb^4$, we have $\un{K}\in\cc^\boc\cb^4$.   

We set $\a=a+3\nu+2\r$. We have canonically 
$$gr_0((\bvt_!\tLL_u)^\a(\a/2))=\op_{y,z\in W}U_{y,z}\ot\TT_{y,z}\tag a$$
where $U_{y,z}$ are well defined mixed $\bbq$ vector spaces of pure weight $0$. 

\proclaim{Lemma 7.2} (a) For any $j\in\ZZ$ we have $(\bvt_!\tLL_u)^j\in\cm^\prq\cb^4$.

(b) If $j>\a$ then $(\bvt_!\tLL_u)^j\in{}'\cm^\prec\cb^4\cap{}''\cm^\prec\cb^4$.

(c) If $y,z\in\boc$, we have canonically $U_{y,z}=\Hom_{\cc^\boc\cb^2}(\LL_y,\LL_u\unb\LL_{z\i})$.

(d) If $y,z\in\boc$, we have canonically $U_{y,z}=\Hom_{\cc^\boc\cb^2}(\LL_z,\LL_{y\i}\unb\LL_u)$.
\endproclaim
The proof of (a) and (b) is given in 7.3 and 7.4. The proof of (c) is given in 7.5. The proof of (d) is 
given in 7.6.

\subhead 7.3\endsubhead
In this subsection we show that

(a) For any $j\in\ZZ$ we have $(\bvt_!\tLL_u)^j\in{}'\cm^\prq\cb^4$.

(b) If $j>\a$ then $(\bvt_!\tLL_u)^j\in{}'\cm^\prec\cb^4$.
\nl
In the setup of 6.1 (with $e=0,f=1,e'=4$ hence $\e=5$) we identify $\cy$, $\fZ$ via the isomorphism
$${}'\fc:\cy@>\sim>>\fZ,\qua((B_0,B_1,B_2,B_3,B_4),g)\m((B_0,B_2,B_3,B_4),g).$$
Then $\bvt$ becomes the composition $\cy@>\vt>>\cb^5@>\th>>\cb^4$ where $\th$ is 
$(B_0,B_1,B_2,B_3,B_4)\m(B_0,B_2,B_3,B_4)$; $\bvt_!\tLL_u$ becomes $\th_!((p_{01}^*\LL_u)\ot\Xi)$.

We have 
$$\th_!((p_{01}^*\LL_u)\ot\Xi)\Bpq\{\th_!((p_{01}^*\LL_u)\ot\Xi^h[-h]);h\le k_5\}$$
where the inequality $h\le k_5$ comes from the fact that $\Xi^h=0$ if $h>k_5$, see 6.2(a). (Recall that 
$k_5=5\nu+2\r$.) Hence it is enough to show:

(c) For any $j,h\in\ZZ$,  we have $(\th_!((p_{01}^*\LL_u)\ot\Xi^h[-h]))^j\in{}'\cm^\prq\cb^4$.

(d) For any $j,h\in\ZZ$ such that $j>\a,h\le k_5$, we have 
$(\th_!((p_{01}^*\LL_u)\ot\Xi^h[-h]))^j\in{}'\cm^\prec\cb^4$.
\nl
Note that in (d) we have $j-h>a-2\nu$. Since $\Xi^h$ is $G^3$-equivariant, its composition factors are of 
the form $\fT_{y,y'}$ with $y,y'\in W$. Hence it is enough to show for any $y,y'\in W$:

(e) For any $j'\in\ZZ$, we have $(\th_!((p_{01}^*\LL_u)\ot\fT_{y,y'}))^{j'}\in{}'\cm^\prq\cb^4$.

(f) For any $j'\in\ZZ$ such that $j'>a-2\nu$, we have 
$(\th_!((p_{01}^*\LL_u)\ot\fT_{y,y'}))^{j'}\in{}'\cm^\prec\cb^4$.
\nl
From the definitions we have 
$$\th_!((p_{01}^*\LL_u)\ot\fT_{y,y'})=\bp_{12}^*(\LL_u\cir\LL_y)[[\nu]]\ot\bp_{34}^*\LL_{y'}.$$
This can be viewed as $(\LL_u\cir\LL_y[[\nu]])\bxt\LL_{y'}\in\cd_m(\cb^2\T\cb^2)$.
Since $\LL_{y'}$ is a perverse sheaf on the second copy of $\cb^2$, we have
$$((\LL_u\cir\LL_y[[\nu]])\bxt\LL_{y'})^{j'}=(\LL_u\cir\LL_y)^{j'+\nu}\bxt\LL_{y'}(\nu/2).$$
It remains to observe that $(\LL_u\cir\LL_y)^{j'+\nu}$ is in $\cm^\prq\cb^2$ for any $j'$ and is in
$\cm^\prec\cb^2$ if $j'+\nu>a-\nu$ (by 3.1). This proves (a),(b).

\subhead 7.4\endsubhead
In this subsection we show that

(a) For any $j\in\ZZ$ we have $(\bvt_!\tLL_u)^j\in{}''\cm^\prq\cb^4$.

(b) If $j>\a$ then $(\bvt_!\tLL_u)^j\in{}''\cm^\prec\cb^4$.
\nl
The arguments are almost a copy of those in 7.3.
In the setup of 6.1 (with $e=1,f=1,e'=5$ hence $\e=5$) we identify $\cy$, $\fZ$ via the isomorphism
$${}''\fc:\cy@>\sim>>\fZ,\qua((B_1,B_2,B_3,B_4,B_5),g)\m((B_1,B_2,B_3,B_5),g).$$
Then $\bvt$ becomes the composition $\cy@>\vt>>\cb^5@>\th>>\cb^4$ where $\th$ is 
$(B_1,B_2,B_3,B_4,B_5)\m(B_1,B_2,B_3,B_5)$; $\bvt_!\tLL_u$ becomes $\th_!((p_{45}^*\LL_u)\ot\Xi)$.
We have $\th_!((p_{45}^*\LL_u)\ot\Xi)\Bpq\{\th_!((p_{45}^*\LL_u)\ot\Xi^h[-h]);h\le k_5\}$
where the inequality $h\le k_5$ comes from the fact that $\Xi^h=0$ if $h>k_5$, see 6.2(a). Hence it is
 enough to show:

(c) For any $j,h\in\ZZ$,  we have $(\th_!((p_{45}^*\LL_u)\ot\Xi^h[-h]))^j\in{}''\cm^\prq\cb^4$.

(d) For any $j,h\in\ZZ$ such that $j>\a,h\le k_5$, we have 
$(\th_!((p_{45}^*\LL_u)\ot\Xi^h[-h]))^j\in{}''\cm^\prec\cb^4$.
\nl
Note that in (d) we have $j-h>a-2\nu$. Since $\Xi^h$ is $G^3$-equivariant, its composition factors are of 
the form $\fT_{y,y'}$ with $y,y'\in W$. Hence it is enough to show for any $y,y'\in W$:

(e) For any $j'\in\ZZ$, we have $(\th_!((p_{45}^*\LL_u)\ot\fT_{y,y'}))^{j'}\in{}''\cm^\prq\cb^4$.

(f) For any $j'\in\ZZ$ such that $j'>a-2\nu$, we have 
$(\th_!((p_{45}^*\LL_u)\ot\fT_{y,y'}))^{j'}\in{}''\cm^\prec\cb^4$.
\nl
From the definitions we have 
$$\th_!((p_{45}^*\LL_u)\ot\fT_{y,y'})=\bp_{34}^*(\LL_{y'}\cir\LL_u)[[\nu]]\ot\bp_{12}^*\LL_y.$$
This can be viewed as $\LL_y\bxt(\LL_{y'}\cir\LL_u[[\nu]])\in\cd_m(\cb^2\T\cb^2)$.
Since $\LL_y$ is a perverse sheaf on the first copy of $\cb^2$, we have
$$(\LL_y\bxt(\LL_{y'}\cir\LL_u[[\nu]]))^{j'}=\LL_y(\nu/2)\bxt(\LL_{y'}\cir\LL_u)^{j'+\nu}.$$
It remains to observe that $(\LL_{y'}\cir\LL_u)^{j'+\nu}$ is in $\cm^\prq\cb^2$ for any $j'$ and is in
$\cm^\prec\cb^2$ if $j'+\nu>a-\nu$ (by 3.1). This proves (a),(b).

Combining (a),(b) with 7.3(a),(b) we see that 7.2(a),(b) hold.

\subhead 7.5\endsubhead
We prove 7.2(c) using the isomorphism ${}'\fc:\cy@>\sim>>\fZ$ in 7.3. (We assume again that we are in the 
setup of 6.1 with $e=0,f=1,e'=4$ hence $\e=5$.) As in 7.3, we have 
$\bvt_!\tLL_u=\th_!((p_{01}^*\LL_u)\ot\Xi)$. Here $\th:\cb^5@>>>\cb^4$ is as in 7.3.

From the exact triangle $(\Xi',\Xi,\Xi^{k_5}[-k_5])$ in 6.2(a) we get an exact triangle
$$(\th_!((p_{01}^*\LL_u)\ot\Xi'), \th_!((p_{01}^*\LL_u)\ot\Xi),\th_!((p_{01}^*\LL_u)\ot\Xi^{k_5}[k_5]))$$
hence an exact sequence
$$\align&(\th_!((p_{01}^*\LL_u)\ot\Xi'))^j@>>>(\th_!((p_{01}^*\LL_u)\ot\Xi))^j\\&@>>>
(\th_!((p_{01}^*\LL_u)\ot\Xi^{k_5}]k_5]))^j@>>>(\th_!((p_{01}^*\LL_u)\ot\Xi'))^{j+1}.\tag a\endalign$$ 
Replacing $\Xi$ by $\Xi'$ in the proof of 7.3(b) given in 7.3 and using that $(\Xi')^h=0$ if $h\ge k_5$ we 
see that
$$(\th_!((p_{01}^*\LL_u)\ot\Xi'))^j\in\cm^\prec\cb^4\text{ for }j\ge\a.$$
Hence the exact sequence (a) implies that
$$(\th_!((p_{01}^*\LL_u)\ot\Xi))^\a@>>>(\th_!((p_{01}^*\LL_u)\ot\Xi^{k_5}[k_5]))^\a$$
has kernel and cokernel in $\cm^\prec\cb^4$. 
This induces a homomorphism
$$\align&gr_0(\th_!((p_{01}^*\LL_u)\ot\Xi))^\a(\a/2))
@>>>gr_0(\th_!((p_{01}^*\LL_u)\ot\Xi^{k_5}[k_5]))^\a(\a/2))\\&
=gr_0(\th_!((p_{01}^*\LL_u)\ot\Xi^{k_5}))^{\a-k_5}(\a/2))\tag b\endalign$$
which has kernel and cokernel in $\cm^\prec\cb^4.$ 

From the exact sequence 6.4(a) we get a distinguished triangle
$$\align&(\th_!((p_{01}^*\LL_u)\ot\cw^{-1}(\Xi^{k_5}(k_5/2)))[-k_5],
\th_!((p_{01}^*\LL_u)\ot\Xi^{k_5}(k_5/2))[-k_5],\\&
\th_!((p_{01}^*\LL_u)\ot gr_0(\Xi^{k_5}(k_5/2)))[-k_5]).\endalign$$
Hence we have an exact sequence  
$$\align&(\th_!((p_{01}^*\LL_u)\ot\cw^{-1}(\Xi^{k_5}(k_5/2)))[-k_5])^\a\\&@>>>
(\th_!((p_{01}^*\LL_u)\ot\Xi^{k_5}(k_5/2))[-k_5])^\a\\&
@>>>(\th_!((p_{01}^*\LL_u)\ot gr_0(\Xi^{k_5}(k_5/2)))[-k_5])^\a\\&@>>>
(\th_!((p_{01}^*\LL_u)\ot\cw^{-1}(\Xi^{k_5}(k_5/2)))[-k_5])^{\a+1}.\tag c\endalign$$
Replacing $\Xi$ by $\cw^{-1}(\Xi^{k_5}(k_5/2))[-k_5/2]$ in the proof of 7.3(b) given in 7.3 and using that
$(\cw^{-1}(\Xi^{k_5}(k_5/2))[-k_5/2])^h=0$ if $h>k_5$ we see that
$$(\th_!((p_{01}^*\LL_u)\ot\cw^{-1}(\Xi^{k_5}(k_5/2)))[-k_5])^{\a+1}\in\cm^\prec\cb^4.\tag d$$
Note that
$$gr_{\a-k_5}(\th_!((p_{01}^*\LL_u)\ot\cw^{-1}(\Xi^{k_5}(k_5/2))))^{\a-k_5}=0.\tag e$$
This follows from the fact that $\cw^{-1}(\Xi^{k_5}(k_5/2))$ has weight $\le-1$ hence \lb
$\th_!((p_{01}^*\LL_u)\ot\cw^{-1}(\Xi^{k_5}(k_5/2)))$ has weight $\le-1$ and \lb
$(\th_!((p_{01}^*\LL_u)\ot\cw^{-1}(\Xi^{k_5}(k_5/2))))^{\a-k_5}$ has weight $\le\a-k_5-1$.

Using (d),(e), we see that (c) induces a morphism
$$gr_{\a-k_5}(\th_!((p_{01}^*\LL_u)\ot\Xi^{k_5}(k_5/2)))^{\a-k_5}
@>>>gr_{\a-k_5}(\th_!((p_{01}^*\LL_u)\ot gr_0(\Xi^{k_5}(k_5/2))))^{\a-k_5}$$
which has kernel $0$ and cokernel in $\cm^\prec\cb^4$, hence  a morphism
$$\align&gr_0(\th_!((p_{01}^*\LL_u)\ot\Xi^{k_5}))^{\a-k_5}(\a/2))\\&
@>>>gr_0(\th_!((p_{01}^*\LL_u)\ot gr_0(\Xi^{k_5}(k_5/2))))^{\a-k_5}((\a-k_5)/2))\endalign$$
which has kernel $0$ and cokernel in $\cm^\prec\cb^4$. Composing this with the morphism (b) we obtain a
morphism
$$\align&gr_0(\th_!((p_{01}^*\LL_u)\ot\Xi))^\a(\a/2))\\&@>>>
gr_0(\th_!((p_{01}^*\LL_u)\ot gr_0(\Xi^{k_5}(k_5/2))))^{\a-k_5}((\a-k_5)/2))\endalign$$
which has kernel and cokernel in $\cm^\prec\cb^4$. 
Using 7.1(a) and 6.4 this becomes a morphism
$$\op_{z,y'\in W}U_{z,y'}\ot\TT_{z,y'}@>>>
\op_{y\in W}gr_0(\th_!((p_{01}^*\LL_u)\ot\fT_{y,y\i}))^{\a-k_5}((\a-k_5)/2))\tag f$$
which has kernel and cokernel in $\cm^\prec\cb^4$. As in 7.3, the right hand side of (f) is
$$\align&\op_{y\in W}gr_0(\bp_{12}^*(\LL_u\cir\LL_y)[[\nu]]\ot\bp_{34}^*\LL_{y\i})^{\a-k_5}((\a-k_5)/2))\\&
=\op_{y\in W}gr_0(\LL_u\cir\LL_y)[[\nu]])^{\a-k_5}((\a-k_5)/2))\bxt \LL_{y\i}\\&
=\op_{y\in W}(\LL_u\cir\LL_y))^{\{a-\nu\}}\bxt \LL_{y\i}\\&
=\op_{y\in\boc,z\in\boc}\Hom_{\cc^\boc\cb^2}(\LL_z,\LL_u\unb\LL_y)\LL_z\bxt\LL_{y\i}
\op\op_{(y,z)\in W\T(W-\boc)}U'_{z,y\i}\LL_z\bxt\LL_{y\i}\endalign$$
where $U'_{z,y\i}$ are well defined mixed $\bbq$-vector spaces. It follows that we have canonically 
$$U_{z,y'}=\Hom_{\cc^\boc\cb^2}(\LL_z,\LL_u\unb\LL_{y'{}\i})$$
whenever $y'\in\boc,z\in\boc$. This completes the proof of 7.2(c).

\subhead 7.6\endsubhead
We prove 7.2(d) using the isomorphism ${}''\fc:\cy@>\sim>>\fZ$ in 7.4. (We assume again that we are in the 
setup of 6.1 with $e=1,f=1,e'=5$ hence $\e=5$.) The arguments will be similar to those in 7.5. As in 7.4, we
have $\bvt_!\tLL_u=\th_!((p_{45}^*\LL_u)\ot\Xi)$. Here $\th:\cb^5@>>>\cb^4$ is as in 7.4.

From the exact triangle $(\Xi',\Xi,\Xi^{k_5}[-k_5])$ in 6.2(a) we get an exact triangle
$$(\th_!((p_{45}^*\LL_u)\ot\Xi'), \th_!((p_{45}^*\LL_u)\ot\Xi),
\th_!((p_{45}^*\LL_u)\ot\Xi^{k_5}[k_5]))$$
hence an exact sequence
$$\align&(\th_!((p_{45}^*\LL_u)\ot\Xi'))^j@>>>(\th_!((p_{45}^*\LL_u)\ot\Xi))^j\\&@>>>
(\th_!((p_{45}^*\LL_u)\ot\Xi^{k_5}]k_5]))^j@>>>(\th_!((p_{45}^*\LL_u)\ot\Xi'))^{j+1}.\tag a\endalign$$ 
Replacing $\Xi$ by $\Xi'$ in the proof of 7.4(b) given in 7.4 and using that $(\Xi')^h=0$ if $h\ge k_5$ we 
see that 
$$(\th_!((p_{45}^*\LL_u)\ot\Xi'))^j\in\cm^\prec\cb^4\text{ for }j\ge\a.$$
Hence the exact sequence (a) implies that
$$(\th_!((p_{45}^*\LL_u)\ot\Xi))^\a@>>>(\th_!((p_{45}^*\LL_u)\ot\Xi^{k_5}[k_5]))^\a$$
has kernel and cokernel in $\cm^\prec\cb^4$. This induces a homomorphism
$$\align&gr_0(\th_!((p_{45}^*\LL_u)\ot\Xi))^\a(\a/2))
@>>>gr_0(\th_!((p_{45}^*\LL_u)\ot\Xi^{k_5}[k_5]))^\a(\a/2))\\&
=gr_0(\th_!((p_{45}^*\LL_u)\ot\Xi^{k_5}))^{\a-k_5}(\a/2))\tag b\endalign$$
which has kernel and cokernel in $\cm^\prec\cb^4.$ 

From the exact sequence 6.4(a) we get a distinguished triangle
$$\align&(\th_!((p_{45}^*\LL_u)\ot\cw^{-1}(\Xi^{k_5}(k_5/2)))[-k_5],
\th_!((p_{45}^*\LL_u)\ot\Xi^{k_5}(k_5/2))[-k_5],\\&
\th_!((p_{45}^*\LL_u)\ot gr_0(\Xi^{k_5}(k_5/2)))[-k_5]).\endalign$$
Hence we have an exact sequence  
$$\align&(\th_!((p_{45}^*\LL_u)\ot\cw^{-1}(\Xi^{k_5}(k_5/2)))[-k_5])^\a@>>> 
(\th_!((p_{45}^*\LL_u)\ot\Xi^{k_5}(k_5/2))[-k_5])^\a\\&
@>>>(\th_!((p_{45}^*\LL_u)\ot gr_0(\Xi^{k_5}(k_5/2)))[-k_5])^\a\\&@>>>
(\th_!((p_{45}^*\LL_u)\ot\cw^{-1}(\Xi^{k_5}(k_5/2)))[-k_5])^{\a+1}.\tag c\endalign$$
Replacing $\Xi$ by $\cw^{-1}(\Xi^{k_5}(k_5/2))[-k_5/2]$ in the proof of 7.4(b) given in 7.4 and using that
$(\cw^{-1}(\Xi^{k_5}(k_5/2))[-k_5/2])^h=0$ if $h>k_5$ we see that
$$(\th_!((p_{45}^*\LL_u)\ot\cw^{-1}(\Xi^{k_5}(k_5/2)))[-k_5])^{\a+1}\in\cm^\prec\cb^4.\tag d$$
Note that
$$gr_{\a-k_5}(\th_!((p_{45}^*\LL_u)\ot\cw^{-1}(\Xi^{k_5}(k_5/2))))^{\a-k_5}=0.\tag e$$
This follows from the fact that $\cw^{-1}(\Xi^{k_5}(k_5/2))$ has weight $\le-1$ hence \lb
$\th_!((p_{45}^*\LL_u)\ot\cw^{-1}(\Xi^{k_5}(k_5/2)))$ has weight $\le-1$  and \lb
$(\th_!((p_{45}^*\LL_u)\ot\cw^{-1}(\Xi^{k_5}(k_5/2))))^{\a-k_5}$ has weight $\le\a-k_5-1$.

Using (d),(e), we see that (c) induces a morphism
$$\align&gr_{\a-k_5}(\th_!((p_{45}^*\LL_u)\ot\Xi^{k_5}(k_5/2)))^{\a-k_5}\\&
@>>>gr_{\a-k_5}(\th_!((p_{45}^*\LL_u)\ot gr_0(\Xi^{k_5}(k_5/2))))^{\a-k_5}\endalign$$
which has kernel $0$ and cokernel in $\cm^\prec\cb^4$, hence  a morphism
$$\align&gr_0(\th_!((p_{45}^*\LL_u)\ot\Xi^{k_5}))^{\a-k_5}(\a/2))\\&
@>>>gr_0(\th_!((p_{45}^*\LL_u)\ot gr_0(\Xi^{k_5}(k_5/2))))^{\a-k_5}((\a-k_5)/2))\endalign$$
which has kernel $0$ and cokernel in $\cm^\prec\cb^4$. Composing this with the morphism (b) we obtain a
morphism
$$\align&gr_0(\th_!((p_{45}^*\LL_u)\ot\Xi))^\a(\a/2))\\&@>>>
gr_0(\th_!((p_{45}^*\LL_u)\ot gr_0(\Xi^{k_5}(k_5/2))))^{\a-k_5}((\a-k_5)/2))\endalign$$
which has kernel and cokernel in $\cm^\prec\cb^4$. Using 7.1(a) and 6.4, this becomes a morphism
$$\op_{y',z\in W}U_{y',z}\ot\TT_{y',z}@>>>
\op_{y\in W}gr_0(\th_!((p_{45}^*\LL_u)\ot\fT_{y,y\i}))^{\a-k_5}((\a-k_5)/2))\tag f$$
which has kernel and cokernel in $\cm^\prec\cb^4$. As in 7.4, the right hand side of (f) is
$$\align&\op_{y\in W}gr_0(\bp_{34}^*(\LL_{y\i}\cir\LL_u)[[\nu]]\ot\bp_{12}^*\LL_y)^{\a-k_5}((\a-k_5)/2))\\&
=\op_{y\in W} \LL_y\bxt gr_0(\LL_{y\i}\cir\LL_y)[[\nu]])^{\a-k_5}((\a-k_5)/2))\\&
=\op_{y\in W}\LL_y\bxt (\LL_{y\i}\cir\LL_u))^{\{a-\nu\}} \\&
=\op_{y\in\boc,z\in\boc}\Hom_{\cc^\boc\cb^2}(\LL_z,\LL_{y\i}\unb\LL_u)\ot(\LL_y\bxt\LL_z)
\op\op_{(y,z)\in W\T(W-\boc)}U''_{y,z}\ot(\LL_y\bxt\LL_z)\endalign$$
where $U''_{y,z}$ are well defined mixed $\bbq$-vector spaces.
It follows that we have canonically $U_{y,z}=\Hom_{\cc^\boc\cb^2}(\LL_z,\LL_{y\i}\unb\LL_u)$
whenever $y\in\boc,z\in\boc$. This completes the proof of 7.2(d). Lemma 7.2 is proved.

\proclaim{Proposition 7.7} For any $y,z,u\in\boc$ we have canonically
$$\Hom_{\cc^\boc\cb^2}(\LL_y,\LL_u\unb\LL_{z\i})=\Hom_{\cc^\boc\cb^2}(\LL_z,\LL_{y\i}\unb\LL_u).\tag a$$
\endproclaim
Indeed both sides of (a) are identified in 7.2(c),(d) with $U_{y,z}$.

\proclaim{Proposition 7.8} Let $u,x\in\boc$. In the setup of 7.1 we have canonically
$$\un{(\bp_{14!}(\bvt_!(\tLL_u)\ot\bp_{23}^*\LL_x))^{\{3a+\nu+2\r\}}}=\op_{y,z\in\boc}\bQ_{y,z}$$
where 
$$\align&\bQ_{y,z}=\Hom_{\cc^\boc\cb^2}(\LL_y,\LL_u\unb\LL_{z\i})\ot(\LL_y\unb\LL_x\unb\LL_z)\\&
=\Hom_{\cc^\boc\cb^2}(\LL_z,\LL_{y\i}\unb\LL_u)(\LL_y\unb\LL_x\unb\LL_z)\in\cc^\boc_0\cb^2.\endalign$$
(The last equality comes from 7.7.)
\endproclaim
Define $\Ph:\cd_m^\prq\cb^4@>>>\cd_m^\prq\cb^2$ by $\Ph(K)=\bp_{14!}(K\ot\bp_{23}^*\LL_x)$. 
This is well defined and maps $\cd_m^\prec\cb^4$ to $\cd_m^\prec\cb^2$. (This can be deduced from 2.2(a),
(e).) Let $(c,c')=(2a-2\nu,a+3\nu+2\r)$.
Let $\XX=\bvt_!(\tLL_u)$. By 7.2(a) we have $\XX^j\in\cm^\prec\cb^4$ for any $j>c'$. Note that $\XX$ has 
weight $\le0$. If $K\in\cd_m^\prq\cb^4$ and $K\in\cm^\prq\cb^4$ then $(\Ph(K))^h\in\cm^\prec\cb^4$ for any 
$h>c$. (This can be deduced from 2.2(a) with $r=3$.)
Now the proof of Lemma 1.12 can be repeated word by word and yield a canonical identification
$$\un{(\Ph(\un{\XX^{\{c'\}}}))^{\{c\}}}=\un{(\Ph(\XX))^{\{c+c'\}}}$$
that is
$$\un{(\bp_{14!}(\un{(\bvt_!(\tLL_u))^{\{a+3\nu+2\r\}}} \ot\bp_{23}^*\LL_x))^{\{2a-2\nu\}}}=
\un{(\bp_{14!}(\bvt_!(\tLL_u)\ot\bp_{23}^*\LL_x))^{\{3a+\nu+2\r\}}}.$$
Replacing here $\un{(\bvt_!(\tLL_u))^{\{a+3\nu+2\r\}}}$ by
$$\op_{y,z\in\boc}U_{y,z}\ot\TT_{y,z}=\op_{y,z\in\boc}U_{y,z}\ot\bp_{12}^*\LL_y\ot\bp_{34}^*\LL_z$$
(see 7.1(a) and 7.2(a)) we obtain
$$\un{(\bp_{14!}(\bvt_!(\tLL_u)\ot\bp_{23}^*\LL_x))^{\{3a+\nu+2\r\}}}=\op_{y,z\in\boc}\bQ_{y,z}$$
where
$$\bQ_{y,z}=U_{y,z}\ot\un{(\bp_{14!}(\bp_{12}^*\LL_y\ot\bp_{34}^*\LL_z \ot\bp_{23}^*\LL_x))^{\{2a-2\nu\}}}=
U_{y,z}\ot(\LL_y\unb\LL_x\unb\LL_z).$$
This completes the proof. (We use 7.2(c),(d).)

\subhead 7.9\endsubhead
Let 
$${}'\cy=\{((B_0,B_1,B_2,B_3,B_4),g)\in\cb^5\T G;gB_1g\i=B_4,gB_2g\i=B_3\},$$
$${}''\cy=\{((B_1,B_2,B_3,B_4,B_5),g)\in\cb^5\T G;gB_1g\i=B_4,gB_2g\i=B_3\}.$$
Note that ${}'\cy$ is what in 6.1 (with $e=0,f=1,e'=4$) was denoted by $\cy$ and ${}''\cy$ is what in 6.1 
(with $e=1,f=1,e'=5$) was denoted by $\cy$.
For $i,j$ in $[0,4]$ define ${}'h_{ij}:{}'\cy@>>>\cb^2$ by $((B_0,B_1,B_2,B_3,B_4),g)\m(B_i,B_j)$.
For $i,j$ in $[1,5]$ define ${}''h_{ij}:{}''\cy@>>>\cb^2$ by $((B_1,B_2,B_3,B_4,B_5),g)\m(B_i,B_j)$.
Let $u,x\in\boc$. Let ${}'\ce={}'h_{04!}({}'h_{01}^*\LL_u\ot{}'h_{23}^*\LL_x)\in\cd_m\cb^2$,
${}''\ce={}''h_{15!}({}''h_{23}^*\LL_x)\ot{}''h_{45}^*\LL_u)\in\cd_m\cb^2$.
From Lemma 6.7 we obtain canonical identifications
$$\un{({}'\ce)^{\{3a+\nu+2\r\}}}=\op_{y\in\boc}{}'Q_y,\qua 
\un{({}''\ce)^{\{3a+\nu+2\r\}}}=\op_{y\in\boc}{}''Q_y,$$
where
$${}'Q_y=\LL_u\unb\LL_y\unb\LL_x\unb\LL_{y\i},\qua {}''Q_y=\LL_y\unb\LL_x\unb\LL_{y\i}\unb\LL_u.$$
Using Theorem 6.8 we have canonically 
$$\op_{y\in\boc}{}'Q_y=\LL_u\unb\unz\unc(\LL_x),\qua \op_{y\in\boc}{}''Q_y=\unz\unc(\LL_x)\unb\LL_u$$
hence 
$$\un{({}'\ce)^{\{3a+\nu+2\r\}}}=\LL_u\unb\unz\unc(\LL_x),\qua 
\un{({}''\ce)^{\{3a+\nu+2\r\}}}=\unz\unc(\LL_x)\unb\LL_u.$$
From the definitions we see that the identification
$$\LL_u\unb\unz\unc(\LL_x)=\unz\unc(\LL_x)\unb\LL_u\tag a$$
in 3.4(a) (with $L=\LL_u,K=\unc(\LL_x)$) is the same as the identification
$$\un{({}'\ce)^{\{3a+\nu+2\r\}}}=\un{({}''\ce)^{\{3a+\nu+2\r\}}}$$
obtained by identifying both sides with $\un{\ce^{\{3a+\nu+2\r\}}}$ 
where $\ce=\bp_{14!}(\bvt_!(\tLL_u)\ot\bp_{23}^*\LL_x)$.
(Note that ${}'\ce=\ce={}''\ce$ via the isomorphisms ${}'\cy@>{}'\fc>>\fZ@<{}''\fc<<{}''\cy$, see 7.3, 7.4,
where ${}'\cy,{}''\cy$ are denoted by $\cy$.)
Using these identifications and Proposition 7.8 we obtain a commutative diagram
$$\CD
\LL_u\unb\unz\unc(\LL_x)@>\si>>\un{\ce^{\{3a+\nu+2\r\}}}@<\si<<\unz\unc(\LL_x)\unb\LL_u\\
@V\si VV                               @V\si VV                           @V\si VV     \\
\op_{y\in\boc}{}'Q_y @>\si>>\op_{y,z\in\boc}\bQ_{y,z}@<\si<<\op_{y\in\boc}{}''Q_y
\endCD$$
where the upper horizontal maps yield the identification (a) and the lower horizontal maps are the obvious 
ones: they map
${}'Q_y$ onto $\op_{z\in\boc}\bQ_{z\i,y\i}$ and ${}''Q_y$ onto $\op_{z\in\boc}\bQ_{y,z}$.

\head 8. Adjunction formula (weak form)\endhead
\proclaim{Proposition 8.1} Let $L\in\cc^\boc_0\cb^2,K\in\cc^\boc_0G$. We have canonically  
$$K\unst\unc(L)=\unc(L\unb\unz(K)).\tag a$$
\endproclaim
Applying 1.12 with $\Ph:\cd_m^\prq G@>>>\cd_m^\prq G$, 
$K_1\m K*K_1$, $\XX=\c(L)$, $(c,c')=(2a+\r,a+\r+\nu)$ (see 4.5, 1.9) we deduce that we have canonically 
$$\un{(K*\un{(\c(L))^{\{a+\r+\nu\}}})^{\{2a+\r\}}}=\un{(K*\c(L))^{\{3a+2\r+\nu\}}}$$
that is,
$$K\unst\unc(L)=\un{(K*\c(L))^{\{3a+2\r+\nu\}}}.\tag b$$
Applying 1.12 with $\Ph:\cd_m^\prq\cb^2@>>>\cd_m^\prq\cb^2$, ${}^1L\m L\cir{}^1L$, $\XX=\z(K)$,
$(c,c')=(a-\nu,a+\nu+\r)$ (see 3.1, 2.8) we deduce that we have canonically 
$$\un{(L\cir\un{(\z(K))^{\{a+\nu+\r\}}})^{\{a-\nu\}}}=\un{(L\cir\z(K))^{\{2a+\r\}}}\tag c$$
and
$$(L\cir\z(K))^j\in\cm^\prec\cb^2\text{ if }j>2a+\r.\tag d$$
Applying 1.12 with $\Ph:\cd_m^\prq\cb^2@>>>\cd_m^\prq G$, ${}^1L\m\c({}^1L)$, $\XX=L\cir\z(K)$,
$(c,c')=(a+\r+\nu,2a+\r)$ (see (d) and 1.9) we deduce that we have canonically 
$$\un{\c(\un{(L\cir\z(K))^{\{2a+\r\}}})^{\{a+\r+\nu\}}}=\un{(\c(L\cir\z(K)))^{\{3a+2\r+\nu\}}}.$$
Combining this with (c) gives
$$\unc(L\unb\unz(K))=\un{(\c(L\cir\z(K)))^{\{3a+2\r+\nu\}}}$$
which together with (b) gives (a). (We use the equality  $K*\c(L)=\c(L\cir\z(K))$, see 4.2.)

\mpb

The following lemma is a variant of 1.12.

\proclaim{Lemma 8.2} Let $c\in\ZZ$ and let $Y$ be one of $G,\cb^2$.
Let $\Ph:\cd_m^\prq Y@>>>\cd_m\pp$ be a functor which takes 
distinguished triangles to distinguished triangles, commutes with shifts and direct sums and maps 
complexes of weight $\le i$ to complexes of weight $\le i$ (for any $i$). Assume that
$$(\Ph(K_0))^h=0\text{ for any }K_0\in\cm_m^\prq Y\text{ and any }h>c.\tag a$$
Then for any $K\in\cd_m^\prq Y$ of weight $\le0$ and any $c'\in\ZZ$ we have canonically
$$(\Ph(\un{K^{\{c'\}}}))^{\{c\}}\sub(\Ph(K))^{\{c+c'\}}.\tag c$$
\endproclaim
As in 1.12 for any $i,h$ we have an exact sequence
$$(\Ph(K^i))^{h-1}@>>>(\Ph(\t_{<i}K))^{i+h}@>>>(\Ph(\t_{\le i}K))^{i+h}
@>>>(\Ph(K^i))^h@>>>(\Ph(\t_{<i}K))^{i+h+1}.$$
Assume first that $i+h=c+c'+1$, $h\ge c+2$ hence $i\le c'-1$. Then 
$(\Ph(K^i))^{h-1}=0,(\Ph(K^i))^h=0$ hence $(\Ph(\t_{<i}K))^{c+c'+1}@>\si>>(\Ph(\t_{\le i}K))^{c+c'+1}.$
Thus we see by induction on $i$ that $(\Ph(\t_{\le i}K))^{c+c'+1}=0$ for $i\le c'-1$; in particular
$$(\Ph(\t_{\le c'-1}K))^{c+c'+1}=0.\tag d$$
Next assume that  $i+h=c+c'$, $h\ge c+2$ hence $i\le c'-2$. Then 
$(\Ph(K^i))^{h-1}=0,(\Ph(K^i))^h=0$ hence $(\Ph(\t_{<i}K))^{c+c'+1}@>\si>>(\Ph(\t_{\le i}K))^{c+c'+1}.$
Thus we see by induction on $i$ that $(\Ph(\t_{\le i}K))^{c+c'}=0$ for $i\le c'-2$; in particular
$(\Ph(\t_{\le c'-2}K))^{c+c'}=0$. Now assume that $i+h=c+c'$, $h=c+1$ hence $i=c'-1$. We have an exact 
sequence $(\Ph(\t_{\le c'-2}K))^{c+c'}@>>>(\Ph(\t_{\le c'-1}K))^{c+c'}@>>>0$ hence
$(\Ph(\t_{\le c'-1}K))^{c+c'}=0$. Now assume that  $i+h=c+c'$, $h=c$ hence $i=c'$. We have an exact sequence
$$0@>>>(\Ph(\t_{\le c'}K))^{c+c'}@>>>(\Ph(K^{c'}))^c@>>>(\Ph(\t_{<c'}K))^{c+c'+1},$$
hence using (d) we have 
$$(\Ph(\t_{\le c'}K))^{c+c'}@>\si>>(\Ph(K^{c'}))^c.\tag e$$
For any $i$, from the exact sequence 
$$(\Ph(K^i))^{c+c'-i-1}@>>>(\Ph(\t_{<i}K))^{c+c'}@>>>(\Ph(\t_{\le i}K))^{c+c'}$$
we deduce an exact sequence
$$gr_{c+c'}(\Ph(K^i))^{c+c'-i-1}@>>>gr_{c+c'}(\Ph(\t_{<i}K))^{c+c'}@>>>gr_{c+c'}(\Ph(\t_{\le i}K))^{c+c'}.$$
Now $(\Ph(K^i))^{c+c'-i-1}$ is mixed of weight $\le c+c'-1$ (by our assumptions) hence
$gr_{c+c'}(\Ph(K^i))^{c+c'-i-1}=0$. Thus for any $i$  we have an imbedding
$$gr_{c+c'}(\Ph(\t_{<i}K))^{c+c'}\sub gr_{c+c'}(\Ph(\t_{\le i}K))^{c+c'}.$$   
Hence each $gr_{c+c'}(\Ph(\t_{<i}K))^{c+c'}$ becomes a subobject of $gr_{c+c'}(\Ph(\t_{\le i}K))^{c+c'}$   
with large $i$, that is of $gr_{c+c'}(\Ph(K))^{c+c'}$. In particular we have 
$$gr_{c+c'}(\Ph(\t_{\le c'}K))^{c+c'}\sub gr_{c+c'}(\Ph(K))^{c+c'}.\tag f$$ 
From the exact sequence  $0@>>>W_{c'-1}K^{c'}@>>>K^{c'}@>>>gr_{c'}K^{c'}@>>>0$ (here we use that $K^{c'}$
has weight $\le c'$) we deduce an exact sequence
$$(\Ph(W_{c'-1}K^{c'}))^c@>>>(\Ph(K^{c'}))^c@>>>(\Ph(gr_{c'}K^{c'}))^c@>>>(\Ph(W_{c'-1}K^{c'}))^{c+1}$$
hence an exact sequence
$$\align&gr_{c+c'}(\Ph(W_{c'-1}K^{c'}))^c@>>>gr_{c+c'}(\Ph(K^{c'}))^c@>>>gr_{c+c'}(\Ph(gr_{c'}K^{c'}))^c
\\&@>>>gr_{c+c'}(\Ph(W_{c'-1}K^{c'}))^{c+1}.\endalign$$
Now $(\Ph(W_{c'-1}K^{c'}))^c$ has weight $\le c+c'-1$ hence 
$gr_{c+c'}(\Ph(W_{c'-1}K^{c'}))^c=0$; by (a) we have $(\Ph(W_{c'-1}K^{c'}))^{c+1}=0$. Hence the previous 
exact sequence yields
$$gr_{c+c'}(\Ph(K^{c'}))^c@>\si>>gr_{c+c'}(\Ph(gr_{c'}K^{c'}))^c.$$
Combining this with 
$gr_{c+c'}(\Ph(\t_{\le c'}K))^{c+c'}=gr_{c+c'}(\Ph(K^{c'}))^c$ obtained from (e) we see that
$$gr_{c+c'}(\Ph(gr_{c'}K^{c'}))^c=gr_{c+c'}(\Ph(\t_{\le c'}K))^{c+c'}.$$
Using this and (f) we obtain an imbedding
$$gr_{c+c'}(\Ph(gr_{c'}K^{c'}))^c\sub gr_{c+c'}(\Ph(K))^{c+c'}.$$ 
Since $\un{gr_{c'}K^{c'}}$ is canonically a direct summand of $gr_{c'}K^{c'}$ we see that the previous
imbedding restricts to an imbedding
$$gr_{c+c'}(\Ph(\un{gr_{c'}K^{c'}}))^c\sub gr_{c+c'}(\Ph(K))^{c+c'}.$$ 
Applying $((c+c')/2)$ to both sides we obtain (c).

\subhead 8.3\endsubhead
Let $\io:\pp@>>>G$ be the map with image $1$. We show:

(a) {\it Let $K\in\cm^\prq_mG$. If $j>-2a-\r$, then $(\io^*(K))^j=0$.}
\nl
We can assume that $K\in CS(G)$. From the cleanness of cuspidal 
character sheaves we see that either $\io^*K=0$ in which case there is nothing to prove, or $K\cong A_E$ for
some $E\in\Irr W$ which we now assume. We have $\ch^i_1A_E=\Hom_W(E,H^{i+\D}(\cb,\bbq))(\D/2)$ where 
$H^{i+\D}(\cb,\bbq)$ has the natural $W$-action. It is known that the polynomial 
$\sum_{k\ge0}\dim\Hom_W(E,H^k(\cb,\bbq))v^k$ has degree $\le2\nu-2\aa(\boc_E)$. Hence 
$\sum_i\dim(\ch^i_1(A_E))v^i\in v^{-2\aa(\boc_E)-\r}\ZZ[v\i]$. Since $\boc_E\prq\boc$, we have 
$\aa(\boc_E)\ge a$ and $\sum_i\dim(\ch^i_1(A_E))v^i\in v^{-2a-\r}\ZZ[v\i]$. This proves (a). 

\mpb

(b) {\it If $K=A_{E_\boc}$, then we have canonically $(\io^*K)^{-2a-\r}=\EE((2a+\r)/2)$ where $\EE$ is a
well defined $1$-dimensional $\bbq$-vector space of pure weight $0$.}
\nl
Equivalently, $\ch^{-2a-\r}_1K$ is a one dimensional mixed $\bbq$-vector space of pure weight $-2a-\r$. (We 
use the fact that $E_\boc$ appears in the $W$-mdule $H^{-2a+2\nu}(\cb,\bbq)(\D/2)$ with multiplicity one and 
that $H^{-2a+2\nu}(\cb,\bbq)(\D/2)$ is pure of weight $-2a-\r$.)

(c) {\it If $K\in\cc^\boc G$ and $\Hom_{\cc^\boc G}(A_{E_\boc},K)=0$  then $(\io^*(K))^{-2a-\r}=0$.}
\nl
We can assume that $K=A_E$ where $E\in\Irr_\boc W$, $E\ne E_\boc$. We then use the fact 
that $E$ does not appear in the $W$-module $H^{-2a+2\nu}(\cb,\bbq)(\D/2)$.

\subhead 8.4\endsubhead
Define $\d:\cb@>>>\cb^2$ by $B\m(B,B)$; let $\o:\cb@>>>\pp$ be the obvious map. From the definitions, for 
any $L\in\cd_m\cb^2$ we have canonically
$$\io^*(\c(L))=\o_!\d^*(L).\tag a$$
We show:

(b) {\it Let $L\in\cm^\prq_m\cb^2$. If $j>-a$ then $(\d^*L)^j=0$.}
\nl
We can assume that $L=\LL_w$ where $w\prq\boc$. It is enough to show that for any $k$ we have 
$(\ch^k(\d^*L)[-k]))^j=0$ that is $(\ch^k(\d^*(L_w^\sh[|w|+\nu])))^{j-k}=0$ or equivalently 
$(\ch^{k+\nu}(\d^*(L_w^\sh[|w|]))[\nu])^{j-k-\nu}=0$. Now $\ch^{k+\nu}(\d^*(L_w^\sh[|w|]))[\nu]$ is a 
perverse sheaf hence we can take $k=j-\nu$ and it is enough to prove that $\ch^j(\d^*(L_w^\sh[|w|]))=0$. Now 
$$\sum_{i\le0}\rk(\ch^i(\d^*L_w^\sh[|w|]))v^i=p_{1,w}\in v^{-\aa(w)}\ZZ[v\i]$$
with $p_{1,w}$ as in \cite{\HEC, 5.3} (see \cite{\HEC, 14.2, P1}). Since $\aa(w)\ge a$ it follows that 
$p_{1,w}\in v^{-a}\ZZ[v\i]$. This proves (b).

We show:

(c) {\it If $L\in\cm^\prq_m\cb^2$ is pure of weight $0$ and $i\in\ZZ$ then $(\d^*L)^i$ is pure of weight 
$i$.}
\nl
We can assume that $L=\LL_w$ where $w\prq\boc$. We have
$(\d^*L)^i=\ch^{i-\nu}(\d^*L)=\ch^{i+|w|}(\d^*L_w^\sh)(|w|/2)$ hence it is enough to show that, setting
$j=i+|w|$, $\ch^j(\d^*L_w^\sh)$ is pure of weight $j$. This follows from the results in \cite{\KL}.

(d) {\it Assume that $w\in\boc$. If $w=d\in\DD_\boc$ then $(\d^*\LL_w)^{-a}=\BB_d[[\nu]](a/2)$ for a well 
defined one dimensional mixed $\bbq$-vector space $\BB_d$ of pure weight $0$, noncanonically isomorphic to 
$\bbq$. If $w\n\DD_\boc$ then $(\d^*\LL_w)^{-a}=0$.}
\nl
In view of (c), an equivalent statement is that the coefficient of $v^{-a}$ in $p_{1,w}$ is $1$ if
$w\in\DD_\boc$ and is $0$ if $w\n\DD_\boc$; this holds by \cite{\HEC, 14.2, P5}.

\subhead 8.5\endsubhead

(a) {\it Assume that $L\in\cm_m\cb$ is $G$-equivariant so that $L=V\ot\bbq[[\nu]]$ where $V$ is a 
mixed $\bbq$-vector space. If $j>\nu$ then $(\o_!L)^j=0$. We have $(\o_!L)^\nu=V(-\nu)$.}
\nl
We have $\ch^j(\o_!L)=V\ot H^{j+\nu}(\cb,\bbq)$. Since $\dim\cb=\nu$, this is zero if $j+\nu>2\nu$ and is
$V(-\nu)$ if $j+\nu=2\nu$. This proves (a).

We show:  

(b) {\it If $L\in\cm_m^\prq\cb^2$ and $j>\nu-a$ then $(\o_!\d^*L)^j=0$. Moreover, we have canonically
$(\o_!\d^*L)^{\nu-a}=(\o_!*((\d^*L)^{-a}))^\nu$.}     
\nl
We set $\XX=\d^*L$. As in the proof of 1.12 we have an exact sequence
$$\o_!(\XX^i))^{h-1}@>>>(\o_!(\t_{<i}\XX))^{i+h}@>>>(\o_!(\t_{\le i}\XX))^{i+h}@>>>(\o_!(\XX^i))^h
@>>>(\o_!(\t_{<i}\XX))^{i+h+1}.$$
From this we see by induction on $i$ (using 8.3 and (a)) that if $j>\nu-a$ then\lb
$(\o_!(\t_{\le i}\XX))^j=0$ for any $i$. Hence the first assertion of (b) holds.
Assume now that $i+h=\nu-a$. From the exact sequence above we see (using 8.3) that
$$(\o_!(\t_{<i}\XX))^{\nu-a}@>\si>>(\o_!(\t_{\le i}\XX))^{\nu-a}$$
when $i>-a$ hence $(\o_!(\t_{\le-a}\XX))^{\nu-a}@>\si>>(\o_!\XX)^{\nu-a}$.
From the same exact sequence we see by induction on $i$ (using (a)) that
$(\o_!(\t_{\le i}\XX))^j=0$ for $i\le-a-1$ hence $(\o_!(\t_{\le-a-1}\XX))^j=0$.
The exact sequence above with $i=-a,h=\nu$ becomes
$$0@>>>(\o_!\XX)^{\nu-a}@>>>(\o_!(\XX^{-a}))^\nu@>>>(\o_!(\t_{<-a}\XX))^{\nu-a+1}.$$
Hence we obtain an isomorphism $(\o_!\XX)^{\nu-a}@>\si>>(\o_!(\XX^{-a}))^\nu$.

\subhead 8.6\endsubhead
Let $L\in\cc^\boc_0\cb^2$. Applying 8.2 with $\Ph:\cd_m^\prq G@>>>\cd_m\pp$, $K_1\m\io^*K_1$, $c=-2a-\r$ 
(see 8.3), $K$ replaced by $\c(L)$ and $c'=a+\nu+\r$ we see that we have canonically
$$(\io^*(\unc(L)))^{\{-2a-\r\}}\sub(\io^*\c(L))^{\{-a+\nu\}}=(\o_!\d^*(L))^{\{-a+\nu\}}.$$
(The last equality comes from 8.4(a).) We set 
$$\bold1'=\op_{d\in\DD_\boc}\BB_d^*\ot\LL_d\in\cc^\boc_0\cb^2$$
where $\BB_d^*$ is the vector space dual to $\BB_d$. From 8.4(d), 8.5, we see that
$$(\o_!\d^*(L))^{\{-a+\nu\}}=(\o_!((\d^*L)^{-a}))^\nu((\nu-a)/2)=\Hom_{\cc^\boc\cb^2}(\bold1',L).\tag a$$
Hence we have canonically
$$(\io^*(\unc(L)))^{\{-2a-\r\}}\sub\Hom_{\cc^\boc\cb^2}(\bold1',L).\tag b$$
We show that the last inclusion is an equality:
$$(\io^*(\unc(L)))^{\{-2a-\r\}}=\Hom_{\cc^\boc\cb^2}(\bold1',L).\tag c$$
To prove this we can assume that $L=\LL_x$ for some $x\in\boc$. If $x\n\DD_\boc$ then the right hand side of
(b) is zero hence the left hand side of (b) is zero and (c) holds. Assume now that $x\in\DD_\boc$. Then the 
right hand side of (b) has dimension $1$; to prove (c) it is enough to show that the left hand side of (b) 
has dimension $1$. By 8.3(b),(c), the left hand side of (b) has dimension
$(A_{E_\boc}:\unc(\LL_x))$ which, as we already know from (b), has dimension $0$ or $1$. Using 1.15(a) we 
see that this dimension is in fact $1$. This proves (c).

The argument above shows also that the assumption of 1.15(a) is satisfied; hence we can now state
unconditionally:

(d) {\it For any $d\in\DD_\boc$ we have $(A_{E_\boc}:\unc(\LL_d))=1$.}
\nl
The argument above shows also:

(e) {\it For any $x\in\boc-\DD_\boc$ we have $(A_{E_\boc}:\unc(\LL_x))=0$.} 

\proclaim{Lemma 8.7} Let $L,L'\in\cc^\boc\cb^2$. We have canonically
$$\Hom_{\cc^\boc\cb^2}(\bold1',L\unb L')=\Hom_{\cc^\boc\cb^2}(\fD(L'{}^\da),L).\tag a$$ 
Here for ${}^1L\in\cc^\boc\cb^2$ or ${}^1L\in\cc^\boc_0\cb^2$ we set $L^\da=h'{}^*L$ where 
$h':\cb^2@>>>\cb^2$ is $(B,B')\m(B',B)$.
\endproclaim
We can assume that $L=\LL_x,L'=\LL_{x'}$ with $x,x'\in\boc$. We view $L,L'$ as objects of
$\cc^\boc_0\cb^2$. Using 8.4(a) we have
$$\Hom_{\cc^\boc\cb^2}(\bold1',L\unb L')=(\o_!\d^*(L\unb L'))^{\{-a+\nu\}}.\tag b$$
Applying 8.2 with $\Ph:\cd^\prq_m\cb^2@>>>\cd_m\pp$, $\tL\m\o_!\d^*\tL$, $(c,c')=(\nu-a,a-\nu)$, see
8.5(b), $K=L\cir L'$, we deduce that we have canonically
$$(\o_!\d^*(L\unb L'))^{\{\nu-a\}}\sub(\o_!\d^*(L\cir L'))^{\{0\}}.\tag c$$
From \cite{\CSII, 7.4} we see that we have canonically
$$(\o_!(L\ot L'{}^\da))^0=(\o_!(L\ot L'{}^\da))^{\{0\}}=\Hom_{\cc^\boc\cb^2}(\fD(L'{}^\da),L).\tag d$$
Note that $\d^*(L\cir L')=L\ot L'{}^\da$. Hence by combining (b),(c),(d) we have
$$\Hom_{\cc^\boc\cb^2}(\bold1',L\unb L')\sub\Hom_{\cc^\boc\cb^2}(\fD(L'{}^\da),L).\tag e$$
The dimension of the left hand side of (e) is the sum over $d\in\DD_\boc$ of the coefficients of $t_d$ in 
$t_xt_{x'}\in\JJ^\boc$ and, by \cite{\HEC, 14.2}, this sum is equal to $1$ if $x'{}\i=x$ and $0$ if 
$x'{}\i\ne x$; hence it is equal to the dimension of the right hand side of (e). It follows that (e) is an
equality and (a) follows.

\subhead 8.8\endsubhead
Let $u:G@>>>\pp$ be the obvious map. 
From \cite{\CSII, 7.4}, we see that for $K,K'\in\cm_m^\prq G$ we have canonically
$$(u_!(K\ot K'))^0=\Hom_{\cm G}(\fD(K),K'),\qua (u_!(K\ot K'))^j=0\text { if }j>0.$$
We deduce that if $K,K'$ are also pure of weight $0$ then $(u_!(K\ot K'))^0$ is pure of weight zero that is
$(u_!(K\ot K'))^0=gr_0(u_!(K\ot K'))^0$.  
From the definitions we see that we have $u_!(K\ot K')=\io^*(K^\da*K')$ where $K^\da=h^*K$ and $h:G@>>>G$ is
given by $g\m g\i$. Hence for $K,K'\in\cc^\boc_0G$ we have
$$\Hom_{\cc^\boc G}(\fD(K),K')=(\io^*(K^\da*K'))^0=(\io^*(K^\da*K'))^{\{0\}}.\tag a$$
Applying 8.2 with $\Ph:\cd_m^\prq G@>>>\cd_m\pp$, $K_1\m\io^*K_1$, $c=-2a-\r$ (see 8.3), $K$ replaced by
$K^\da*K'$ and $c'=2a+\r$ we see that we have canonically
$$(\io^*(K^\da\unst K'))^{\{-2a-\r\}}\sub(\io^*(K^\da*K'))^{\{0\}}.$$
In particular if $L,L'\in\cc^\boc_0\cb^2$ then we have canonically  
$$(\io^*(\unc(L')\unst\unc(L)))^{\{-2a-\r\}}\sub(\io^*(\unc(L')*\unc(L)))^{\{0\}}.$$
Using the equality 
$$(\io^*(\unc(L')\unst\unc(L)))^{\{-2a-\r\}}=(\io^*(\unc(L\unb\unz(\unc(L')))))^{\{-2a-\r\}}$$
which comes from 8.1 we deduce that we have canonically
$$(\io^*(\unc(L\unb\unz(\unc(L')))))^{\{-2a-\r\}}\sub(\io^*(\unc(L')*\unc(L)))^{\{0\}}$$
or equivalently, using (a) with $K,K'$ replaced by $\unc(L')^\da$, $\unc(L)$:
$$\align&(\io^*(\unc(L\unb\unz(\unc(L')))))^{\{-2a-\r\}}\sub\Hom_{\cc^\boc G}(\fD(\unc(L')^\da),\unc(L))\\&=
\Hom_{\cc^\boc G}(\fD(\unc(L)^\da),\unc(L')).\endalign$$
Using now 8.6(c) we deduce that we have canonically
$$\Hom_{\cc^\boc\cb^2}(\bold1',L\unb\unz\unc L')\sub\Hom_{\cc^\boc G}(\fD(\unc(L)^\da),\unc(L'))$$
or equivalently (see 8.7)
$$\Hom_{\cc^\boc\cb^2}(\fD(L^\da),\unz\unc L')\sub\Hom_{\cc^\boc G}(\fD(\unc(L)^\da),\unc(L')).$$
We now set ${}^1L=\fD(L^\da)$ and note that 
$$\fD(\unc(L)^\da)=\fD(\unc(L^\da))=\unc(\fD(L^\da))=\unc({}^1L),$$
see 1.13(a). We obtain
$$\Hom_{\cc^\boc\cb^2}({}^1L,\unz\unc L')\sub\Hom_{\cc^\boc G}(\unc({}^1L),\unc(L'))\tag b$$
for any ${}^1L,L'\in\cc^\boc_0\cb^2$.

We have the following result which is a weak form of an adjunction formula, of which the full form will be 
proved in 9.8.

\proclaim{Proposition 8.9}For any ${}^1L,L'\in\cc^\boc_0\cb^2$ we have canonically
$$\Hom_{\cc^\boc\cb^2}({}^1L,\unz\unc(L'))=\Hom_{\cc^\boc G}(\unc({}^1L),\unc(L'))\tag a$$
\endproclaim
We can assume that ${}^1L=\LL_z,L'=\LL_u$ where $z,u\in\boc$. By 6.9(a) and 1.10(b), both sides of the 
inclusion 8.8(b) have dimension $\sum_{y\in\boc}\t(t_{y\i}t_zt_yt_{u\i})$. Hence that inclusion is an 
equality. The proposition is proved.

\head 9. Equivalence of $\cc^\boc G$ with the centre of $\cc^\boc\cb^2$\endhead
\subhead 9.1\endsubhead
In this section we assume that the $\FF_q$-rational structure on $G$ in 0.1 is such that 

(a) {\it any $A\in CS(G)$ admits a mixed structure of pure weight $0$.}
\nl
(This can be achieved by replacing if necessary $q$ by a power of $q$.)

The bifunctor $\cc^\boc_0G\T\cc^\boc_0G@>>>\cc^\boc_0G$, $K,K'\m K\unst K'$ in 4.6 defines a bifunctor 
$\cc^\boc G\T\cc^\boc G@>>>\cc^\boc G$ denoted again by $K,K'\m K\unst K'$ as follows. Let $K\in\cc^\boc G$,
$K'\in\cc^\boc G$; we choose mixed structures of pure weight $0$ on $K,K'$ (this is possible by (a)), we 
define $K\unst K'\in\cc^\boc_0G$ as in 4.6 in terms of these mixed structures and we then disregard the 
mixed structure on $K\unst K'$. The resulting object of $\cc^\boc G$ is denoted again by $K\unst K'$; it is 
independent of the choices made.

In the same way, the bifunctor $\cc^\boc_0\cb^2\T\cc^\boc_0\cb^2@>>>\cc^\boc_0\cb^2$, $L,L'\m L\unb L'$ gives
rise to a bifunctor $\cc^\boc\cb^2\T\cc^\boc\cb^2@>>>\cc^\boc\cb^2$ denoted again by $L,L'\m L\unb L'$; the 
functor $\unc:\cc^\boc_0\cb^2@>>>\cc^\boc_0G$ gives rise to a functor $\cc^\boc\cb^2@>>>\cc^\boc G$ denoted 
again by $\unc$ (it is again called {\it truncated induction}); the functor 
$\unz:\cc^\boc_0G@>>>\cc^\boc_0\cb^2$ gives rise to a functor $\cc^\boc G@>>>\cc^\boc\cb^2$ denoted again by 
$\unz$ (it is again called {\it truncated restriction}).

The operation $K\unst K'$ is again called {\it truncated convolution}. It has a canonical associativity 
isomorphism (deduced from that in 4.7) which again satisfies the pentagon property. Thus $\cc^\boc G$
becomes a monoidal category; it has a braiding coming from 4.6(a). 

The operation $L\unb L'$ makes $\cc^\boc\cb^2$ into a monoidal abelian category (see also \cite{\CEL}). 

\subhead 9.2\endsubhead
We set 
$$\bold1=\op_{d\in\DD_\boc}\BB_d\ot\LL_d.$$
Here $\BB_d$ is as in 8.4(d). 

Let $u,z\in\boc$. From 7.7(a) we have canonically for any $d\in\DD_\boc$:
$$\Hom_{\cc^\boc\cb^2}(\LL_d,\LL_u\unb\LL_{z\i})=\Hom_{\cc^\boc\cb^2}(\LL_z,\LL_d\unb\LL_u).$$
Hence
$$\Hom_{\cc^\boc\cb^2}(\bold1',\LL_u\unb\LL_{z\i})=\Hom_{\cc^\boc\cb^2}(\LL_z,\bold1\unb\LL_u).\tag a$$
From 8.7(a) we have
$$\Hom_{\cc^\boc\cb^2}(\bold1',\LL_u\unb\LL_{z\i})=\Hom_{\cc^\boc\cb^2}(\LL_z,\LL_u).\tag b$$
From (a),(b) we deduce
$$\Hom_{\cc^\boc\cb^2}(\LL_z,\bold1\unb\LL_u)=\Hom_{\cc^\boc\cb^2}(\LL_z,\LL_u).$$
Since this holds for any $z\in\boc$, we have canonically 
$\bold1\unb\LL_u=\LL_u.$ Since this holds for any $u\in\boc$, we have 
canonically $\bold1\unb L=L$ for any $L\in\cc^\boc\cb^2$. Applying ${}^\da$, we 
deduce that we have canonically $L\unb\bold1=L$
for any $L\in\cc^\boc\cb^2$. We see that 

{\it $\bold1$ is a unit object of the monoidal category $\cc^\boc\cb^2$.}

\subhead 9.3\endsubhead
For $L\in\cc^\boc\cb^2$ let $L^*=\fD(L^\da)$. Note that $L^{**}=L$. According to \cite{\BOF}, the monoidal 
category $\cc^\boc\cb^2$ is rigid and the dual of an object $L$ is $L^*$. (I thank V. Ostrik for pointing
out the reference \cite{\BOF}.) The proof of rigidity given in \cite{\BOF} relies on the use of the 
geometric Satake isomorphism. Below we will give another approach to proving the rigidity of 
$\cc^\boc\cb^2$, which is more self contained. 

\mpb

For each $d\in\BB_d$ we choose an identification $\BB_d=\bbq$, so that $\bold1=\bold1'=\fD(\bold1)$.

As a special case of 8.7(a), for any $L\in\cc^\boc\cb^2$ we have canonically
$$\Hom_{\cc^\boc\cb^2}(\bold1,L\unb\fD(L^\da))=\Hom_{\cc^\boc\cb^2}(L,L).\tag a$$
Let $\x_L\in\Hom_{\cc^\boc\cb^2}(\bold1,L\unb\fD(L^\da))$ be the element corresponding under (a) to the 
identity homomorphism in $\Hom_{\cc^\boc\cb^2}(L,L)$. Using 3.3(a) we have 
$$\Hom_{\cc^\boc\cb^2}(\bold1,\fD(L)\unb L^\da)
=\Hom_{\cc^\boc\cb^2}(\fD(\fD(L)\unb L^\da),\fD(\bold1))=
\Hom_{\cc^\boc\cb^2}(L\unb\fD(L^\da),\bold1).$$
Under these identifications, the element $\x_{\fD(L)}\in\Hom_{\cc^\boc\cb^2}(\bold1,\fD(L)\unb L^\da)$ 
corresponds to an element $\x'_L\in\Hom_{\cc^\boc\cb^2}(L\unb\fD(L^\da),\bold1)$.
The elements $\x_L,\x'_L$ define the rigid structure on $\cc^\boc\cb^2$.

\subhead 9.4\endsubhead
Let $\cz^\boc$ be the centre of the monoidal abelian category $\cc^\boc\cb^2$. (The notion of centre of a 
monoidal abelian category was introduced by Joyal and Street \cite{\JS}, Majid \cite{\MA} and Drinfeld,
unpublished.) 

If $K\in\cc^\boc G$ then the isomorphisms 3.4(a) provide a central structure on \lb
$\unz(K)\in\cc^\boc\cb^2$ so
that $\unz(K)$ can be naturally viewed as an object of $\cz^\boc$ denoted by $\ov{\unz(K)}$. (Note that 3.4 
is stated in the mixed category but, as above, it implies the corresponding result in the unmixed 
category.) Then $K\m\ov{\unz(K)}$ is a functor $\cc^\boc G@>>>\cz^\boc$. The following result will be proved
in 9.7.

\proclaim{Theorem 9.5}The functor $\cc^\boc G@>>>\cz^\boc, K\m\ov{\unz(K)}$ is an equivalence of categories.
\endproclaim
Note that the existence of an equivalence of categories $\cc^\boc G@>>>\cz^\boc$ was conjectured by 
Bezrukavnikov, Finkelberg and Ostrik \cite{\BFO}, who constructed such an equivalence in characteristic 
zero. 

\subhead 9.6\endsubhead
By a general result on semisimple rigid monoidal categories in \cite{\ENO, Proposition 5.4}, for any
$L\in\cc^\boc\cb^2$ one can define directly a central structure on the object 
$I(L):=\op_{y\in\boc}\LL_y\unb L\unb\LL_{y\i}$ of $\cc^\boc\cb^2$ such that, denoting by $\ov{I(L)}$ the
corresponding object of $\cz^\boc$, we have canonically
$$\Hom_{\cc^\boc\cb^2}(L,L')=\Hom_{\cz^\boc}(\ov{I(L)},L')\tag a$$
for any $L'\in\cz^\boc$. (We use that for $y\in\boc$, the dual of the simple object $\LL_y$ of 
$\cc^\boc\cb^2$ is $\LL_{y\i}$.) The central structure on $I(L)$ can be described as follows: for any 
${}^1L\in\cc^\boc\cb^2$ we have canonically
$$\align&{}^1L\unb I(L)=\op_{y\in\boc}{}^1L\unb\LL_y\unb L\unb\LL_{y\i}=
\op_{y,z\in\boc}\Hom_{\cc^\boc\cb^2}(\LL_z,{}^1L\unb\LL_y)\ot\LL_z\unb L\unb\LL_{y\i}\\&
=\op_{y,z\in\boc}\Hom_{\cc^\boc\cb^2}(\LL_{y\i},\LL_{z\i}\unb{}^1L)\ot\LL_z\unb L\unb\LL_{y\i}\\&
=\op_{z\in\boc}\LL_z\unb L\unb\LL_{z\i}\unb{}^1L=I(L)\unb{}^1L.\endalign$$

\subhead 9.7\endsubhead
For $x\in\boc$ we have canonically $\unz\unc\LL_x=I(\LL_x)$ as objects of $\cc^\boc\cb^2$, see Theorem 6.8.
From the last commutative diagram in 7.9 we see that this identification is compatible with the central 
structures (see 9.4, 9.6), so that 
$$\ov{\unz\unc \LL_x}=\ov{I(\LL_x)}.\tag a$$
Using this and 9.6(a) with $L'=\ov{\unz\unc\tL}$, $\tL\in\cc^\boc\cb^2$, we see that
$$\Hom_{\cc^\boc\cb^2}(\LL_x,\unz\unc\tL)=\Hom_{\cz^\boc}(\ov{\unz\unc\LL_x},\ov{\unz\unc\tL}).$$
Combining this with 8.9 we obtain for    $\tL=\LL_{x'}$ (with $x'\in\boc$):
$$\AA_{x,x'}=\AA'_{x,x'}\tag b$$
where
$$\AA_{x,x'}=\Hom_{\cc^\boc G}(\unc(\LL_x),\unc(\LL_{x'})),
\AA'_{x,x'}=\Hom_{\cz^\boc}(\ov{\unz\unc\LL_x},\ov{\unz\unc\LL_{x'}}).$$
Note that the identification (b) is induced by the functor $K\m\ov{\unz(K)}$.
Let $\AA=\op_{x,x'\in\boc}\AA_{x,x'},\AA'=\op_{x,x'\in\boc}\AA'_{x,x'}$. Then from (b) we have $\AA=\AA'$. 
Note that this identification is compatible with the obvious algebra structures of $\AA,\AA'$.

For any $A\in CS_\boc$ we denote by $\AA_A$ the set of all $f\in\AA$ such that for any $x,x'$, the 
$(x,x')$-component of $f$ maps the $A$-isotypic component of $\unc(\LL_x)$ to the $A$-isotypic component of 
$\unc(\LL_{x'})$ and any other isotypic component of $\unc(\LL_x)$ to $0$. Then $\AA=\op_{A\in CS_\boc}\AA_A$
is the decomposition of $\AA$ into simple algebras (each $\AA_A$ is $\ne0$ since, by 1.7(b) and 1.10(a), any
$A$ is a summand of some $\unc(\LL_x)$).

From \cite{\MU}, \cite{\ENO}, we see that $\cz^\boc$ is a semisimple abelian category with finitely many 
simple objects up to isomorphism. Let $\fS$ be a set of representatives for the isomorphism classes of 
simple objects of $\cz^\boc$. For any $\s\in\fS$ we denote by $\AA'_\s$ the set of all $f'\in\AA'$ such that
for any $x,x'$, the $(x,x')$-component of $f'$ maps the $\s$-isotypic component of $\ov{\unz\unc(\LL_x)}$ to
the $\s$-isotypic component of $\ov{\unz\unc(\LL_{x'})}$ and any other isotypic component of 
$\ov{\unz\unc(\LL_x)}$ to $0$. Then $\AA'=\op_{\s\in\fS}\AA'_\s$ is the decomposition of $\AA'$ into a sum
of simple 
algebras (each $\AA'_\s$ is $\ne0$ since any $\s$ is a summand of some $\ov{\unz\unc(\LL_z)}$; indeed, if 
$\LL_z$ is a summand of $\s$ viewed as an object of $\cc^\boc\cb^2$ then by 9.6(a), $\s$ is a summand of
$\ov{I(\LL_z)}$ hence of $\ov{\unz\unc(\LL_z)}$).

Since $\AA=\AA'$, from the uniqueness of decomposition of a semisimple algebra as a direct sum of simple
algebras, we see that there is a unique bijection $CS_\boc\lra\fS$, $A\lra\s_A$ such that
$\AA_A=\AA'_{\s_A}$ for any $A\in CS_\boc$. From the definitions we now see that for any $A\in CS_\boc$ we 
have $\ov{\unz A}\cong\s_A$. Therefore Theorem 9.5 holds.

\proclaim{Theorem 9.8} Let $L\in\cc^\boc\cb^2$, $K\in\cc^\boc G$. We have canonically
$$\Hom_{\cc^\boc\cb^2}(L,\unz(K))=\Hom_{\cc^\boc G}(\unc(L),K).\tag a$$
Moreover, in $\cc^\boc\cb^2$ we have $\unz(K)\cong\op_{z\in\boc^0}\LL_z^{\op m_z}$ where $\boc^0$ is as in 
1.7 and $m_z\in\NN$.
\endproclaim
From 9.5, 9.7, we see that 
$$\Hom_{\cc^\boc G}(\unc(L),K)=\Hom_{\cz^\boc}(\ov{\unz\unc(L)},\ov{\unz K})=
\Hom_{\cz^\boc}(\ov{I(L)},\ov{\unz K}).$$
Using 9.6(a) we see that $\Hom_{\cz^\boc}(\ov{I(L)},\ov{\unz K})=\Hom_{\cc^\boc\cb^2}(L,\unz(K))$ and (a) 
follows. To prove the second assertion of the theorem it is enough to show that for any $z\in\boc-\boc^0$ we
have $\Hom_{\cc^\boc\cb^2}(\LL_z,\unz(K))=0$; by (a), it is enough to show that $\unc(\LL_z)=0$ and this
follows from 1.7(c).

\subhead 9.9\endsubhead
We show that for $K\in\cc^\boc G$ we have canonically
$$\fD(\unz(\fD(K)))=\unz(K).\tag a$$
It is enough to show that for any $L\in\cc^\boc\cb^2$ we have canonically
$$\Hom_{\cc^\boc\cb^2}(L,\fD(\unz(\fD(K))))=\Hom_{\cc^\boc\cb^2}(L,\unz(K)).$$
Here the left side equals 
$$\align&\Hom_{\cc^\boc\cb^2}(\unz(\fD(K)),\fD(L))=\Hom_{\cc^\boc G}(\fD(K),\unc(\fD(L)))\\&=
\Hom_{\cc^\boc G}(\fD(K),\fD(\unc(L)))\endalign$$
(we have used 9.8(a) and 1.13(a)) and the right hand side equals 
$$\Hom_{\cc^\boc G}(\unc(L),K)=\Hom_{\cc^\boc G}(\fD(K),\fD(\unc(L))).$$
(We have again used 9.8(a)). This proves (a).

\subhead 9.10\endsubhead
The monoidal structure on $\cc^\boc\cb^2$ induces a monoidal structure on $\cz^\boc$. Using 5.2 and the 
definitions we see the equivalence of categories in 9.5 is compatible with the monoidal structures. Since 
$\cz^\boc$ has a unit object, it follows that the monoidal category $\cc^\boc G$ also has a unit object, say
$A$. We show:
$$A\cong A_{E_\boc}.\tag a$$
From 8.6(d),(e) we see that for $x\in\boc$, $(A_{E_\boc}:\unc(\LL_x))$ is $1$ if $x\in\DD_\boc$ and is $0$
if $x\n\DD_\boc$. Using 9.8(a) we deduce that for $x\in\boc$, $\dim\Hom_{\cd\cb^2}(\LL_x,\unz(A_{E_\boc}))$ 
is $1$ if $x\in\DD_\boc$ and is $0$ if $x\n\DD_\boc$. Thus $\unz(A_{E_\boc})$ is isomorphic in 
$\cc^\boc\cb^2$ to the unit object $\bold1$ of the monoidal category $\cc^\boc\cb^2$. Then
$\unz(A_{E_\boc})$ viewed as an object of $\cz^\boc$ is also the unit object of $\cz^\boc$ hence is
isomorphic in $\cz^\boc$ to $\unz(A)$. Using Theorem 9.5 we deduce that (a) holds.

\subhead 9.11\endsubhead
Let $z,u\in\boc$. We have canonically
$$\unc(\LL_z)\unst\unc(\LL_u)=\op_{y\in\boc}\unc(\LL_u\unb\LL_y\unb\LL_z\unb\LL_{y\i}).\tag a$$
Indeed, by 8.1(a), it is enough to prove that we have canonically
$$\unc(\LL_u\unb\unz\unc(\LL_z))=\op_{y\in\boc}\unc(\LL_u\unb\LL_y\unb\LL_z\unb\LL_{y\i})$$
and this follows from 6.8(a). We see that
$$\unc(\LL_z)\unst\unc(\LL_u)\cong\op_{r\in\boc^0}\unc(\LL_r)^{\op\ps(r)}$$
in $\cc^\boc G$ where $\ps(r)\in\NN$ are given by the following equation in $\JJ^\boc$:
$$\sum_{y\in\boc}t_{u}t_yt_zt_{y\i}=\sum_{r\in\boc}\ps(r)t_z.$$

\subhead 9.12\endsubhead
Let $\JJ_0^\boc$ be the subgroup of $\JJ^\boc$ spanned by $\{t_z;z\in\boc^0\}$. For $\x,\x'\in\JJ_0^\boc$ we
set 
$$\x\circ\x'=\sum_{y\in\boc}\x t_y\x't_{y\i}\in\JJ^\boc.$$
We show that $\x\circ\x'\in\JJ_0^\boc$. We can assume that $\x=t_w,\x=t_{w'}$ with $w,w'\in\boc^0$. If 
$t_z (z\in\boc)$ appears with nonzero coefficient in $\x\circ\x'$ then $t_{z\i}t_wt_yt_{w'}t_{y\i}\ne0$ for 
some $y\in\boc$ and $t_wt_{y'}t_{w'}t_{y'{}\i}t_{z\i}\ne0$ for some $y'\in\boc$. Using \cite{\HEC, P8} we 
deduce: $z\i\si_Lw\i$, $w\si_Ly'{}\i$, $y'{}\i\si_Lz$. Since $w\si_Lw\i$, it follows that $z\si_L z\i$, as 
claimed.

For $\x,\x',\x''$ in $\JJ_0^\boc$ we show that $(\x\circ\x')\circ\x''=\x\circ(\x'\circ\x'')$.
We can assume that $\x=t_w,\x'=t_{w'},\x''=t_{w''}$ where $w,w',w''$ are in $\boc^0$. We must show:
$$\sum_{y,u\in\boc}t_wt_yt_{w'}t_{y\i}t_ut_{w''}t_{u\i}=
\sum_{y,s\in\boc}t_wt_yt_{w'}t_st_{w''}t_{s\i}t_{y\i},$$
or equivalently
$$\sum_{y,u,s\in\boc}h^*_{y\i,u,s}t_wt_yt_{w'}t_st_{w''}t_{u\i}=
\sum_{y,s,u\in\boc}h^*_{s\i,y\i,u\i}t_wt_yt_{w'}t_st_{w''}t_{u\i}.$$
It remains to use the identity $h^*_{y\i,u,s}=h^*_{s\i,y\i,u\i}$ for $y,u,s\in\boc$ (see \cite{\HEC, P7}).

We see that $(\JJ_0^\boc,\circ)$ is an associative ring (without $1$ in general). Let $\cg$ be the 
Grothendieck group of the category $\cc^\boc G$; this is an associative and commutative ring under truncated 
convolution (see 9.1) and 9.11 shows that $t_w\m\unc(\LL_z)$ is a ring homomorphism $\JJ_0^\boc@>>>\cg$.

\head 10. Remarks on the noncrystallographic case\endhead
\subhead 10.1\endsubhead
In this subsection we consider a not necessarily crystallographic Coxeter group $W'$ with a fixed
two-sided cell $\boc'$. The following discussion assumes the truth of Soergel's conjecture for $W'$, 
recently proved by Elias and Williamson \cite{\EW}. Let $w\m|w|$ be the length function of $W'$.
For any $w\in W'$ we define $\aa(w)\in\NN$ as in \cite{\HEC, 13.6}. (The assumption in {\it loc.cit.} that
$W'$ with $w\m|w|$ is bounded in the sense of \cite{\HEC, 13.2} is not necessary for the definition of 
$\aa(w)$; to show that $\aa(w)$ is well defined we use instead the inequality $\aa(w)\le|w|$ which is proved
by the argument in \cite{\HEC, 15.2}, applicable in view of the positivity results of \cite{\EW}.) Now the 
properties of $\aa(w)$ stated in \cite{\HEC, 14.2} hold by the arguments in \cite{\HEC, \S15}, using again 
the positivity results in \cite{\EW}. The ring $J$ and its subring $J^{\boc'}$ is defined as in 
\cite{\HEC, 18.3} in terms of the $\aa$-function (the assumptions in \cite{\HEC, 18.1} are not needed in our
case.) Note that the ring $J^{\boc'}$ does not have in general a unit element (unless $\boc'$ contains only 
finitely many left cells.)

We show that the definition of the monoidal category in 3.2 can be adapted to the more general case of 
$W',\boc'$ by using Soergel bimodules \cite{\SO} instead of perverse sheaves.

Let $R$ be the algebra of regular real valued functions on a fixed (real) reflection representation of $W'$.
Then 
for each $x\in W'$, the indecomposable Soergel graded $R$-bimodule $B_x$ is defined as in \cite{\SO, 6.16}. 
Let $\tC$ (resp. $C$) be the category of graded $R$-bimodules wich are isomorphic to finite direct sums of 
graded $R$-bimodules of the form $B_x$ with shift (resp. without shift). As shown by Soergel, $\tC$ is a 
monoidal category under the usual tensor product $L,L'\m L\cir L'$.
If $L\in\tC$ and $j\in\ZZ$ we write $L^j\in C$ for what in \cite{\EW, 6.2} is denoted by $\ch^j(L)$. (The 
fact that $L^j$ is well defined follows from the results of \cite{\SO} and \cite{\EW}.)
Let $C_{\boc'}$ be the category of graded $R$-bimodules wich are isomorphic to finite direct 
sums of graded $R$-bimodules of the form $B_x$ ($x\in\boc'$) without shift. For any $L\in C$ there is a 
unique direct sum decomposition $L=\un{L}\op L'$ where $\un{L}\in C_{\boc'}$ and $L'$ is a direct sum of 
graded $R$-bimodules of the form $B_x$ ($x\n\boc'$). (The uniqueness of this direct sum decomposition 
follows from the results of \cite{\SO} and \cite{\EW}.) Let $a'$ be the value on $\boc'$ of the 
$\aa$-function $W'@>>>\NN$. By arguments parallel to those in \cite{\CEL} and making use of the results of 
\cite{\EW} and the properties of the $\aa$-function we see that for $L,L'\in C_{\boc'}$ we have 
$\un{(L\cir L')^j}=0$ if $j>a'$ and 
$L,L'\m L\unb L':=\un{(L\cir L')^{a'}}$ defines a monoidal structure on $C_{\boc'}$ (without a unit object
unless $\boc'$ is a union of finitely many left cells) such that the induced 
ring structure on the Grothendieck group 
of $C_{\boc'}$ is isomorphic to the ring $J^{\boc'}$. (For three objects $L,L',L''$ in $C_{\boc'}$ we 
have $(L\unb L')\unb L''=L\unb(L'\unb L'')=\un{(L\cir L'\cir L'')^{2a'}}$.)
(Note that in the finite crystallographic case, the objects of $C_{\boc'}$ should be thought of as
perverse sheaves on $\cb$ rather than on $\cb^2$ as in 3.2; this accounts for our usage of $a'$ instead
of the $a-\nu$ in 3.2.) Let $Z^{\boc'}$ be the centre of the monoidal category
$(C_{\boc'},\unb)$. (We use the definition of the center given in \cite{\MU, 3.4} in which the property
(iv) of a half braiding given in \cite{\MU, 3.1} is omitted so that the definition makes sense without
requiring the existence of a unit object; note that if a unit object exists then property (iv) in
\cite{\MU, 3.1} is a consequence of properties (i)-(iii) in \cite{\MU, 3.1}, see \cite{\MU, 3.2}.) By
\cite{\MU, 3.5, 3.6}, $Z^{\boc'}$ is an $\RR$-linear category.
Let $\fS_{\boc'}$ be the set of isomorphism classes of objects of $Z^{\boc'}$
which are indecomposable with respect to direct sum.
The objects of $\fS_{\boc'}$ can be called the {\it character sheaves} of $W',\boc'$;
this is justified by Theorem 9.5.

Now assume further that $W'$ is finite, of type $H_3$ or $H_4$ or a dihedral group. In this case $\boc'$ is 
uniquely determined by the number $a'$. Recall that in \cite{\COX} the ``unipotent characters'' associated 
to $W'$ were ``described''. The unipotent characters whose degree polynomial is divisible by $q^{a'}$ but 
not by $q^{a'+1}$ can be viewed as unipotent characters associated to $\boc'$; they form a set 
$\cu_{\boc'}$. We expect that $\cu_{\boc'}$ and $\fS_{\boc'}$ are in a natural bijection. This predicts for 
example that, if $\boc'$ in type $H_4$ has $a'=6$, then $\fS_{\boc'}$ has exactly $74$ elements; if $\boc'$ 
for the dihedral group of order $4k+2$ (resp. $4k+4$) has $a'=1$ then $\fS_{\boc'}$ has exactly $k^2$ (resp.
$k^2+k+2$) elements. We expect that the monoidal category $C_{\boc'}$ is rigid, so that (by a result of 
\cite{\MU}, \cite{\ENO}), $Z^{\boc'}$ is a semisimple abelian category and $\fS_{\boc'}$ is the same as the 
set of isomorphism classes of simple objects of $Z^{\boc'}$. We also expect that $Z^{\boc'}$ is a modular 
tensor category whose $S$-matrix is the matrix described in \cite{\EXO}, \cite{\MG}, which transforms the
 fake degrees polynomials of $W'$ corresponding to $\boc'$ to
the unipotent character degrees corresponding to $\boc'$.

\enddocument